\documentclass[12pt]{article}
\usepackage{amsfonts}
\usepackage{amssymb}
\usepackage{amsmath}
\usepackage{graphicx}
\usepackage{subfigure}
\usepackage{theorem}

\def\date{9 March 2012}  
\newtheorem{proposition}{proposition}[section]
\newtheorem{lemma}[proposition]{Lemma}
\newtheorem{corollary}[proposition]{Corollary}
\newtheorem{theorem}[proposition]{Theorem}
\newtheorem{claim}[proposition]{Claim}
\newtheorem{conjecture}[proposition]{Conjecture}
\newcommand{\qed}{$\square$\bigskip}

\newcommand{\proof}{{\noindent\bf Proof. }}
{\theorembodyfont{\rmfamily}

\newtheorem{hypothesis}[proposition]{Hypothesis}
}
\def\mytextindent#1{\indent\indent\llap{\rm#1\enspace}\ignorespaces}
\def\myitem{\par\hangindent0pt\mytextindent}
\def\claim#1#2{\bigskip\noindent\rlap{\rm(#1)}\ignorespaces
   \hangindent=30pt\hskip30pt{\sl#2}\bigskip}
\def\ifc{internally $4$-connected}

\def\myclaim#1#2\par{{\medbreak\noindent\rlap{\rm(#1)}\ignorespaces
 \rightskip20pt
 \hangindent=20pt\hskip20pt{\ignorespaces\sl#2}\smallskip}}
\def\junk#1{}

\def\restriction{|}

\font\smallrm=cmr8
\font\smallrm=cmr8
\def\minimal{minimal}
\def\uom{\Omega}
\def\lam{\lambda}
\def\calw{{\cal W}}
\def\udel{\Delta}
\def\ugam{\Gamma}
\begin{document}
\newif\ifproofmode
\baselineskip=12pt
\phantom{a}\vskip .25in 
\centerline{{\bf  $K_6$ MINORS IN LARGE 6-CONNECTED GRAPHS}}
\vskip.4in
\centerline{\bf Ken-ichi Kawarabayashi}
\centerline{National Institute of Informatics}
\centerline{2-1-2 Hitotsubashi, Chiyoda-ku, Tokyo 101-8430, Japan}
\medskip
\centerline{{\bf Serguei Norine}\footnote{\smallrm
Partially supported by NSF under Grants No.~DMS-0200595 and DMS-0701033.}}
\centerline{Department of Mathematics}
\centerline{Princeton University}
\centerline{Princeton, NJ 08544, USA}
\medskip
\centerline{{\bf Robin Thomas}%
\footnote{\smallrm Partially supported by NSF under
Grants No.~DMS-0200595, DMS-0354742, and DMS-0701077. }}
\centerline{School of Mathematics} 
\centerline{Georgia Institute of Technology} 
\centerline{Atlanta, Georgia  30332-0160, USA}
\medskip
\centerline{and}
\medskip
\centerline{\bf Paul Wollan}
\centerline{Mathematisches Seminar der Universit\"at Hamburg}
\centerline{Bundesstrasse 55}
\centerline{D-20146 Hamburg, Germany}

\vskip 1in \centerline{\bf ABSTRACT}
\bigskip
\parshape=1.5truein 5.5truein
J\o rgensen conjectured that every 6-connected graph $G$ with no
$K_6$ minor has a vertex whose deletion makes the graph planar.
We prove the conjecture for all sufficiently large graphs.

\vfill \baselineskip 11pt \noindent 8 April 2005, revised \date.
\vfil\eject

\ifproofmode
  \def\sectionbreak{\vfil\eject}
   \newcount\remarkno
   \def\REMARK#1{{%
      \footnote{\baselineskip=11pt #1
      \vskip-\baselineskip}\global\advance\remarkno by1}}
    \def\mylabel#1{{\label{#1}\marginpar{#1}}}
   \def\showlabel#1{{\marginpar{#1}}}
   \def\showfiglabel#1{\marginpar{#1}}
\else
  \def\sectionbreak{}
  \def\REMARK#1{}
  \def\mylabel#1{\label{#1}}
  \def\showlabel#1{}
  \def\showfiglabel#1{}
\fi

\baselineskip 18pt

\section{Introduction}

{\em Graphs} in this paper are allowed to have loops and multiple edges.
A graph is a {\em minor} of another if the first can be obtained from 
a subgraph of the second by contracting edges. An {\em $H$ minor}
is a minor isomorphic to $H$.
A graph $G$ is {\em apex} if it has a vertex $v$ such that
$G\backslash v$ is planar.
(We use $\backslash$ for deletion.)
J\o rgensen~\cite{Jor} made the following beautiful conjecture.

\begin{conjecture}
\label{con:jorgensen}
\showlabel{con:jorgensen}
Every $6$-connected graph with no $K_6$ minor is apex.
\end{conjecture}

\noindent
This is related to Hadwiger's conjecture~\cite{Had}, the following.

\begin{conjecture}
\label{con:hadwiger}
\showlabel{con:hadwiger}
For every integer $t\ge1$, if a loopless graph has no $K_t$ minor, then it
is $(t-1)$-colorable.
\end{conjecture}

Hadwiger's conjecture is known for $t\le6$. 
For $t=6$ it has been proven in~\cite{RobSeyThoHad} by showing
that a minimal counterexample to Hadwiger's conjecture for $t=6$
is apex.
The proof uses an earlier result of Mader~\cite{MadHomkrit}
that every minimal counterexample to Conjecture~\ref{con:hadwiger}
is $6$-connected.
Thus Conjecture~\ref{con:jorgensen}, if true, would give more
structural information.
Furthermore, the structure of all graphs with no $K_6$ minor is not known,
and appears complicated and difficult.
On the other hand, Conjecture~\ref{con:jorgensen} provides a nice
and clean statement for $6$-connected graphs.
Unfortunately, it, too, appears to be a difficult problem.
In this paper we prove Conjecture~\ref{con:jorgensen} for all
sufficiently large graphs, as follows.

\begin{theorem}
\label{thm:largejorgensen}
\showlabel{thm:largejorgensen}
There exists an absolute constant $N$ such that every $6$-connected
graph on at least $N$ vertices with no $K_6$ minor is apex.
\end{theorem}

The second and third author recently announced
a generalization~\cite{NorTho} 
of Theorem~\ref{thm:largejorgensen}, where $6$ is
replaced by an arbitrary integer $t$.
The result states that for every integer $t$ there exists 
an integer $N_t$ such that every $t$-connected graph on at least
$N_t$ vertices with no $K_t$ minor has a set of at most $t-5$ vertices
whose deletion makes the graph planar.
The proof follows a different strategy, but makes use of several
ideas developed in this paper and its companion~\cite{KawNorThoWolbdtw}.

We use a number of results from the Graph Minor series of Robertson
and Seymour, and also three results of our own that are proved
in~\cite{KawNorThoWolbdtw}.
The first of those is a version of Theorem~\ref{thm:largejorgensen}
for graphs of bounded tree-width, the following. 
(We will not define tree-width here, because it is sufficiently well-known,
and because we do not need the concept {\it per se}, 
only several theorems that use it.)

\begin{theorem}
\label{thm:bdedtwjorgensen}
\showlabel{thm:bdedtwjorgensen}
For every integer $w$ there exists an integer $N$
such that every $6$-connected
graph of tree-width at most $w$ on at least $N$ vertices and with no 
$K_6$ minor is apex.
\end{theorem}

Theorem~\ref{thm:bdedtwjorgensen} reduces the proof of 
Theorem~\ref{thm:largejorgensen} to graphs of large tree-width.
By a result of Robertson and Seymour~\cite{RobSeyGM5} those graphs
have a large grid minor.
However, for our purposes it is more convenient to work with walls
instead.
Let $h\ge2$ be even.
An {\em elementary wall of height $h$} has vertex-set
$$\{(x,y): 0\le x\le 2h+1, 0\le y\le h\}-\{(0,0),(2h+1,h)\}$$
and an edge between any vertices $(x,y)$ and $(x',y')$ if either
\myitem{$\bullet$} $|x-x'|=1$ and $y=y'$, or 
\myitem{$\bullet$} $x=x'$, $|y-y'|=1$ and $x$ and $\max\{y,y'\}$
have the same parity.

\begin{figure}
\begin{center}
\leavevmode
\includegraphics[scale = 1]{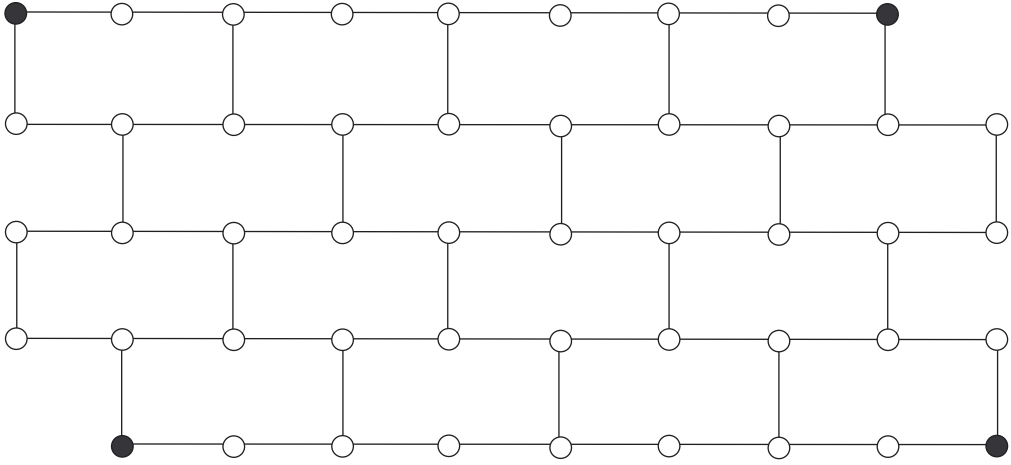}
\end{center}
\caption{
An elementary wall of height $4$.
}
\label{fig:wall}
\end{figure}

\noindent
Figure~\ref{fig:wall} shows an elementary wall of height $4$.
A {\em wall of height $h$} is a subdivision of an elementary wall
of height $h$. The result of~\cite{RobSeyGM5} 
(see also~\cite{DieGorJenTho,ReedBCC,RobSeyThoQuickPlanar}) 
can be restated as follows.
\nocite{DieGorJenTho}
\nocite{ReedBCC}
\nocite{RobSeyThoQuickPlanar}

\showfiglabel{fig:wall}

\begin{theorem}
\label{thm:grid}
\showlabel{thm:grid}
For every even integer $h\ge2$ there exists an integer $w$ such that
every graph of tree-width at least $w$ has a subgraph isomorphic to
a wall of height $h$.
\end{theorem}

The {\em perimeter} of a wall is the cycle that bounds the infinite face
when the wall is drawn as in Figure~\ref{fig:wall}.
Now let $C$ be the perimeter of a wall $H$ in a graph $G$.
The {\em compass of $H$ in $G$} is the restriction of $G$ to $X$, where
$X$ is the union of $V(C)$ and the vertex-set of the unique component
of $G\backslash V(C)$ that contains a vertex of $H$.
Thus $H$ is a subgraph of its compass, and the compass is connected.
A wall $H$ with perimeter $C$ in a graph $G$ is {\em planar} if its compass
can be drawn in the plane with $C$ bounding the infinite face.
In Section~\ref{sec:planarwall} we prove the following.

\begin{theorem}
\label{thm:plwall}
\showlabel{thm:plwall}
For every even integer $t\ge2$ there exists an even integer $h\ge2$
such that if a $5$-connected graph $G$ with no $K_6$ minor
has a wall of height at least $h$, then either it is apex, or
has a planar wall of height $t$.
\end{theorem}

Actually, in the proof of Theorem~\ref{thm:plwall} we need 
Lemma~\ref{lem:pinwheel} that is proved in~\cite{KawNorThoWolbdtw}. The lemma says
that if a $5$-connected graph with no $K_6$ minor has a 
subgraph isomorphic to subdivision
of a pinwheel with sufficiently many vanes 
 (see Figure~\ref{fig:pinwheel}), then it is apex.

By Theorem~\ref{thm:plwall} we may assume that our graph $G$ has
an arbitrarily large planar wall $H$. Let $C$ be the perimeter of
$H$, and let $K$ be the compass of $H$. Then $C$ separates $G$ into
$K$ and another graph, say $J$, such that $K\cup J=G$, 
$V(K)\cap V(J)=V(C)$ and $E(K)\cap E(J)=\emptyset$.
Next we study the graph $J$. Since the order of the vertices on $C$
is important, we are lead to the notion of a ``society", introduced
by Robertson and Seymour in~\cite{RobSeyGM9}.

Let $\Omega$ be a cyclic permutation of the elements of some set;
we denote this set by $V(\Omega)$.
A {\em society} is a pair $(G,\Omega)$, where $G$ is a graph, and $\Omega$
is a cyclic permutation with $V(\Omega)\subseteq V(G)$. 
Now let $J$ be as above, and let $\Omega$ be one of the cyclic permutations
of $V(C)$ determined by the order of vertices on $C$.
Then $(J,\Omega)$ is a society that is of primary interest to us.
We call it the {\em anticompass  society of $H$ in $G$}.

We say that $(G,\uom,\uom_0)$ is a {\em neighborhood} if $G$ is a graph
and $\uom,\uom_0$ are cyclic permutations, where both $V(\uom)$
and $V(\uom_0)$ are subsets of $V(G)$.
Let $\Sigma$ be a plane, with some orientation called ``clockwise."
We say that a neighborhood $(G,\uom,\uom_0)$ is {\em rural} if $G$ has
a drawing $\Gamma$ in $\Sigma$ without crossings (so $G$ is planar)
and there are closed discs $\Delta_0\subseteq \Delta\subseteq \Sigma$, such
that

\myitem{(i)} the drawing $\Gamma$ uses no point of $\Sigma$ outside $\Delta$,
and none in the interior of $\Delta_0$, and
\myitem{(ii)} for $v\in V(G)$, the point of $\Sigma$ representing $v$ in the
drawing $\Gamma$ lies in $bd(\Delta)$ (respectively, $bd(\Delta_0)$)
if and only if $v\in V(\uom)$ (respectively, $v\in V(\uom_0))$, and
the cyclic permutation of $V(\uom)$ (respectively, $V(\uom_0))$
obtained from the clockwise orientation of $bd(\Delta)$ (respectively,
$bd(\Delta_0)$) coincides (in the natural sense) with $\uom$ (respectively,
$\uom_0$).

\noindent We call $(\Sigma,\Gamma,\Delta,\Delta_0)$ a {\em presentation}
of $(G,\uom,\uom_0)$.

Let $(G_1,\Omega,\Omega_0)$ be a neighborhood, let $(G_0,\Omega_0)$
be a society with $V(G_0)\cap V(G_1)=V(\Omega_0)$, and let $G=G_0\cup G_1$.
Then $(G,\Omega)$ is a society, and we say that $(G,\Omega)$
is the {\em composition} of the society $(G_0,\Omega_0)$ with
the neighborhood $(G_1,\Omega,\Omega_0)$.
If the neighborhood $(G_1,\Omega,\Omega_0)$ is rural, then
we say that $(G_0,\Omega_0)$ is a {\em planar truncation} of $(G,\Omega)$.
We say that a society $(G,\Omega)$ is {\em $k$-cosmopolitan}, where $k\ge 0$
is an integer, if for every planar truncation $(G_0,\Omega_0)$ of
$(G,\Omega)$ at least $k$ vertices in $V(\Omega_0)$ have at least two
neighbors in $V(G_0)$.
At the end of Section~\ref{sec:planarwall} we deduce

\begin{theorem}
\label{cosmopolitan}
\showlabel{cosmopolitan}
For every integer $k\ge 1$ there exists an even integer $t\ge 2$ such that
if $G$ is a simple graph of minimum degree at least six and $H$ is a
planar wall of height $t$ in $G$, then the anticompass society of
$H$ in $G$ is $k$-cosmopolitan.
\end{theorem}

For a fixed presentation $(\Sigma,\Gamma,\Delta,\Delta_0)$
of a neighborhood $(G,\Omega,\Omega_0)$ and an integer $s\ge0$ we define an
{\em $s$-nest} for  $(\Sigma,\Gamma,\Delta,\Delta_0)$
to be a sequence $(C_1,C_2,\ldots,C_s)$
of pairwise disjoint cycles of $G$ such that
$\Delta_0\subseteq\Delta_1\subseteq\cdots\subseteq\Delta_s\subseteq\Delta$,
where $\Delta_i$ denotes the closed disk in $\Sigma$ bounded by
the image under $\Gamma$ of $C_i$.
We say that a society $(G,\Omega)$  is {\em $s$-nested} if it
is the composition of a society $(G_1,\Omega_0)$ with a rural
neighborhood $(G_2,\Omega,\Omega_0)$ that has
an $s$-nest for some presentation of $(G_2,\Omega,\Omega_0)$.

Let $\Omega$ be a cyclic permutation.
For $x\in V(\Omega)$ we denote the image of $x$ under $\Omega$
by $\Omega(x)$.
If $X\subseteq V(\Omega)$, then we denote by $\Omega|X$ the restriction
of $\Omega$ to $X$.
That is, $\Omega|X$ is the permutation $\Omega'$ defined by saying
that $V(\Omega')=X$ and $\Omega'(x)$ is the first term of the
sequence $\Omega(x),\Omega(\Omega(x)),\ldots$ which belongs to $X$.
Let $v_1,v_2,\ldots,v_k\in V(\Omega)$ be distinct. 
We say that $(v_1,v_2,\ldots,v_k)$ is {\em clockwise} in $\Omega$ (or simply
{\em clockwise} when $\Omega$ is understood from context) if
$\Omega'(v_{i-1})=v_{i}$ for all $i=1,2,\ldots,k$, where $v_0$ means $v_k$
and $\Omega'=\Omega|\{v_1,v_2,\ldots,v_k\}$.
For $u,v\in V(\Omega)$ we define $u\Omega v$ as the set of all
$x\in V(\Omega)$ such that either $x=u$ or $x=v$ or $(u,x,v)$ is clockwise
in $\Omega$.

A separation of a graph is a pair $(A,B)$ such that $A\cup B=V(G)$
and there is no edge with one end in $A-B$ and the other end in $B-A$.
The order of $(A,B)$ is $|A\cap B|$.
We say that a society $(G,\Omega)$ is {\em $k$-connected}
if there is no separation $(A,B)$ of $G$ of order at most $k-1$ 
with $V(\Omega)\subseteq A$ and $B-A\ne\emptyset$.
A {\em bump} in $(G,\Omega)$ is a path in $G$ with at least one edge,
both ends in $V(\Omega)$ and otherwise disjoint from $V(\Omega)$.

\begin{figure}
\begin{center}
\leavevmode
\includegraphics[scale = 0.75]{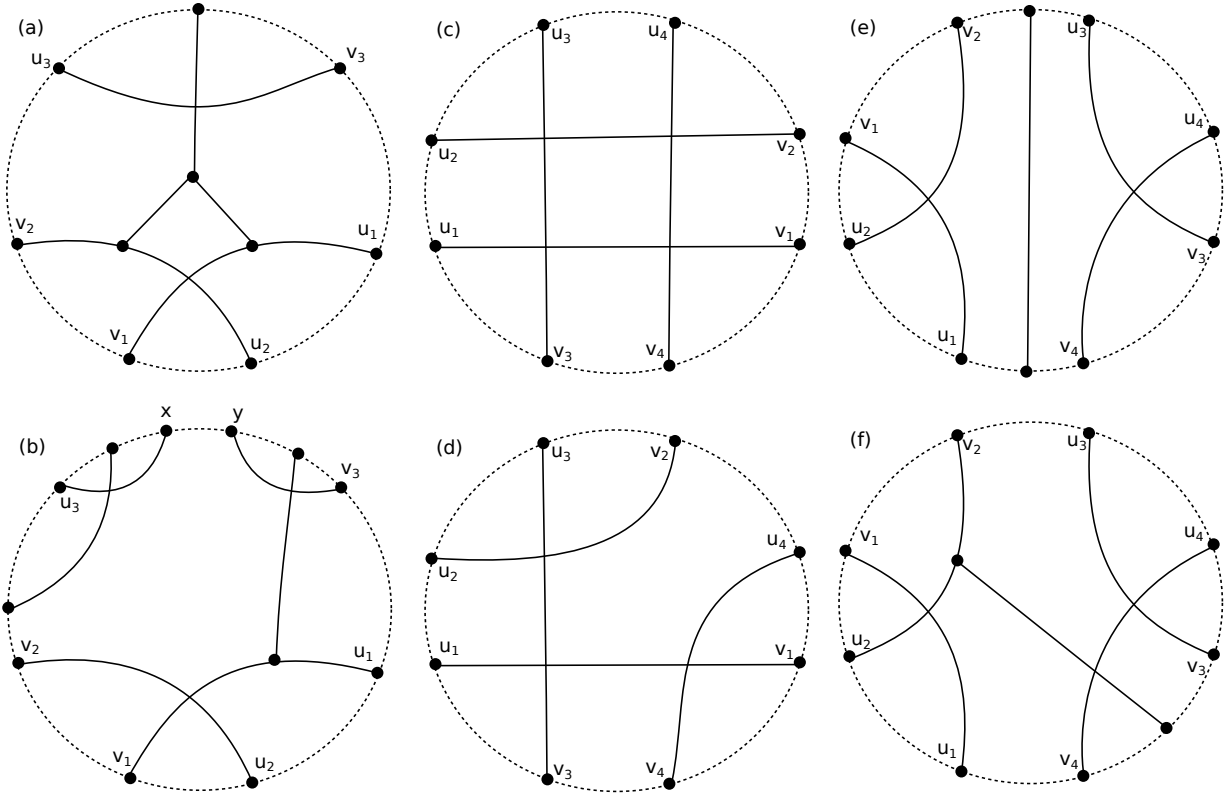}
\end{center}
\caption{
(a),(b) A turtle. (c),(d) A gridlet. (e),(f) A separated doublecross.
}
\label{fig:bumps}
\end{figure}
\showfiglabel{fig:bumps}

Let  $(G,\Omega)$ be a society and let $(u_1,u_2,v_1,v_2,u_3,v_3)$ be
clockwise in $\Omega$. For $i=1,2$ let $P_i$ be a bump in $G$
with ends $u_i$ and $v_i$,
and let $L$ be either a bump with ends $u_3$ and $v_3$,
or the union of two internally disjoint bumps, one with
ends $u_3$ and $x\in u_3\Omega v_3$ and the other with ends $v_3$ and
$y\in u_3\Omega v_3$. In the former case let $Z=\emptyset$, and in
the latter case let $Z$ be the subinterval of $u_3\Omega v_3$ with
ends $x$ and $y$, including its ends. 
Assume that $P_1,P_2,L$ are pairwise disjoint.
Let $q_1,q_2\in V(P_1)\cup V(P_2)\cup v_3\Omega u_3-\{u_3,v_3\}$ be
distinct such that neither of the sets $V(P_1)\cup v_3\Omega u_1$,
$V(P_2)\cup v_2\Omega u_3$ includes both $q_1$ and $q_2$.
Let $Q_1$ and $Q_2$ be  two not necessarily disjoint paths with
one end in $u_3\Omega v_3-Z-\{u_3,v_3\}$ and the other end $q_1$ and $q_2$,
respectively, both internally disjoint from $V(P_1\cup P_2\cup L)\cup V(\Omega)$.
In those circumstances we
say that $P_1\cup P_2\cup L\cup Q_1\cup Q_2$ is a {\em turtle}
in $(G,\Omega)$.
We say that $P_1,P_2$ are the {\em legs},
$L$ is the {\em neck}, and $Q_1\cup Q_2$ is the {\em body} of the turtle. (See Figure~\ref{fig:bumps}(a),(b).)

Let $(G,\Omega)$ be a society, let $(u_1,u_2,u_3, v_1,v_2,v_3)$ be
clockwise in $\Omega$, and let $P_1,P_2,P_3$ be disjoint bumps
such that $P_i$ has ends $u_i$ and $v_i$.  In those 
circumstances we say that $P_1,P_2,P_3$ are {\em three crossed paths}
in $(G,\Omega)$.

Let $(G,\Omega)$ be a society, and let 
$u_1,u_2,u_3,u_4,v_1,v_2,v_3,v_4\in V(\Omega)$ be such that either
$(u_1,u_2,u_3,v_2,u_4,v_1,v_4,v_3)$ or
$(u_1,u_2,u_3,u_4,v_2,v_1,v_4,v_3)$ or
$(u_1,u_2,u_3,v_2=u_4,v_1,v_4,v_3)$ is clockwise.
For $i=1,2,3,4$ let $P_i$ be a bump with ends $u_i$ and $v_i$ such that
these bumps are pairwise disjoint, except possibly for $v_2=u_4$.
In those circumstances we say that $P_1,P_2,P_3,P_4$ is 
a {\em gridlet}. (See Figure~\ref{fig:bumps}(c),(d).)

Let $(G,\Omega)$ be a society and let $(u_1,u_2,v_1,v_2,u_3,u_4,v_3,v_4)$ 
be clockwise or counter-clockwise in $\Omega$. 
For $i=1,2,3,4$ let $P_i$ be a bump with ends $u_i$ and $v_i$ such that
these bumps are pairwise disjoint, and let $P_5$ be a path with one
end in $V(P_1)\cup v_4\Omega u_2-\{u_2,v_1,v_4\}$, the other
end in $V(P_3)\cup v_2\Omega u_4-\{v_2,v_3,u_4\}$, and otherwise
disjoint from $P_1\cup P_2\cup P_3\cup P_4$.
In those circumstances we say that $P_1,P_2,\ldots,P_5$ is
a {\em separated doublecross}.(See Figure~\ref{fig:bumps}(e),(f).)

A society $(G,\Omega)$ is {\em rural} if $G$ can be drawn in a disk
with $V(\Omega)$ drawn on the boundary of the disk in the order given
by $\Omega$.
A society $(G,\Omega)$ is {\em nearly rural} if there exists a vertex
$v\in V(G)$ such that the society $(G\backslash v,\Omega\backslash v)$
obtained from $(G,\Omega)$ by deleting $v$ is rural.

In Sections~\ref{sec:leap}--\ref{sec:lack} we prove the following.
The proof strategy is explained in Section~\ref{sec:transactions}.
It uses a couple of theorems from~\cite{RobSeyGM9} and  
Theorem~\ref{thm:leap} that we prove in Section~\ref{sec:leap}.

\begin{theorem}
\label{thm:society}
\showlabel{thm:society}
There exists an integer $k\ge1$ such that for every 
integer $s\ge0$ and every $6$-connected $s$-nested
$k$-cosmopolitan society $(G,\Omega)$ either $(G,\Omega)$ is  nearly rural, 
 or $G$ has a triangle $C$ such that
$(G\backslash E(C),\uom)$ is  rural,
or $(G,\Omega)$ has an $s$-nested planar
truncation that has a turtle,  three crossed paths, a gridlet, 
or a separated doublecross.
\end{theorem}

Finally, we need to convert a turtle, three crossed paths,
gridlet and a separated double-cross 
into a $K_6$ minor.
Let $G$ be a $6$-connected graph, let $H$ be a sufficiently high planar wall
in $G$, and let $(J,\Omega)$ be the anticompass society of $H$ in $G$.
We wish to apply to Theorem~\ref{thm:society} to $(J,\Omega)$. We can,
in fact, assume that $H$ is a subgraph of a larger planar wall $H'$ that
includes $s$ concentric cycles $C_1,C_2,\ldots,C_s$
surrounding $H$ and disjoint from $H$,
for some suitable integer $s$, and hence $(J,\Omega)$ is
$s$-nested.
Theorem~\ref{thm:society} guarantees a turtle 
or paths in $(J,\Omega)$ forming three crossed paths, a gridlet,
or a separated double-cross, but it 
does not say how the turtle or paths might intersect the cycles $C_i$.
In Section~\ref{sec:nest} we prove a theorem that says that the cycles and
the turtle (or paths) can be changed such that after possibly sacrificing 
a lot of
the cycles, the remaining cycles and the new turtle (or paths)
intersect nicely.
Using that information it is then easy to find a $K_6$ minor in $G$.
We complete the proof of Theorem~\ref{thm:largejorgensen} in
Section~\ref{sec:turtle}.

\section{Finding a planar wall}
\label{sec:planarwall}
\showlabel{sec:planarwall}

\begin{figure}[ht!]
\begin{center}
\leavevmode
\includegraphics[trim = 0mm 5mm 65mm 0mm, clip]{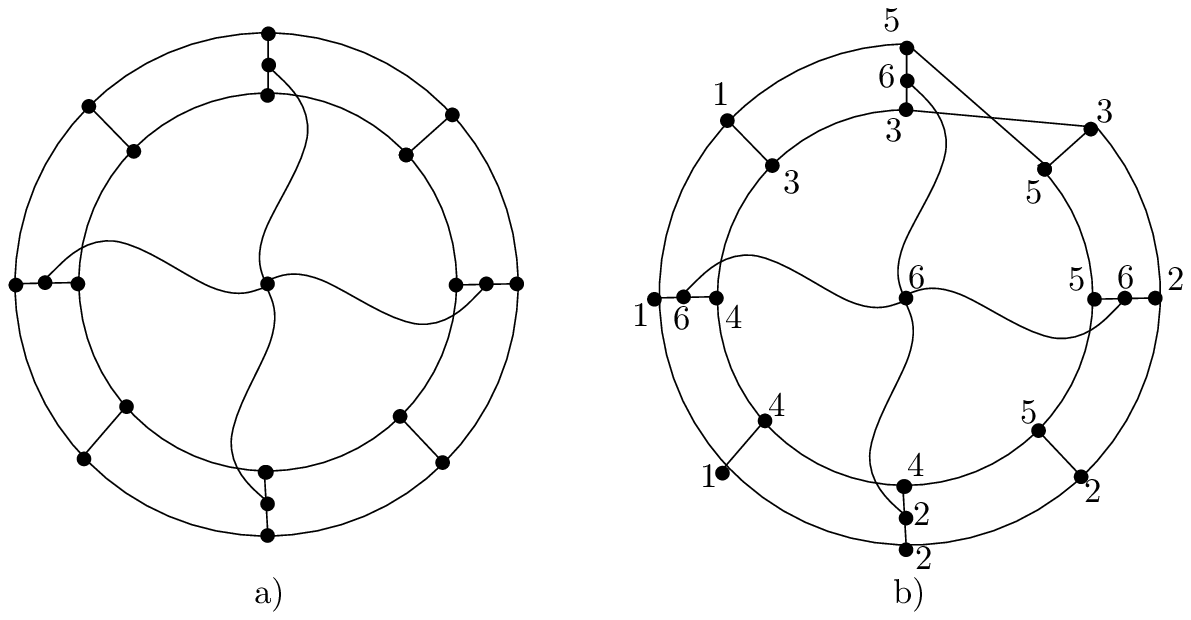}
\end{center}
\caption{
A pinwheel with four vanes.
}
\label{fig:pinwheel}
\end{figure}

Let a {\em pinwheel with four vanes} be the graph pictured in
Figure~\ref{fig:pinwheel}. We define a pinwheel with $k$ vanes analogously.
A graph $G$ is {\em \ifc} if it is simple, $3$-connected, has at least
five vertices, and for every separation $(A,B)$ of $G$ of
order three, one of $A,B$ induces a graph with at most three edges.

\showfiglabel{fig:pinwheel}

The objective of this section is to prove the following theorem.

\begin{theorem}
\label{plsepwall}
\showlabel{plsepwall}
For every even integer $t\ge2$ there exists an even integer $h$ such that if
$H$ is a wall of height at least $h$ in an \ifc\ graph $G$, then either
\myitem{(1)}$G$ has a $K_6$ minor, or
\myitem{(2)}$G$ has a subgraph isomorphic to a subdivision of a
pinwheel with $t$ vanes, or
\myitem{(3)}$G$ has a planar wall of height $t$.
\end{theorem}

In the proof we will be using several results from~\cite{RobSeyGM13}.
Their statements require the following terminology:
distance function,  $(l,m)$-star over $H$,
external $(l,m)$-star over $H$, subwall, dividing subwall,
flat subwall, cross over a wall.
We refer to~\cite{RobSeyGM13} for precise definitions, but we offer
the following informal descriptions.
The distance of two distinct vertices $s,t$ of a wall is the minimum number
of times a curve in the plane joining $s$ and $t$ intersects the drawing
of the wall, when the wall is drawn as in Figure~\ref{fig:wall}.
An $(l,m)$-star over a wall $H$ in $G$  is a subdivision of a star
with $l$ leaves such that only the leaves and possibly the center
belong to $H$, and the leaves are pairwise at distance at least $m$.
The star is external if the center does not belong to $H$.
A subwall of a wall is dividing if its perimeter separates the subwall from
the rest of the wall. A cross over a wall is a set of two disjoint paths
joining the diagonally opposite pairs of ``corners" of the wall,
the vertices represented by solid circles in Figure~\ref{fig:wall}.
A subwall $H$ is flat in $G$ if there is no cross $P,Q$ over $H$
such that $P\cup Q$ is a subgraph of the compass of $H$ in $G$.

We begin with the following easy lemma. We leave the proof to the reader.

\begin{lemma}
\label{usestar}
\showlabel{usestar}
For every integer $t$ there exist integers $l$ and $m$ such that
if a graph $G$ has a wall $H$ with an external $(l,m)$-star, then
it has a subgraph isomorphic to a subdivision of a pinwheel with $t$ vanes.
\end{lemma}

We need one more lemma, which follows immediately
from~\cite[Theorem~8.6]{RobSeyGM13}.

\begin{lemma}
\label{crossinwall}
\showlabel{crossinwall}
Every flat wall in an \ifc\ graph is planar.
\end{lemma}

\begin{figure}[ht!]
\begin{center}
\leavevmode
\includegraphics[scale = 1]{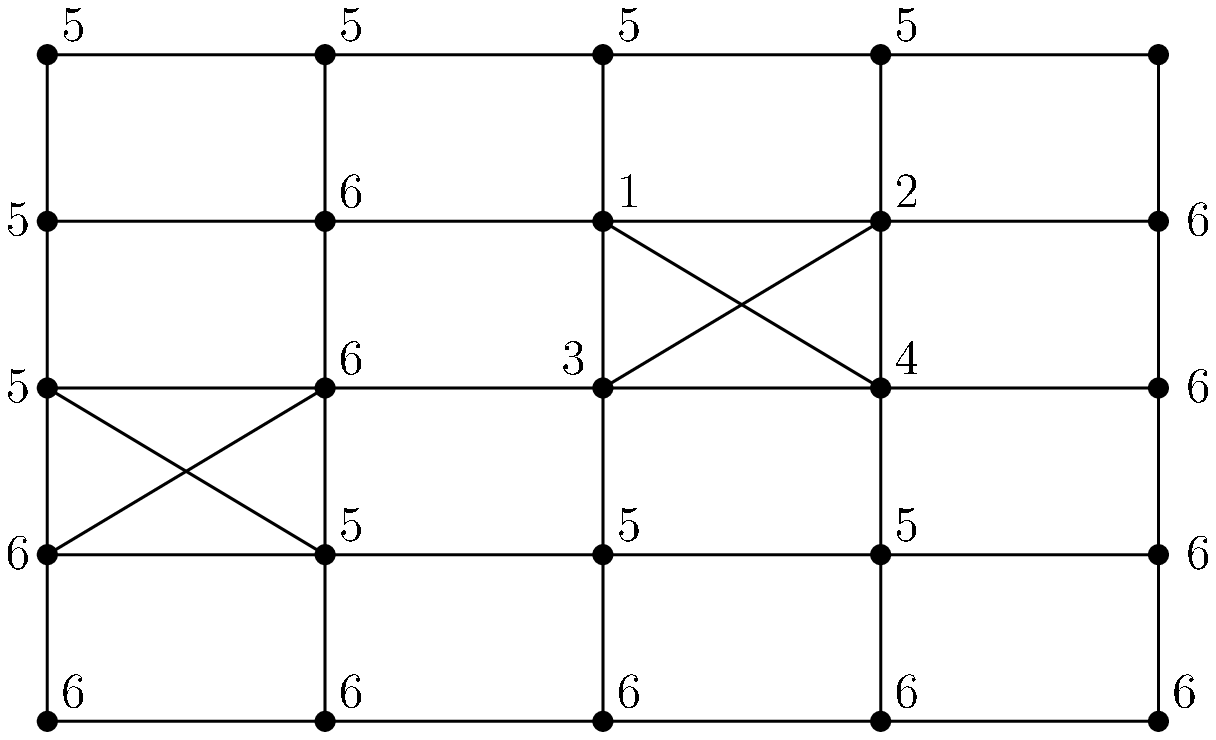}
\end{center}
\caption{
A $K_6$ minor in a grid with two crosses.
}
\label{fig:gridk6}
\end{figure}
\showfiglabel{fig:pinwheel}

\noindent
{\bf Proof of Theorem~\ref{plsepwall}}.
Let $t\ge1$ be given,  let $l,m$ be as in Lemma~\ref{usestar},
let $p=6$, and let $k,r$ be as in~\cite[Theorem 9.2]{RobSeyGM13}.
If $h$ is sufficiently large, then $H$ has
$k+1$ subwalls of height at least $t$, pairwise at distance at least $r$.
If at least $k$ of these subwalls are non-dividing, then
by \cite[Theorem 9.2]{RobSeyGM13} $G$ either has a $K_6$ minor, or
an $(l,m)$-star over $H$, in which case it has a subgraph isomorphic
to a pinwheel with $t$ vanes by Lemma~\ref{usestar}. In either case
the theorem holds, and so we may assume that at least two of the
subwalls, say $H_1$ and $H_2$, are dividing.
We may assume that $H_1$ and $H_2$ are not planar,
for otherwise the theorem holds.
Let $i\in\{1,2\}$.
By Lemma~\ref{crossinwall} the wall $H_i$  is not flat, and hence
its compass has a cross  $P_i\cup Q_i$.
Since the subwalls $H_1$ and $H_2$ are dividing, it follows that
the paths $P_1,Q_1,P_2,Q_2$ are pairwise disjoint.
Thus $G$ has a minor isomorphic to the graph shown in Figure~\ref{fig:gridk6},
but that graph has a minor isomorphic to a minor of $K_6$,
as indicated by the numbers in the figure.
Thus $G$ has a $K_6$ minor, and the theorem holds.~\qed

To deduce Theorem~\ref{thm:plwall} we need the following lemma, proved
in~\cite[Lemma~5.3]{KawNorThoWolbdtw}.

\begin{lemma}
\label{lem:pinwheel}
\showlabel{lem:pinwheel}
If a $5$-connected graph $G$ with no $K_6$ minor has a subdivision
isomorphic to a pinwheel with $20$ vanes, then $G$ is apex.
\end{lemma}

\noindent
{\bf Proof of Theorem~\ref{thm:plwall}.}
Let $t\ge2$ be an even integer. We may assume that $t\ge 20$.
Let $h$ be as in Theorem~\ref{plsepwall}, and let $G$ be a $5$-connected
graph with no $K_6$ minor. From Theorem~\ref{plsepwall} we deduce that $G$ either
satisfies the conclusion of Theorem~\ref{thm:plwall}, or has a
subdivision isomorphic to a pinwheel with $t$ vanes. 
In the latter case the theorem follows from Lemma~\ref{lem:pinwheel}.~\qed


We need the following theorem of DeVos and Seymour~\cite{DevSeyExt3col}.

\begin{theorem}
\label{devosseymour}
\showlabel{devosseymour}
Let $(G,\Omega)$ be a rural society such that $G$ is a simple graph
and every vertex of $G$ not in $V(\Omega)$ has degree at least six.
Then $|V(G)|\le |V(\Omega)|^2/{12}+|V(\Omega)|/2+1$.
\end{theorem}



\noindent
{\bf Proof of Theorem~\ref{cosmopolitan}.}
 Let $k\ge 1$ be an integer, and let $t$ be an even integer such that
if $W$ is the elementary wall of height $t$ and $|V(W)|\le\ell^2/12+\ell/2+1$, then
$\ell >6k-6$.  Let $K$ be the compass of $H$ in $G$, let $(J,\uom)$
be the anticompass society of $H$ in $G$, let $(G_0,\uom_0)$ be a planar
truncation of $(J,\uom)$, and let $\ell=|V(\uom_0)|$.  
Thus $(J,\uom)$ is the composition of $(G_0,\uom_0)$ with a rural
neighborhood $(G',\uom,\uom_0)$.  Then $|V(H)|\le\ell^2/12 + \ell/2+1$
by Theorem~\ref{devosseymour} applied to the society $(K\cup G',\uom_0)$,
and hence $\ell >6k-6$.  Let $L$ be the graph obtained from $K\cup G'$ 
by adding
a new vertex $v$ and joining it to every vertex of $V(\uom_0)$ and
by adding an edge joining every pair of nonadjacent vertices of $V(\uom_0)$
that are consecutive in $\uom_0$.  Then $L$
is planar.  Let $s$ be the number of vertices of $V(\Omega_0)$ with at least
two neighbors in $G_0$.  Then all but $s$ vertices of $K\cup G'$ have
degree in $L$ at least six.  Thus the sum of the degrees of vertices of
$L$ is at least $6|V(K\cup G')|-6s+\ell$.  On the other hand, the sum of
the degrees is at most $6|V(L)|-12$, because $L$ is planar, and hence
$s\ge k$, as desired.~\qed

\section{Rural societies}

If $P$ is a path and $x,y\in V(P)$, we denote by $xPy$ the unique
subpath of $P$ with ends $x$ and $y$.
Let $(G,\Omega)$ be a society.
An {\em orderly transaction}  in  $(G,\Omega)$ is a sequence of $k$
pairwise disjoint bumps ${\cal T}=(P_1,\dots, P_k)$ such that $P_i$ has
ends $u_i$ and $v_i$ and $u_1,u_2,\dots, u_k, v_k,v_{k-1},\dots,
v_1$ is clockwise.
Let $M$ be the graph obtained from
$P_1\cup P_2\cup\cdots\cup P_k$ by adding
the vertices of $V(\Omega)$ as isolated vertices.
We say that $M$ is the {\em frame} of $\cal T$.
We say that a path $Q$ in $G$ is {\em $\cal T$-coterminal} if
$Q$ has both ends in $V(\Omega)$ and is otherwise disjoint from it and
for every $i=1,2,\ldots,k$ the following holds: if $Q$ intersects $P_i$,
then their intersection is a path whose one end is a common end of
$Q$ and $P_i$.

Let $(G,\Omega)$ be a society, and let $M$ and $\cal T$ be as in
the previous paragraph.
Let $i\in\{1,2,\ldots,k\}$ and let $Q$ be
a $\cal T$-coterminal path in  $G\backslash V(P_i)$ with one end in
$v_i\Omega u_i$  and the other end in $u_i\Omega v_i$.
In those circumstances we say that $Q$ is a {\em $\cal T$-jump over $P_i$},
or simply a {\em $\cal T$-jump}.

Now let $i\in\{0,1,\ldots,k\}$ and let $Q_1,Q_2$ be two disjoint
$\cal T$-coterminal paths such that $Q_j$ has ends $x_j,y_j$
and $(u_{i},x_1,x_2,u_{i+1},v_{i+1},y_1,y_2,v_{i})$ is clockwise in $\Omega$,
where possibly $u_{i}=x_1$, $x_2=u_{i+1}$, $v_{i+1}=y_1$, or $y_2=v_{i}$,
and $u_0$ means $x_1$, $u_{k+1}$ means $x_2$, $v_{k+1}$ means $y_1$,
and $v_0$ means $y_2$.
In those circumstances we say that $(Q_1,Q_2)$ is a {\em $\cal T$-cross
in region $i$},
or simply a {\em $\cal T$-cross}.

Finally, let $i\in\{1,2,\ldots,k\}$ and let  $Q_0$, $Q_{1}$, $Q_2$
be three paths such that $Q_j$ has ends $x_j,y_j$ and is otherwise
disjoint from all members of $\cal T$,
$x_0,y_0\in V(P_i)$,
the vertices $x_1,x_2$ are internal vertices of $x_0P_i y_0$,
$y_1,y_2\not\in V(P_i)$,
$y_1\in u_{i-1}\Omega u_i\cup  v_{i}\Omega v_{i-1}$,
$y_2\in u_{i}\Omega u_{i+1}\cup  v_{i+1}\Omega v_{i}$,
and the paths $Q_0$, $Q_{1}$, $Q_2$ are pairwise disjoint, except possibly
$x_1=x_2$.
In those circumstances we say that $(Q_0,Q_{1}, Q_2)$ is a
{\em $\cal T$-tunnel
under $P_i$}, or simply a {\em $\cal T$-tunnel}.

Intuitively, if we think of the paths in $\cal T$ as dividing the
society into ``regions", then a $\cal T$-jump arises from a $\cal T$-path
whose ends do not belong to the same region.
A $\cal T$-cross arises from two $\cal T$-paths with ends in the
same region that cross inside that region, and furthermore,
each path in $\cal T$ includes at most two ends of those crossing paths.
Finally, a $\cal T$-tunnel can be converted into a $\cal T$-jump
by rerouting $P_i$ along $Q_0$.
However, in some applications such rerouting will be undesirable,
and therefore we need to list $\cal T$-tunnels as outcomes.

Let $M$ be a subgraph of a graph $G$.  An {\em $M$-bridge} in $G$ is a
connected subgraph $B$ of $G$ such that $E(B)\cap E(M)=\emptyset$
and either $E(B)$ consists of a unique edge with both ends in $M$, or for
some component $C$ of $G\backslash V(M)$ the set $E(B)$ consists of all
edges of $G$ with at least one end in $V(C)$.  The vertices in
$V(B)\cap V(M)$ are called the {\em attachments} of $B$.  
Now let $M$ be  such that no block of $M$ is  a cycle.
By a {\em segment of $M$} we mean a maximal subpath $P$ of $M$ such that
every internal vertex of $P$ has degree two in $M$.
It follows that the segments 
of $M$ are uniquely determined.  Now if $B$ is an $M$-bridge of $G$, then we
say that $B$ is {\em unstable} if some segment of $M$ includes all the
 attachments of $B$, and otherwise we say that $B$ is
{\em stable}.

A society $(G,\Omega)$ is {\em rurally $4$-connected} if for
every separation $(A,B)$ of order at most three with
$V(\Omega)\subseteq A$ the graph $G[B]$ can be drawn in a disk
with the vertices of $A\cap B$ drawn on the boundary of the disk.
A society is {\em cross-free} if it has no cross.
The following, a close relative of Lemma~\ref{crossinwall},
follows from \cite[Theorem~2.4]{RobSeyGM9}.

\begin{theorem}
\label{2paththm}
\showlabel{2paththm}
Every cross-free rurally $4$-connected society is rural.
\end{theorem}

\begin{lemma}
\mylabel{orderlyrural1}
Let $(G,\Omega)$ be a rurally $4$-connected society,
let ${\cal T}=(P_1,\dots, P_k)$ be an orderly transaction in  $(G,\Omega)$,
and let $M$ be the frame of $\cal T$.
If every $M$-bridge of $G$ is stable and $(G,\Omega)$  is not
rural, then $(G,\Omega)$ has a $\cal T$-jump, a $\cal T$-cross, or a $\cal T$-tunnel.
\end{lemma}

\proof
For $i=1,2,\ldots,k$ let $u_i$ and $v_i$ be the ends of $P_i$ numbered
as in the definition of orderly transaction, and for
convenience let $P_0$ and $P_{k+1}$ be null graphs.
We define $k+1$
cyclic permutations $\Omega_0,\Omega_1,\ldots,\Omega_k$ as follows.
For $i=1,2,\ldots,k-1$ let $V(\Omega_i):= V(P_{i})\cup V(P_{i+1})\cup
u_i\Omega u_{i+1}\cup v_{i+1}\Omega v_i$ with the cyclic order
defined by saying that $u_i\Omega u_{i+1}$ is followed by $V(P_{i+1})$
in order from $u_{i+1}$ to $v_{i+1}$, followed by $v_{i+1}\Omega v_i$
followed by $V(P_{i})$ in order from $v_i$ to $u_i$.
The cyclic permutation $\Omega_0$ is defined by letting $v_1\Omega u_1$
be followed by $V(P_1)$ in order from $u_1$ to $v_1$,
and $\Omega_k$ is defined by letting $u_k\Omega v_k$
be followed by $V(P_k)$ in order from $v_k$ to $u_k$.

Now if for some $M$-bridge $B$ of $G$ there is no index $i\in\{0,1,\ldots,k\}$
such that all attachments of $B$ belong to $V(\Omega_i)$, then
$(G,\Omega)$ has a $\cal T$-jump.
Thus we may assume that such index exists for every $M$-bridge $B$,
and since $B$ is stable that index is unique. Let us denote it by $i(B)$.
For $i=0,1,\ldots,k$
let $G_i$ be the subgraph of $G$ consisting of $P_i\cup P_{i+1}$,
the vertex-set $V(\Omega_i)$ and all $M$-bridges $B$ of $G$ with
$i(B)=i$.
The society $(G_i,\Omega_i)$ is rurally $4$-connected.
If each $(G_i,\Omega_i)$ is cross-free,
then each of them is rural by Theorem~\ref{2paththm}
and it follows that  $(G,\Omega)$ is rural.
Thus we may assume that for some $i=0,1,\ldots,k$ the society
$(G_i,\Omega_i)$ has a cross $(Q_1,Q_2)$.
If neither $P_i$ nor $P_{i+1}$ includes three or four ends of the
paths $Q_1$ and $Q_2$, then $(G,\Omega)$ has a $\cal T$-cross.
Thus we may assume that $P_i$ includes both ends of $Q_1$ and
at least one end of $Q_2$.
Let $x_j,y_j$ be the ends of $Q_j$. Since the $M$-bridge of $G$ containing
$Q_2$ is stable, it has an attachment outside $P_i$,
and so if needed, we may replace $Q_2$ by a path with an end outside $P_i$
(or conclude that $(G,\Omega)$ has a $\cal T$-jump). 
Thus we may assume that $u_i,x_1,x_2,y_1,v_i$ occur on $P_i$ in
the order listed, and $y_2\not\in V(P_i)$.

The $M$-bridge of $G$ containing $Q_1$ has an attachment outside $P_i$.
If it does not include $Q_2$ and has an attachment outside $V(P_i)\cup\{y_2\}$, then $(G,\Omega)$ has a $\cal T$-jump or $\cal T$-cross,
and so we may assume not. Thus there exists a path $Q_3$ with
one end $x_3$ in the interior of $Q_1$ and the other end  $y_3\in V(Q_2)-\{x_2\}$
with no internal vertex in $M\cup Q_1\cup Q_2$.
We call the triple $(Q_1,Q_2,Q_3)$ a {\em tripod}, and the path
$y_3Q_2y_2$ the {\em leg} of the tripod.
If $v$ is an internal vertex of $x_1P_iy_1$, then we say that
$v$ is {\em sheltered} by the tripod $(Q_1,Q_2,Q_3)$.
Let $L$ be a path that is the leg of some tripod, and subject to that
$L$ is minimal.
From now on we fix $L$ and  will consider different tripods with leg $L$;
thus the vertices $x_1,y_1,x_2,x_3$ may change, but $y_2$ and $y_3$
will remain fixed as the ends of $L$.

Let $x_1',y_1'\in V(P_i)$ be such that they are sheltered by no tripod
with leg $L$, but every internal vertex of $x_1'P_iy_1'$ is sheltered
by some tripod with leg $L$.
Let $X'$ be the union of $x_1'P_iy_1'$ and all tripods with leg $L$
that shelter some internal vertex of $x_1'P_iy_1'$,
let $X:=X'\backslash V(L)\backslash\{x_1',y_1'\}$
and
let $Y:=V(M\cup L)- x_1'P_iy_1'-\{y_3\}$.
Since $(G,\Omega)$ is rurally $4$-connected we deduce that
the set $\{x_1',y_1',y_3\}$ does not separate $X$ from $Y$ in $G$.
It follows that there exists a path $P$ in $G\backslash \{x_1',y_1',y_3\}$
with ends $x\in X$ and $y\in Y$.
We may assume that $P$ has no internal vertex in $X\cup Y$.
Let $(Q_1,Q_2,Q_3)$ be a tripod with leg L such that either $x$ is sheltered
by it, or $x\in V(Q_1\cup Q_2\cup Q_3)$.
If $y\not\in V(L\cup P_i)$, then by considering the paths $P,Q_1,Q_2,Q_3$
it follows that either $(G,\Omega)$ has a $\cal T$-jump or $\cal T$-tunnel.
If $y\in V(L)$, then there is a tripod whose leg is a proper subpath of
$L$, contrary to the choice of $L$.
Thus we may assume that $y\in V(P_i)$, and that $y\in V(P_i)$ for every
choice of the path $P$ as above.
If $x\in V(Q_1\cup Q_2\cup Q_3)$ then there is a tripod with leg $L$
that shelters $x_1'$ or $y_1'$, a contradiction. Thus $x\in V(P_i)$.
Let $B$ be the $M$-bridge containing $P$.
Since $y\in V(P_i)$ for all choices of $P$ it follows that the attachments
of $B$ are a subset of $V(P_i)\cup\{y_2\}$.
But $B$ is stable, and hence $y_2$ is an attachment of $B$.
The minimality of $L$ implies that $B$ includes a path from $y$ to $y_3$,
internally disjoint from $L$.
Using that path and the paths $P,Q_1,Q_2,Q_3$ it is now easy to
construct a tripod that shelters either $x_1'$ or $y_1'$, a contradiction.~\qed


\section{Leap of length five}
\label{sec:leap}
\showlabel{sec:leap}

A {\em leap of length} $k$ in a society $(G,\uom)$ is a sequence of $k+1$
pairwise disjoint bumps $P_0,P_1,\dots, P_k$ such that $P_i$ has
ends $u_i$ and $v_i$ and $u_0,u_1,u_2,\dots, u_k, v_0, v_k,v_{k-1},\dots,
v_1$, is clockwise. In this section we prove the following.

\begin{theorem}\label{thm:leap}
\showlabel{thm:leap}
Let $(G,\uom)$ be a $6$-connected society with a leap of length five.
Then $(G,\uom)$ is nearly rural, or $G$ has a triangle $C$ such that
$(G\backslash E(C),\uom)$ is  rural, or $(G,\uom)$ has three crossed paths,
a gridlet, a separated doublecross, or a turtle.
\end{theorem}

The following is a hypothesis that will be common to several lemmas
of this section, and so we state it separately to avoid repetition.

\begin{hypothesis}
\label{hyp:leap}
\showlabel{hyp:leap}
Let  $(G,\uom)$ be a society with no three crossed paths,
a gridlet, a separated doublecross, or a turtle,
let $k\ge1$ be an integer, let
$$(u_0,u_1,u_2,\ldots,u_k,v_0,v_k,v_{k-1},\ldots,v_1)$$
be clockwise, and let $P_0,P_1,\ldots,P_k$ be pairwise disjoint bumps
such that $P_i$ has ends $u_i$ and $v_i$.
Let $\cal T$ be the orderly transaction $(P_1,P_2,\ldots,P_k)$, let $M$ be the frame of $\cal T$ and let
$$Z=u_1\uom u_k\cup v_k\uom v_1\cup V(P_2)\cup V(P_3)\cup\cdots\cup
V(P_{k-1})-\{u_1,u_k,v_1,v_k\}.$$
Let $Z_1=v_1\uom u_1-\{u_0,u_1,v_1\}$ and $Z_2=u_k\uom v_k-\{v_0,u_k,v_k\}$.
\end{hypothesis}

If $H$ is a subgraph of $G$, then an \emph{$H$-path} is a (possibly trivial) path with both ends in $V(H)$ and otherwise disjoint from $H$.
This is somewhat non-standard, typically an $H$-path is required to have at least one edge, but we use our definition for convenience.
We say that a  vertex $v$ of $P_0$ is {\em exposed} if
there exists an $(M\cup P_0)$-path $P$ with one end $v$ and the other in $Z$.

\begin{lemma}
\label{lem:2pathstoZ}
\showlabel{lem:2pathstoZ}
Assume Hypothesis~\ref{hyp:leap} and let $k\ge 3$.
Let $R_1,R_2$ be two disjoint $(M\cup P_0)$-paths in $G$ such that $R_i$
has ends $x_i\in V(P_0)$ and $y_i\in V(M)-\{u_0,v_0\}$, and assume that
$u_0,x_1,x_2,v_0$ occur on $P_0$ in the order listed, where possibly $u_0=x_1$, or $v_0=x_2$, or both.
Then either $y_1\in V(P_1)\cup v_1\Omega u_1$, or
$y_2\in V(P_k)\cup u_k\Omega v_k$, or both.
In particular, there do not  exist two disjoint $(M\cup P_0)$-paths from
$V(P_0)$ to $Z$.
\end{lemma}

\proof The second statement follows immediately from the first,
and so it suffices to prove the first statement.
Suppose for a contradiction that there exist paths  $R_1,R_2$
satisfying the hypotheses but not the conclusion of the lemma.
By using the paths $P_2, P_3, \ldots, P_{k-1}$ we conclude that
there exist two disjoint paths
$Q_1,Q_2$  in $G$ such that $Q_i$ has ends
$x_i\in V(P_0)$ and $z_i \in V(\Omega)$,
and is otherwise disjoint from $V(P_0)\cup V(\Omega)$,
and if $Q_i$ intersects some $P_j$ for $j\in\{1,2,\ldots,k\}$,
then $j\in\{2,\ldots,k-1\}$ and $Q_i \cap P_j$ is a path one of whose ends is a common end
of $Q_i$ and $P_j$. Furthermore, $z_1 \in u_1 \uom v_1-\{u_1,v_1\} $ and $z_2 \in v_k \uom u_k-\{u_k,v_k\}$.
From the symmetry we may assume that either
$(u_0,v_0,z_2,z_1)$, or $(u_0,z_1,v_0,z_2)$
or $(u_0,v_0,z_1,z_2)$ is clockwise.
In the first two cases  $(G,\uom)$ has a separated doublecross
(the two pairs of crossing bumps are $P_1$ and
$Q_1\cup u_0P_0x_1$, and $P_k$ and $Q_2\cup v_0P_0x_2$, and
the fifth path is a subpath of $P_2$), unless the second case holds and
$z_1 \in u_k \uom v_0$ or $z_2 \in v_1 \uom u_0$, or both. By symmetry we may assume that
$z_1 \in u_k \uom v_0$. Then, if $z_2 \in v_{k-2} \uom u_0$, $(G,\uom)$ has a gridlet formed by the paths $P_k,P_{k-1},
u_0P_0x_1 \cup Q_1$ and $v_0P_0x_2 \cup Q_2$. Otherwise, $z_2 \in v_k\uom v_{k-2} - \{v_k,v_{k-2}\}$ and $(G,\uom)$ has a turtle
with legs $P_k$ and $v_0P_0x_2 \cup Q_2$, neck $P_1$ and body $u_0P_0x_2 \cup Q_1$.

Finally, in the third case
$(G,\uom)$ has a turtle  or three crossed paths.
More precisely, if $z_2 \in v_0 \uom v_1 - \{v_1\}$, then $(G,\uom)$ has a turtle described in the paragraph above.
Otherwise, by symmetry, we may assume that $z_2 \in v_1 \uom u_0$ and $z_1 \in v_0 \uom v_k$, in which case $v_0P_0x_2 \cup Q_2$,
$u_0P_0x_1 \cup Q_1$ and $P_2$ are the three crossed paths.~\qed

\begin{lemma}
\mylabel{prejump}
Assume Hypothesis~\ref{hyp:leap} and let $k\ge2$.
Then $(G\backslash V(P_0),\Omega\backslash V(P_0))$ has no $\cal T$-jump.
\end{lemma}

\proof Suppose for a contradiction that
$(G\backslash  V(P_0),\Omega\backslash  V(P_0))$
has a $\cal T$-jump.
Thus there is an index $i\in\{1,2,\ldots,k\}$ and a $\cal T$-coterminal
path $P$ in $G\backslash V(P_0\cup P_i)$ with ends
$x\in v_i\Omega u_i$ and $y\in u_i\Omega v_i$.
Let $j\in\{1,2,\ldots,k\}-\{i\}$.
Then using the paths $P_0,P_i,P_j$ and $P$ we deduce that $(G,\Omega)$ has either
three crossed paths or a gridlet,
in either case a contradiction.~\qed

\begin{lemma}
\mylabel{lem:nojump}
Assume Hypothesis~\ref{hyp:leap} and let $k\ge2$.
Let $v\in V(P_0)$ be such that there is no $(M\cup P_0)$-path in $G\backslash v$
from $vP_0v_0$ to $vP_0u_0\cup V(P_1\cup P_2\cup\cdots\cup P_{k-1})
\cup v_k\Omega u_k-\{v_k,u_k\}$
and none from $vP_0u_0$ to
$V(P_2\cup P_3\cup\cdots\cup P_{k})
\cup u_1\Omega v_1-\{u_1,v_1\}$.
Then $(G\backslash v,\Omega\backslash v)$ has no $\cal T$-jump.
\end{lemma}

\proof
The hypotheses of the lemma imply that every $\cal T$-jump in $(G\backslash v,\Omega\backslash v)$ is disjoint from $P_0$.
Thus the lemma follows from Lemma~\ref{prejump}.~\qed

\begin{lemma}
\mylabel{lem:nocross}
Assume Hypothesis~\ref{hyp:leap},  let $k\ge3$, and let
$v\in V(P_0)$ be such that no vertex in $V(P_0)-\{v\}$ is exposed.
Let $i\in\{0,1,\ldots,k\}$ be such that
$(G\backslash v,\Omega\backslash v)$ has
a $\cal T$-cross $(Q_1,Q_2)$ in region $i$.
Then $i\in\{0,k\}$ and $v$ is not exposed.
Furthermore, assume that $i=0$, and that there exists an $(M\cup P_0)$-path $Q$
with one end $v$ and the other end in $P_1\cup v_1\Omega u_1-\{u_0\}$,
and that $v_0P_0v$ is disjoint from $Q_1\cup Q_2$. Then for some $j\in\{1,2\}$
there exist $p \in V(Q_j\cap u_0P_0v)$ and $q \in V(Q_j\cap Q)$
such that $pP_0v$ and $qQv$ are internally disjoint from $Q_1\cup Q_2$.
\end{lemma}

\proof
If $i\not\in\{0,k\}$, then the $\cal T$-cross is disjoint from
$P_0$ by the choice of $v$, and hence the $\cal T$-cross and $P_0$ give
rise to three crossed paths.
To complete the proof of the first assertion we may assume that
$i=0$ and that $v$ is exposed. Subject to these assumptions we choose $Q_1$ and $Q_2$ so that
$Q_1 \cup Q_2 \cup P_0$ is minimal.  Since $v$ is exposed there exists a $\cal T$-coterminal path $Q'$ from $v$ to
$y \in Z\cap V(\Omega)$ disjoint from $P_0\cup P_1\cup P_k\backslash v$. Let $Q''=Q'\cup vP_0v_0$.
If $Q'' \cap (Q_1 \cup Q_2) = \emptyset$ then $(G,\Omega)$ has a separated doublecross, where
one pair of crossed paths is obtained from the $\cal T$-cross,
the other pair is $P_k$ and $Q''$, and the fifth path is
a subpath of $P_2$. Thus we may assume that there exists $x \in V(Q'') \cap V(Q_j)$ for some $j \in \{1,2\}$ and that
$x$ is chosen so that $xQ''y$ is internally disjoint from $Q_1 \cup Q_2$.  For $r=1,2$ let $z_r \in v_1\uom u_1 - \{v_1,u_1\}$
be an end of $Q_r$ such that $Q_{3-r}$ has one end in $z_r\Omega v_0$ and another in $v_0\Omega z_r$.  If $x \in V(Q')$, then $Q_j$ is disjoint from $P_0$, because $v$ 
is the only exposed vertex and $v \not \in V(Q_1) \cup V(Q_2)$. Thus $z_jQ_jx \cup xQ'y$ is a $\cal T$-jump
disjoint from $P_0$, contrary to Lemma~\ref{prejump}. It follows that $x \in V(v_0P_0v)$, and 
$Q'$ is disjoint from $Q_1 \cup Q_2$. 
 
\begin{figure}
\begin{center}
\leavevmode
\includegraphics[scale = 0.8]{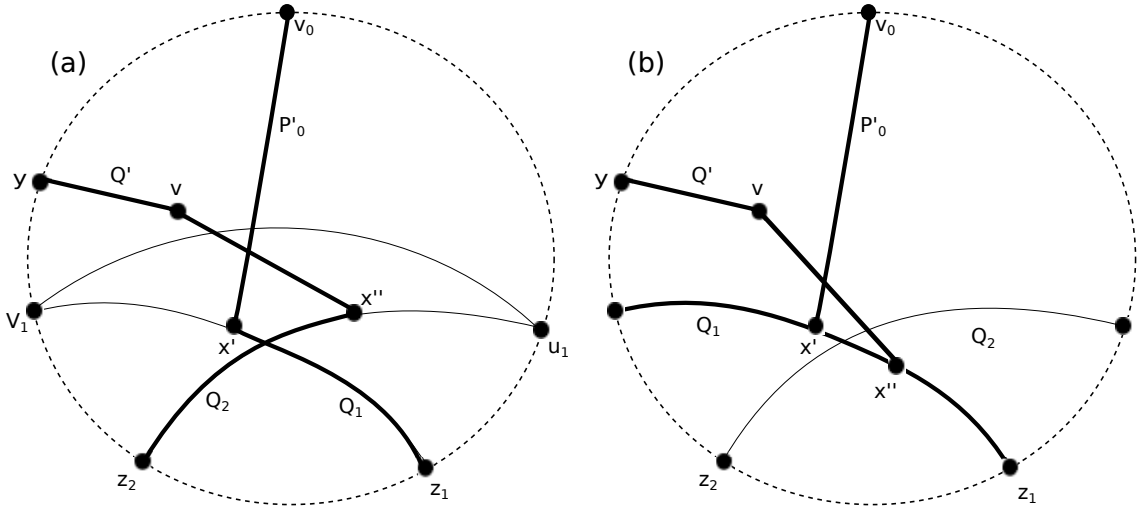}
\end{center}
\caption{
Configurations considered in the proof of the first assertion of Lemma~\ref{lem:nocross}.
}
\label{fig:leaplem46}
\end{figure}

Let $x' \in V(P_0) \cap (V(Q_1) \cup V(Q_2))$ be chosen so that  $x'P_0v_0$ is internally disjoint from $Q_1 \cup Q_2$.   Without loss of generality, we assume that $x' \in V(Q_1)$.  Define $P_0'=v_0P_0x' \cup x'Q_1z_1$. Let $x'' \in vP_0u_0 \cap (V(P'_0) \cup V(Q_2) \cup \{u_0\})$ be chosen so that $vP_0x''$ is internally
disjoint from $P'_0 \cup Q_2$. If $x'' \not \in V(P'_0)$  then the path $Q' \cup vP_0u_0$, if $x'' =u_0$, or the path $Q' \cup vP_0x'' \cup x''Q_2z_2$, if $x'' \in V(Q_2)$ is a $\cal T$-jump, disjoint from $P_0'$, contradicting Lemma~\ref{prejump}. (See Figure~\ref{fig:leaplem46}(a).)  If $x'' \in V(P'_0)$ then $x''P_0v \cup Q'$ and $Q_1 \backslash (V(P_0)-\{x'\})$ are paths with one end in $V(P_0')$ and another in $V(\Omega)$, contradicting Lemma~\ref{lem:2pathstoZ}, after we replace $P_0$ by $P_0'$ and $P_1$ by $Q_2$ in $M$. (See Figure~\ref{fig:leaplem46}(b).) This proves the first assertion of the lemma.


\showfiglabel{fig:leaplem46}

To prove the second statement of the lemma we assume that $i=0$
and that $Q$ is a path from $v$
to $v'\in v_1\Omega u_1-\{u_0\}$, disjoint from $M\cup P_0\backslash v$,
except that $P_1\cap Q$ may be  a path with one end $v'$.
Let
the ends of  $Q_1,Q_2$ be labeled as in the definition of $\cal T$-cross.
If $P_0$ is disjoint from $Q_1\cup Q_2$, then  $(G,\uom)$ has
three crossed paths (if $(y_2,u_0,x_1)$ is clockwise) or
a gridlet with paths $Q_1,Q_2,P_0,P_2$ (if $(x_1,u_0,x_2)$
or $(y_1,u_0,y_2)$ is clockwise), or a separated doublecross
with paths $Q_1,Q_2,P_0,P_2,P_k$ (if $(v_1,u_0,y_1)$
or $(x_2,u_0,u_1)$ is clockwise).
Thus we may assume that  $P_0$ intersects $Q_1\cup Q_2$.
(Please note that $v_0P_0v$ is disjoint from $Q_1\cup Q_2$ by hypothesis.)
Similarly we may assume that  $Q$ intersects  $Q_1\cup Q_2$,
for otherwise we apply the previous argument with $P_0$ replaced by
$Q\cup vP_0v_0$.
Let $p\in V(Q_1\cup Q_2)\cap u_0P_0v$ and $q\in V(Q_1\cup Q_2)\cap V(Q)$
be chosen to minimize $pP_0v$ and $qQv$.
If $p$ and $q$ belong to different paths $Q_1,Q_2$, then
$(G,\uom)$ has a turtle with legs $Q_1,Q_2$, neck
$P_k$ and body $pP_0v_0\cup qQv$.
Thus $p$ and $q$ belong to the same $Q_j$ and the lemma holds.~\qed

In the proof of the following lemma we will be applying Lemma~\ref{orderlyrural1}.
To guarantee that the conditions of Lemma~\ref{orderlyrural1} are satisfied, 
we will need a result from~\cite{KawNorThoWolbdtw}.
We need to precede the statement of this result by a few definitions.

Let $M$ be a subgraph of a graph $G$, such that no block of $M$ is a cycle.
Let $P$ be a segment of $M$ of length at least two, and let $Q$ be a path
in $G$ with ends $x,y\in V(P)$ and otherwise disjoint from $M$.
Let $M'$ be obtained from $M$ by replacing the path $xPy$ by $Q$;
then we say that $M'$ was obtained from $M$ by {\em rerouting} $P$ along $Q$, or
simply that $M'$ was obtained from $M$ by {\em rerouting}.  Please
note that $P$ is required to have length at least two,
and hence this relation is not symmetric.  We say that the
rerouting is {\em proper} if all the attachments
of the $M$-bridge that contains $Q$  belong to $P$.
The following is proved in~\cite[Lemma~2.1]{KawNorThoWolbdtw}.

\begin{lemma}
\mylabel{prestable}
Let $G$ be a graph,
and let $M$ be a subgraph of $G$
such that no block of $M$ is a cycle.
Then there exists a subgraph $M'$ of $G$
 obtained from $M$ by a sequence of proper reroutings
such that if an $M'$-bridge $B$ of $G$ is unstable, say all its attachments
belong to a segment $P$ of $M'$, then there exist
vertices
$x,y\in V(P)$ such that
some component of $G\backslash\{x,y\}$ includes a vertex of $B$
and is disjoint from $M\backslash V(P)$.
\end{lemma}

\begin{lemma}
\mylabel{lem:1exposed}
Assume Hypothesis~\ref{hyp:leap}, and let $k\ge 4$.
If every leap of length $k-1$ has at most one exposed vertex,
 $(G,\uom)$ is $4$-connected and
$(G\backslash v,\uom\backslash v)$ is rurally $4$-connected for
every $v\in V(P_0)$,
then $(G,\uom)$ is nearly rural.
\end{lemma}

\proof
Since  $(G,\uom)$ has no separated doublecross it follows that
it does not have a $\cal T$-cross both in region $0$ and region $k$.
Thus we may assume that it has no $\cal T$-cross in region $k$.
Similarly, it follows that it does not have a $\cal T$-tunnel under both
$P_1$ and $P_k$, or a $\cal T$-cross in region $0$ and a
$\cal T$-tunnel under $P_k$.
Thus we may also assume that  $(G,\uom)$ has no $\cal T$-tunnel under $P_k$.
If some leap of length $k$ in $(G, \uom)$ has an exposed vertex, then we
may assume that $v$ is an exposed vertex.
Otherwise, let the leap $(P_0,P_1,\ldots,P_k)$ and $v\in V(P_0)$
be chosen such that either $v=u_0$ or
there exists an $(M\cup P_0)$-path with one end
$v$ and the other end in $P_1\cup v_1\Omega u_1-\{u_0\}$,
and, subject to that, $vP_0v_0$ is as short as possible.

By Lemma~\ref{prestable} we may assume,
by properly rerouting $M$ if necessary, that every
$M$-bridge of $G\backslash v$ is stable.
Since the reroutings are proper the new paths $P_i$ will still be
disjoint from $P_0$, and the property that defines $v$ will continue to hold.
Similarly, the facts that there is no $\cal T$-cross in region $k$
and no $\cal T$-tunnel under $P_k$ remain unaffected.
We claim that $(G\backslash v,\uom\backslash v)$ is rural.

We apply Lemma~\ref{orderlyrural1} to the society
$(G\backslash v,\uom\backslash v)$ and orderly transaction $\cal T$.
We may assume that $(G\backslash v,\uom\backslash v)$
is not rural, and hence
by Lemma~\ref{orderlyrural1} the society $(G\backslash v,\uom\backslash v)$
has a $\cal T$-jump, a $\cal T$-cross or a $\cal T$-tunnel.
By the choice of $v$ there exists a path $Q$ from $v$ to
$v'\in v_k\Omega u_k-\{v_k,u_k\}$ such that $Q$ does not intersect
$P_k\cup P_0\backslash v$ and intersects at most one of
$P_1,P_2,\ldots,P_{k-1}$. Furthermore, if it intersects $P_i$ for some $i \in \{1,2,\ldots,k-1\}$  then $P_i \cap Q$
is a path with one end a common end of both. (If $v=u_0$ then we can choose $Q$ to be a one vertex path.)

%
We claim that $v$ satisfies the hypotheses of Lemma~\ref{lem:nojump}.
To prove this claim suppose for a contradiction that $P$ is an
$(M\cup P_0)$-path violating that hypothesis.
Suppose first that $P$ and $Q$ are disjoint.
Then $P$ joins different components of $P_0\backslash v$
by Lemma~\ref{lem:2pathstoZ}. But then
changing $P_0$ to the unique path in $P_0\cup P$ that does
not use $v$ either produces a leap with at least two exposed vertices,
or contradicts the minimality of $vP_0v_0$.
Thus $P$ and $Q$ intersect. Since no leap of length $k$ has two or
more exposed vertices, it follows that $v$ is not exposed.
Thus $P$ has one end in $u_0P_0v$ by the minimality of $vP_0v_0$,
and the other end in $P_k\cup u_k\Omega v_k$, because $v$ is not
exposed.
But then $P\cup Q$ includes a $\cal T$-jump disjoint from $P_0$,
contrary to Lemma~\ref{prejump}.
This proves our claim that $v$ satisfies the hypotheses of
Lemma~\ref{lem:nojump}.
We conclude that $(G\backslash v,\uom\backslash v)$ has no
 $\cal T$-jump.

Assume now that
  $(G\backslash v,\uom\backslash v)$ has a $\cal T$-cross $(Q_1,Q_2)$
in region $i$ for some integer $i\in\{0,1,\ldots,k\}$. 
By the first part of Lemma~\ref{lem:nocross} and the fact that there
is no $\cal T$-cross in region $k$ it follows that $i=0$ and $v$ is not exposed.
We have $v \neq u_0$, for otherwise 
$V(P_0) \cap (V(Q_1) \cup V(Q_2))=\emptyset$ and
either $Q_1,Q_2,P_0$ are three crossed paths, or
$Q_1,Q_2,P_0,P_3,P_2$ is a separated double cross in $(G, \Omega)$. 
Since $v$ is not exposed we deduce that $Q$ satisfies the requirements of
Lemma~\ref{lem:nocross}. By the first part of
Lemma~\ref{lem:nocross} and the assumption made earlier it follows
that $i=0$ and $v$ is not exposed.
But the existence of $Q$ and the second statement of Lemma~\ref{lem:nocross}
imply that some leap of length $k$ has at least two exposed vertices,
a contradiction. (To see that let $j,p,q$ be as in Lemma~\ref{lem:nocross}. Replace $P_1$ by $Q_{3-j}$ and replace
$P_0$ by a suitable subpath of $Q_j \cup pP_0v_0\cup qQv$.)

We may therefore
assume that  $(G\backslash v,\uom\backslash v)$ has a
$\cal T$-tunnel $(Q_0,Q_1,Q_2)$ under $P_i$ for some $i\in\{1,2,\ldots,k\}$.
Then the leap $L'=(P_0,P_1,\ldots,P_{i-1},P_{i+1},\ldots,P_k)$
of length $k-1\ge3$ has a $\cal T'$-cross, where $\cal T'$ is the
corresponding orderly society, and the result follows in the same
way as above.~\qed

\begin{lemma}
\mylabel{lem:cycleC}
Assume Hypothesis~\ref{hyp:leap} and let $k \geq 3$. If there exist at least two
exposed vertices, then there exists a cycle $C$ and three disjoint
$(M\cup C)$-paths $R_1,R_2,R_3$ such that $R_i$ has ends $x_i\in V(C)$
and $y_i\in V(M)$, $C\backslash \{x_1,x_2,x_3\}$ is disjoint from $M$,
 $y_1=u_0$, $y_2=v_0$ and $y_3\in Z$.
\end{lemma}

\proof
Let $x_1$ be the closest exposed vertex  to $u_0$ on $P_0$,
and let $x_2$ be the closest exposed vertex to $v_0$.
Let $R_1=x_1P_0u_0$ and let $R_2=x_2P_0v_0$.
For $i=1,2$ let $S_i$ be an $(M\cup P_0)$-path with one end $x_i$
and the other end in $Z$.
By Lemma~\ref{lem:2pathstoZ} $S_1$ and $S_2$ intersect, and so we may
assume that $S_1\cap S_2$ is a path $R_3$ containing an end of both
$S_1$ and $S_2$, say $y_3$. Let $x_3$ be the other end of $R_3$.
Then $P_0\cup S_1\cup S_2$ includes a unique cycle $C$.
The cycle $C$ and paths $R_1,R_2,R_3$ are as desired for the lemma.~\qed

If the cycle $C$ in Lemma~\ref{lem:cycleC} can be chosen to
have at least four vertices,
then we say that the leap $(P_0,P_1,\ldots,P_k)$ is {\em diverse}.

\begin{lemma}
\label{lem:triangle}
\showlabel{lem:triangle}
Assume Hypothesis~\ref{hyp:leap}, let $k\ge4$, and let there be
no diverse leap of length $k$. If $C$ is as in Lemma~{\rm\ref{lem:cycleC}}
and $(G\backslash E(C),\uom)$ is rurally $4$-connected,
then $(G\backslash E(C),\uom)$ is rural.
\end{lemma}

\proof Since the leap $(P_0,P_1,\ldots,P_k)$ is not diverse,
it follows that $C$ is a triangle. Let $R_1,R_2,R_3$ and their ends
be numbered as in Lemma~\ref{lem:cycleC}.
We may assume that $P_0=R_1\cup R_2+x_1x_2$.
Since there is no diverse leap, Lemma~\ref{lem:2pathstoZ} implies that
there is no path in $G\backslash E(C)\backslash V(P_k)$ from $x_2$
to $v_k\Omega u_k$,
and none in $G\backslash E(C)\backslash V(P_1)$ from $x_1$
to $u_1\Omega v_1$. It also implies that no vertex on $P_0$ is exposed in $G \backslash x_1x_3 \backslash x_2x_3$.

As in Lemma~\ref{lem:1exposed}, we can apply Lemma~\ref{prestable} and assume,
by properly rerouting $M$ if necessary, that the conditions of Lemma~\ref{orderlyrural1} are satisfied.
We assume that
the society $(G\backslash E(C),\uom)$ has a
$\cal T$-jump, a  $\cal T$-cross, or a $\cal T$-tunnel, as otherwise by Lemma~\ref{orderlyrural1}  $(G\backslash E(C),\uom)$ is rural.
By the observation at the end of the previous paragraph this
$\cal T$-jump,  $\cal T$-cross, or  $\cal T$-tunnel
cannot use both $x_1$ and $x_2$;
say it does not use $x_2$.
But that contradicts
Lemma~\ref{lem:nojump} or the first part of Lemma~\ref{lem:nocross},
applied to $v=x_2$ and the graph $G \backslash x_1x_3$, in case of a $\cal T$-jump or a $\cal T$-cross.

Thus we may assume that $(G\backslash E(C)\backslash x_2,\Omega\backslash x_2)$
has a $\cal T$-tunnel $(Q_0,Q_1,Q_2)$
under $P_i$ for some $i\in\{1,2,\ldots,k\}$.
But then the leap $L'=(P_0,P_1,\ldots,P_{i-1},P_{i+1},\ldots,P_k)$
of length $k-1\ge3$ has a $\cal T'$-cross $(Q_1',Q_2')$,
where $\cal T'$ is the
corresponding orderly transaction, $Q_1'$ is obtained from $P_i$
by rerouting along $Q_0$ and $Q_2'$ is the  union of $Q_1\cup Q_2$
with the subpath of $P_i$ joining the ends of $Q_1$ and $Q_2$.
By the first half of Lemma~\ref{lem:nocross} applied to the graph $G \backslash x_1x_3$, the leap $L'$,
$v:=x_2$ and the $\cal T'$-cross $(Q_1',Q_2')$ we may assume that $i=1$ and that $y_3 \in v_2\Omega u_2 - \{u_0\}$.
By the second half of Lemma~\ref{lem:nocross} applied to the same entities
and $Q:=R_3+x_3x_2$ there exist $j\in\{1,2\}$,
$p\in V(Q_j'\cap R_1)$ and $q\in V(Q_j'\cap Q)$
such that $pP_0x_2$ and $qQx_2$ are internally disjoint from $Q_1'\cup Q_2'$.
If $j=1$, then $p,q$ belong to the interior of $Q_0$, and the
leap $(P_0,P_1,\ldots,P_k)$ is diverse, as a subpath of $Q_0$ joins a vertex of $R_1$ to a vertex of $Q$ in
$G \backslash x_1x_3$. If $j=2$ then we obtain a diverse leap from  $(P_0,P_1,\ldots,P_k)$
by replacing $P_1$ by $Q_1'$ and replacing $P_0$ by
a suitable subpath of $Q\cup v_0P_0p\cup Q_2'$.~\qed

\begin{lemma}
\mylabel{lem:q4}
Assume Hypothesis~\ref{hyp:leap}, let $k\ge3$,
let $(G,\Omega)$ be $4$-connected,
let $C,R_1,R_2,R_3$ be as in Lemma~{\rm\ref{lem:cycleC}},
and assume that $C$ is not a triangle.
Then there exist four disjoint $(M\cup C)$-paths, each with one end in $V(C)$
and the other end respectively in the sets $\{u_0\}$, $\{v_0\}$, $Z$
and $V(P_1\cup P_k)$.
\end{lemma}

\proof By an application of the proof of the max-flow min-cut theorem
there exist four disjoint $(M\cup C)$-paths, each with one end in $V(C)$
and the other end respectively in the sets $\{u_0\}$, $\{v_0\}$, $Z$
and $V(M)$.
By Lemma~\ref{lem:2pathstoZ} the fourth path does not end in
$V(M)-V(P_1)-V(P_k)$.
The result follows.~\qed

\begin{lemma}
\mylabel{lem:y4}
Assume Hypothesis~\ref{hyp:leap}, let $k\ge3$, 
let $C,R_1,R_2,R_3$ be as in Lemma~{\rm\ref{lem:cycleC}},
let $D:=M\cup C\cup R_1\cup R_2\cup R_3$,
and let $R_4$ be a $D$-path with ends
$x_4\in V(C)-\{x_1,x_2,x_3\}$ and $y_4\in V(P_1)$.
Then $x_1,x_2,x_3,x_4$ occur on $C$ in the order listed.
Furthermore, if $R$ is a $D$-path with ends
$x\in V(C)-\{x_1,x_2,x_3\}$ and $y \in V(M)$, then
$x_1,x_2,x_3,x$ occur on $C$ in the order listed and $y\in V(P_1)$.
\end{lemma}

\proof
The vertices $x_1,x_2,x_3,x_4$ occur on $C$ in the order listed by
Lemma~\ref{lem:2pathstoZ}. Now let $R$ be as stated.
By Lemma~\ref{lem:2pathstoZ} we have $y\in V(P_1\cup P_k)$,
and so by the first part of the lemma we may assume that $y\in V(P_k)$.
By the symmetric statement to the first half of the lemma
it follows that $x_1,x_2,x,x_3$ occur on $C$ in the order listed.
We may assume that $P_0$ is the unique path from $u_0$ to $v_0$ in
$R_1\cup R_2\cup C\backslash x_3$. Then $R_4\cup R\cup C\backslash V(P_0)$
includes a $\cal T$-jump disjoint from $P_0$,
contrary to Lemma~\ref{prejump}.~\qed

We need to further upgrade the assumptions of Hypothesis~\ref{hyp:leap},
as follows.

\begin{figure}
\begin{center}
\leavevmode
\includegraphics[scale = 1]{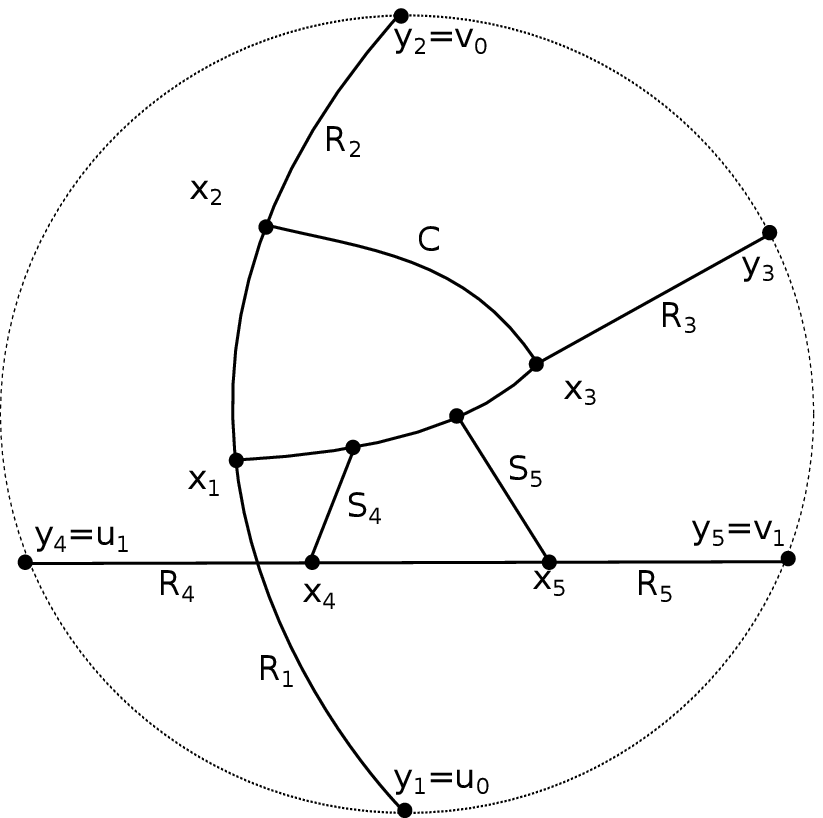}
\end{center}
\caption{Hypothesis~\ref{hyp:ray}.
}
\label{fig:hypothesis13}
\end{figure}
\showfiglabel{fig:hypothesis13}

\begin{hypothesis}
\mylabel{hyp:ray}
Assume  Hypothesis~\ref{hyp:leap}.
Let $C$ be a cycle with distinct vertices $x_1,x_2,x_3$
such that $C\backslash\{x_1,x_2,x_3\}$ is disjoint from $M$.
Let $R_1,R_2,R_3$ be pairwise disjoint $(M\cup C)$-paths such that
$R_i$ has ends $x_i$ and $y_i$, where $y_1=u_0$, $y_2=v_0$, and $y_3\in Z$.
By a {\em ray} we mean an $(M\cup C)$-path from $C$ to $M$,
disjoint from $R_1\cup R_2\cup R_3$.
We say that a vertex $v\in V(P_1)$ is {\em illuminated} if there
is a ray with end $v$.
Let $x_4,x_5\in V(P_1)$ be illuminated vertices such that either $x_4=x_5$, or
$u_1,x_4,x_5,v_1$ occur on $P_1$ in the order listed,
and $x_4P_1x_5$ includes all illuminated vertices.
Let $R_4:=u_1P_1x_4$ and $R_5:=v_1P_1x_5$, and let $y_4:=u_1$
and $y_5:=v_1$.
Let $S_4$ and $S_5$ be rays with ends $x_4$ and $x_5$, respectively,
and let
$A_0:=V(M)-V(P_1)$ and  $B_0:=V(C\cup S_4\cup S_5\cup x_4P_1x_5)$. (See Figure~\ref{fig:hypothesis13}.)
\end{hypothesis}

\begin{lemma}
\mylabel{lem:x4x5}
Assume Hypothesis~\ref{hyp:ray}, let $k\ge 3$, and let
$(G,\Omega)$ be $6$-connected.
Then $x_4\ne x_5$, and the path $x_4P_1x_5$ has at least one internal
vertex.
\end{lemma}

\proof If $x_4=x_5$ or $x_4P_1x_5$ has no internal vertex,
then by Lemma~\ref{lem:y4} the set $\{x_1,x_2,\ldots,x_5\}$ is a cutset
separating $C$ from $M\backslash V(P_1)$, contrary to the
$6$-connectivity of $(G,\Omega)$. Note that 
$V(C) - \{x_1,x_2,\ldots,x_5\}$ is non-empty
as it includes an end of a ray.~\qed

Assume Hypothesis~\ref{hyp:ray}. By Lemma~\ref{lem:x4x5} the
paths $R_1,R_2,\ldots,R_5$ are disjoint paths from $A_0$ to $B_0$.
The following lemma follows by a standard ``augmenting path" argument.

\begin{lemma}
\mylabel{lem:augmentation}
Assume Hypothesis~\ref{hyp:ray}, and let $k\ge 2$.
If there is no separation $(A,B)$ of order at most five with
$A_0\subseteq A$ and $B_0\subseteq B$, then there exist
an integer $n$ and internally
disjoint paths $Q_1,Q_2,\ldots,Q_n$ in $G$, where $Q_i$ has distinct
ends $a_i$ and $b_i$ such that
\myitem{(i)}$a_1\in A_0-\{y_1,y_2,\ldots,y_5\}$ and $
b_n\in B_0-\{x_1,x_2,\ldots,x_5\}$,
\myitem{(ii)}for all $i=1,2,\ldots,n-1$, $a_{i+1},b_{i}\in V(R_t)$
for some $t\in\{1,2,\ldots,5\}$, and $y_t,a_{i+1},b_i,x_t$ are
pairwise distinct and occur on $R_t$ in the order listed,
\myitem{(iii)}if $a_i,b_j\in V(R_t)$ for some $t\in\{1,2,\ldots,5\}$
and $i,j\in\{1,2,\ldots,5\}$ with $i>j+1$, then either $a_i=b_j$, or
$y_t,b_j,a_i,x_t$ occur on $R_t$ in the order listed, and
\myitem{(iv)}for $i=1,2,\ldots,n$, if a vertex of $Q_i$
belongs to $A_0\cup B_0\cup V(R_1\cup R_2\cup\cdots\cup R_5)$,
then it is an end of $Q_i$.
\end{lemma}

The sequence of paths $(Q_1,Q_2,\ldots,Q_n)$ as in Lemma~\ref{lem:augmentation}
will be called an {\em augmenting sequence}.

\begin{lemma}
\mylabel{lem:endinP2}
Assume Hypothesis~\ref{hyp:ray}, and let $k\ge 3$.
Then there is no augmenting sequence $(Q_1,Q_2,\ldots,Q_n)$, where
$Q_1$ is disjoint from $P_2$.
\end{lemma}

\proof Suppose for a contradiction that there is an augmenting
sequence $(Q_1,Q_2,\ldots,Q_n)$, where
$Q_1$ is disjoint from $P_2$, and let us assume that the leap
$(P_0,P_1,\ldots,P_k)$, cycle $C$, paths $R_1,R_2,R_3,S_4,S_5$ and
augmenting sequence $(Q_1,Q_2,\ldots,Q_n)$ are chosen with $n$
minimum. Let the ends of the paths $Q_i$ be labeled as in
Lemma~\ref{lem:augmentation}.
We may assume that
$P_0$ is the unique path from $u_0$ to $v_0$ in
$R_1\cup R_2\cup C\backslash x_3$.
We proceed in a series of claims.

\claim{1}{The vertex $b_n$ belongs to the interior of $x_4P_1x_5$.}

\noindent
To prove (1) suppose for a contradiction that $b_n\in V(C\cup S_4\cup S_5)$.
By Lemma~\ref{lem:y4}, the choice of $x_4,x_5$ and the fact that
$a_n\ne x_4,x_5$ by Lemma~\ref{lem:augmentation}(ii) we deduce
that $a_n\in V(R_i)$ for some $i\in\{1,2,3\}$. Then we can use $Q_n$ to modify
$C$ to include $a_nR_ix_i$ (and modify $R_1,R_2,R_3$ accordingly),
in which case $(Q_1,Q_2,\ldots,Q_{n-1})$ is an augmentation
contradicting the choice of $n$. This proves (1).

\claim{2}{$a_i,b_i\in V(R_j)$ for no $i\in\{1,2,\ldots,n\}$ and no
         $j\in\{1,2,\ldots,5\}$.}

\noindent
To prove (2) suppose to the contrary that $a_i,b_i\in V(R_j)$.
Then $1<i<n$ and
by rerouting $R_j$ along $Q_i$ we obtain an augmentation
$(Q_1,Q_2,\ldots,Q_{i-2},Q_{i-1}\cup b_{i-1}R_ja_{i+1}\cup Q_{i+1},Q_{i+2},\ldots,Q_n)$,
contrary to the minimality of $n$. This proves (2).

\claim{3}{$a_i,b_i\in V(R_1\cup R_2\cup R_3)$ for no $i\in\{1,2,\ldots,n\}$.}

\noindent
Using (2) the
proof of (3) is analogous to the argument at the end of the proof
of Claim (1).

\claim{4}{$a_i,b_i\in V(R_4\cup R_5)$ for no $i\in\{1,2,\ldots,n\}$.}

\noindent
By (2) one of $a_i,b_i$ belongs to
$R_4$ and the other to $R_5$. We can reroute $P_1$ along $Q_i$,
and then $(Q_1,Q_2,\ldots,Q_{i-1})$ becomes an augmentation,
contrary to the minimality of $n$.


\claim{5}{For $i=1,2,\ldots,n-1$, the graph $Q_i\cup R_1\cup R_2\cup R_3$
includes no $\cal T$-jump.}

\noindent
This claim follows from (3),
Lemma~\ref{lem:2pathstoZ}  and Lemma~\ref{prejump} applied to $P_0$.

\claim{6}{$a_1\not\in v_1\Omega u_1$.}

\noindent
To prove (6) suppose for a contradiction that $a_1\in v_1\Omega u_1$.
Since $a_1\ne y_1$, we may assume from the symmetry that
$a_1\in v_1\Omega y_1-\{y_1\}$.
Then $b_1\in V(P_1\cup R_1)$ by (5). But if $b_1\in V(R_i)$, where $i=1$
or $i=5$,
then by rerouting $R_i$  along $Q_1$
we obtain an augmenting sequence $(Q_2\cup x_1R_ia_2,Q_3,Q_4,\ldots,Q_n)$,
contrary to the choice of $n$.
Thus $b_1\in u_1P_1x_5$. By replacing $P_1$ by the path
$Q_1\cup u_1P_1b_1$ and considering the paths $R_3$ and $S_5\cup R_5$
we obtain contradiction to Lemma~\ref{lem:2pathstoZ}.
This proves (6).

\claim{7}{$a_1\not\in u_k\Omega v_k$.}

Similarly as in the proof of (6), if $a_1\in u_k\Omega v_k$,
then $b_1\in V(R_2)$ by (5),
and we reroute $R_2$ along $Q_1$ to obtain a contradiction to the minimality
of $n$. This proves (7).

\claim{8}{$a_1\in V(P_k)$.}

To prove (8) we may assume by (6) and (7) that $a_1\in Z$.
Then $b_1\in V(R_3\cup P_1)$ by (5).
If $b_1\in V(R_3)$, then
we reroute $R_3$ along $Q_1$ as before.
Thus $b_1\in V(P_1)$. It follows from (5) and the hypothesis $V(P_2) \cap V(Q_1)= \emptyset$ that $a_1\in u_1\Omega u_2-\{u_1,u_2\}$
or $a_1\in v_2\Omega v_1-\{v_1,v_2\}$, and so from the symmetry we
may assume the latter.

Let us assume for a moment that  $y_3 \in a_1\Omega v_1$.
We reroute $P_1$ along $Q_1 \cup b_1P_1v_1$. 
The union of $R_3$, $R_2$ and a path in $C$ between $x_2$ and $x_3$, 
avoiding $x_1,x_4,x_5$, will play the role of $P_0$ after rerouting. 
If $b_1 \in x_4P_1v_1 - \{x_4\}$, then
$R_1\cup C\cup S_4 \cup R_4$ includes two disjoint paths that contradict 
Lemma~\ref{lem:2pathstoZ} applied to the new frame and new path $P_0$. 
Therefore $b_1 \in V(R_4)$, and hence $(u_1P_1a_2 \cup Q_2, Q_3, \ldots, Q_n)$ 
is an augmenting sequence after the
rerouting, contrary to the choice of $n$.

It follows that $y_3 \not \in a_1\Omega v_1$.
If $b_1\in V(R_5)$, we replace $P_1$ by $Q_1\cup u_1P_1b_1$;
then $(v_1P_1a_2\cup Q_2,Q_3,\ldots,Q_n)$ is an augmenting sequence that
contradicts the choice of $n$.
So it follows  that $b_1\in u_1P_1x_5$.
But now $(G,\Omega)$ has a gridlet using the paths $P_0, P_k$,
$Q_1\cup u_1P_1b_1$ and a subpath of $R_5\cup S_5\cup R_3\cup C\backslash V(P_0)$.
This proves (8).

\claim{9}{$n>1$.}

\noindent
To prove (9) suppose for a contradiction that $n=1$. 
Thus $b_1$ belongs to the interior of $x_4Px_5$ by (1), and
$a_1\in V(P_k)$ by (8). But then $Q_1$ is a $\cal T$-jump, contrary to (5).

\claim{10}{$b_1\in V(R_3)$.}

\noindent
To prove (10)  we first notice that $b_1\in V(R_2\cup R_3)$ by (5), (9) and (1).
Suppose for a contradiction that $b_1\in V(R_2)$.
Then $a_2\in V(R_2)$, but $b_2\not\in V(R_1\cup R_2\cup R_3)$ by (3)
and $b_2\not\in V(P_1)$ by (5), a contradiction. This proves (10).

\medskip
Let $P_{12}$ and $P_{34}$ be two disjoint subpaths of $C$, where the
first has ends $x_1,x_2$, and the second has ends $x_3,x_4$.
By (8) and (10) the path $Q_1 \cup b_1R_3x_3 \cup P_{34} \cup
S_4$ is a $\cal T$-jump disjoint from $R_1 \cup P_{12} \cup R_2$,
contrary to Lemma~\ref{prejump}.~\qed

We are now ready to prove Theorem~\ref{thm:leap}.
\medskip

\noindent {\bf Proof of Theorem~\ref{thm:leap}.}
Let $(G,\Omega)$ be a $6$-connected society with a leap of length five.
Thus we may assume that Hypothesis~\ref{hyp:leap} holds for $k=5$.
By Lemma~\ref{lem:1exposed} either $(G,\Omega)$ is nearly rural,
in which case the theorem holds,
or there exists a leap of length at least four with at least two exposed vertices.
Thus we may assume that there exists a leap of length four with at least two exposed vertices.
Let $C$ be a cycle as in Lemma~\ref{lem:cycleC}. If there is no
diverse leap, then $C$ is a triangle, $(G\backslash E(C),\Omega)$
is rurally $4$-connected and hence rural by Lemma~\ref{lem:triangle}, and the theorem holds.
Thus we may assume that the cycle $C$ is not a triangle,
and so by Lemma~\ref{lem:q4} we may assume that Hypothesis~\ref{hyp:ray} for $k=4$
holds.
By Lemma~\ref{lem:x4x5} and the $6$-connectivity of $G$ there is
no separation $(A,B)$ as described in Lemma~\ref{lem:augmentation},
and hence by that lemma there exists an augmenting sequence
$(Q_1,Q_2,\ldots,Q_n)$.
By Lemma~\ref{lem:endinP2} the path $Q_1$ intersects $P_2$,
and hence $Q_1$ is disjoint from $P_3$, contrary to
Lemma~\ref{lem:endinP2} applied to the leap $(P_0,P_1,P_3,P_4)$
of length three and an augmenting sequence $(Q_1',Q_2,\ldots, Q_n)$, where $Q_1'$ is the union of $Q_1$ and $a_1P_2u_2$ or $a_1P_2v_2$.~\qed 


\section{Societies of bounded depth}
\label{sec:transactions}
\showlabel{sec:transactions}

Let $(G,\uom)$ be a society.  A {\em linear decomposition} of $(G,\uom)$
is an enumeration $\{t_1,\dots,t_n\}$ of $V(\uom)$ where
$(t_1,\dots, t_n)$ is clockwise, together with a family $(X_i:1\le i\le n)$
of subsets 
 of $V(G)$, with the following properties:
\myitem{(i)} $\bigcup (G[X_i]:1\le i\le n)=G$,
\myitem{(ii)} for $1\le i\le n$, $t_i\in X_i$, and
\myitem{(iii)} for $1\le i\le i'\le i''\le n$, $X_i\cap X_{i''}\subseteq X_{i'}$.

\noindent The {\em depth} of such a linear decomposition is
$$\max(|X_i\cap X_{i'}|:1\le i <i'\le n),$$
and the depth of $(G,\uom)$ is the minimum depth of a linear decomposition
of $(G,\uom)$.
Theorems~(6.1), (7.1) and (8.1) of~\cite{RobSeyGM9} imply the following.

\begin{theorem}\label{gm9}
\showlabel{gm9}
There exists an integer $d$ such that every $4$-connected society
$(G,\uom)$ either has a separated doublecross, 
three crossed paths or a leap of length
five, or some planar truncation of $(G,\uom)$ has depth at most $d$.
\end{theorem}

In light of Theorems~\ref{thm:leap} and~\ref{gm9},
in the remainder of the paper we concentrate on societies of
bounded depth.  We need a few definitions.  
Let $(G,\Omega)$ be a society,  let 
$u_1,u_2,\ldots,u_{4t}$ be clockwise in $\Omega$, and let
$P_1,P_2,\ldots,P_{2t}$ be disjoint bumps in $G$ such that
for $i=1,2,\ldots,2t$ 
the path $P_{2i-1}$ has ends $u_{4i-3}$ and $u_{4i-1}$, and
the path $P_{2i}$ has ends $u_{4i-2}$ and $u_{4i}$.
In those circumstances we say that  $(G,\Omega)$ has 
{\em $t$ disjoint consecutive crosses}.

Now let $u_1,v_1,w_1,u_2,v_2,w_2,\ldots,u_t,v_t,w_t$ be clockwise
in $\Omega$, let 
$x\in V(G)-\{u_1,v_1,\allowbreak w_1,\ldots,u_t,\allowbreak v_t,w_t\}$,
for $i=1,2,\ldots,t$ let $P_i$ be a path in
$G\backslash x$ with ends $u_i$ and $w_i$
and otherwise disjoint from $V(\Omega)$, let $Q_i$ be a path with
ends $x$ and $v_i$ and  otherwise disjoint from $V(\Omega)$, 
and assume that the paths $P_i$ and $Q_i$ are
pairwise disjoint, except that the paths $Q_i$ meet at $x$.
Let $W$ be the union of all the paths $P_i$ and $Q_i$.
We say that $W$ is 
a {\em windmill with $t$ vanes}, and that the graph $P_i\cup Q_i$ is
a {\em vane} of the windmill.

Finally, let  $u_1,u_2,\ldots,u_t$ and  $v_1,v_2,\ldots,v_t$ be
vertices of $V(\Omega)$ such that for all $x_i\in\{u_i,v_i\}$
the sequence  $x_1,x_2,\ldots,x_t$ is clockwise in $\Omega$.
Let $z_1,z_2\in V(G)-\{u_1,v_1,\ldots,u_t,v_t\}$ be distinct,
for $i=1,2,\ldots,t$ let $P_i$ be a path in $G\backslash z_2$ with ends $z_1$
and $u_i$ and otherwise disjoint from $V(\Omega)$, and let 
$Q_i$ be a path in $G\backslash z_1$ with ends $z_2$
and $v_i$ and otherwise disjoint from $V(\Omega)$. 
Assume that the paths $P_i$ and $Q_j$ are disjoint, except that the
$P_i$ share $z_1$, the $Q_i$ share $z_2$ and $P_i$ and $Q_i$ are allowed to
intersect.
Let $F$ be the union of all  the paths $P_i$ and $Q_i$.
Then we say that $F$ is a {\em fan with $t$ blades}, and we
say that $P_i\cup Q_i$ is a {\em blade} of the fan.
The vertices $z_1$ and $z_2$ will be called the {\em hubs} of the fan.
In Section~\ref{usewars} we prove the following theorem.

\begin{theorem}
\label{thm1}
\showlabel{thm1}
For every two integers $d$ and $t$ there exists an integer $k$ such
that every $6$-connected%
\REMARK{We could lower 6 to 5, but the proof will have to be changed}
$k$-cosmopolitan society  $(G,\Omega)$ of depth at most $d$
contains one of the following:
\myitem{(1)}$t$ disjoint consecutive crosses, or
\myitem{(2)}a windmill with $t$ vanes, or
\myitem{(3)}a fan with $t$ blades.
\end{theorem}

Unfortunately, windmills and fans are nearly rural, and so for our 
application we need to improve Theorem~\ref{thm1}.  We need more definitions.

Let $x,u_i,v_i,w_i,P_i,Q_i$ be as in the definition of a windmill $W$
with $t$ vanes,
let $a,b,c,d\in V(G)$ be such that 
$u_1,v_1,w_1,\ldots,u_t,v_t,w_t, a, b, c, d$ is clockwise
in $\Omega$, and let $(P,Q)$ be a cross disjoint from $W$ whose 
paths have ends in $\{a,b,c,d\}$. In those circumstances we say that
$W\cup P\cup Q$ is a {\em windmill with $t$ vanes and a cross}.

Now let $u_i,v_i,P_i,Q_i$ be as in the definition of a fan $F$ with
$t$ blades, and let $a,b,c,d\in V(\Omega)$ be such that 
all $x_i\in\{u_i,v_i\}$
the sequence  $x_1,x_2,\ldots,x_t,a,b,c,d$ is clockwise in $\Omega$.
Let  $(P,Q)$ be a cross disjoint from $F$ whose
paths have ends in $\{a,b,c,d\}$. In those circumstances we say that
$W\cup P\cup Q$ is a {\em fan with $t$ blades and a cross}.

Let $z_1,z_2,u_i,v_i,P_i,Q_i$ be as in the definition of a fan $F$ with
$t$ blades, and let $a_1,b_1,c_1,a_2,\allowbreak b_2,c_2\in V(G)$ be such that
all $x_i\in\{u_i,v_i\}$
the sequence  $x_1,x_2,\ldots,x_t,a_1,b_1,c_1,\allowbreak a_2,b_2,c_2$ 
is clockwise in $\Omega$, except that we permit $c_1=a_2$.
 For $i=1,2$ let $L_i$ be a path in $G\backslash V(F)$
with ends $a_i$ and $c_i$ and otherwise disjoint from $V(\Omega)$,
and let $S_i$ be a path with ends $z_i$ and $b_i$ and otherwise disjoint
from $V(F)\cup V(\Omega)$.
If the paths $L_1,L_2,S_1,S_2$ are pairwise disjoint, except possibly
for $L_1$ intersecting $L_2$ at $c_1=a_2$, then we say that 
$F\cup L_1\cup L_2\cup S_1\cup S_2$ is a {\em fan with $t$ blades and two jumps}.

Now let  $u_i,v_i,P_i,Q_i$ be as in the definition of a fan $F$ with
$t+1$ blades, and let $a,b\in V(\Omega)$ be such that
all $x_i\in\{u_i,v_i\}$
the sequence  $x_1,x_2,\ldots,x_t,a,x_{t+1},b$ is clockwise in $\Omega$.
Let $P$ be a path in $G\backslash V(F)$ with ends $a$ and $b$, and
otherwise disjoint from $V(F)$. We say that $F\cup P$ is a
{\em fan with $t$ blades and a jump}.
In Section~\ref{sec:lack} we improve Theorem~\ref{thm1} as follows.

\begin{theorem}
\label{thm2}
\showlabel{thm2}
For every two integers $d$ and $t$ there exists an integer $k$ such
that every $6$-connected
$k$-cosmopolitan society  $(G,\Omega)$ of depth at most $d$
is either nearly rural, or
contains one of the following:
\myitem{(1)}$t$ disjoint consecutive crosses, or
\myitem{(2)}a windmill with $t$ vanes and a cross, or
\myitem{(3)}a fan with $t$ blades and a cross, or
\myitem{(4)}a fan with $t$ blades and a jump, or
\myitem{(5)}a fan with $t$ blades and two jumps.
\end{theorem}

\noindent
For $t=4$ each of the above outcomes gives a turtle, and hence we
have the following immediate corollary.

\begin{corollary}
\label{thm3}
\showlabel{thm3}
For every integer $d$  there exists an integer $k$ such
that every $6$-connected
$k$-cosmopolitan society  $(G,\Omega)$ of depth at most $d$
is either nearly rural, or has a turtle.
\end{corollary}

The next four sections are devoted to proofs of Theorems~\ref{thm1}
and~\ref{thm2}. The proof of Theorem~\ref{thm1} will be completed in
Section~\ref{sec:wars} and the proof of  Theorem~\ref{thm2} will be 
completed in Section~\ref{sec:lack}.
At that time we will be able to deduce Theorem~\ref{thm:society}.

\section{Crosses and goose bumps}

In this section we prove that a society $(G,\Omega)$ either 
satisfies Theorem~\ref{thm1}, or it has 
many disjoint bumps.
If $X$ is a set and $\Omega$ is a cyclic permutation,
we define $\Omega\backslash X$ to be $\Omega\restriction(V(\Omega)-X)$.
Let $P_1,P_2,\ldots,P_k$ be a set of pairwise disjoint bumps in
$(G,\Omega)$, where $P_i$ has ends $u_i$ and $v_i$ and
$u_1,v_1,u_2,v_2,\ldots,u_k,v_k$ is clockwise in $\Omega$. 
In those circumstances we say that $P_1,P_2,\ldots,P_k$ is a 
{\em goose bump in $(G,\Omega)$ of strength $k$.}

\begin{lemma}
\label{EPbumps}
\showlabel{EPbumps}
Let $b,d$ and $t$ be positive integers, and let $(G,\Omega)$ be
a society of depth at most $d$. Then either $(G,\Omega)$ has a goose
bump of strength $b$, or there is a set $X\subseteq V(G)$
of size at most $(b-1)d$
such that the society $(G\backslash X, \Omega\backslash X)$ has no bump.
\end{lemma}

\proof
Let $(t_1,t_2,\dots,t_n)$ and $(X_1,X_2,\dots, X_n)$ be a linear
decomposition of $(G,\uom)$ of depth at most $d$, and for $i=1,2,\dots, n-1$
let $Y_i=X_i\cap X_{i+1}$.  If $P$ is a bump in $(G,\uom)$, then the
axioms of a linear decomposition imply that
$$I_P:=\{i\in \{1,2,\dots, n-1\}: Y_i\cap V(P)\ne \emptyset\}$$
is a nonempty subinterval of $\{1,2,\dots, n-1\}$.  It follows that
either there exist bumps $P_1,P_2,\dots, P_b$ such that $I_{P_1},
I_{P_2},\dots, I_{P_b}$ are pairwise disjoint, or there exists a set
$I\subseteq \{1,2,\dots, n-1\}$ of size at most $b-1$ such that 
$I\cap I_P\ne\emptyset$ for every bump $P$.  In the former case
$P_1,P_2,\dots, P_b$ is a desired goose bump, and in the latter 
case the set $X:=\bigcup_{i\in I} Y_i$ is as desired.~\qed

The proof of the following lemma is similar and is omitted.

\begin{lemma}
\label{EPcrosses}
\showlabel{EPcrosses}
Let $t$ and $d$ be positive integers, and let $(G,\Omega)$ be
a society of depth at most $d$. Then either $(G,\Omega)$ has 
$t$ disjoint consecutive crosses, or there is a set $X\subseteq V(G)$
of size at most $(t-1)d$
such that the society $(G\backslash X, \Omega\backslash X)$ is cross-free.
\end{lemma}

\begin{lemma}
\label{existgoose}
\showlabel{existgoose}
Let $d,b,t$ be positive integers, let $k\ge (b-1)d+(t-1)\binom{(b-1)d}2+1$
and let $(G,\Omega)$ be a  $3$-connected 
society of depth at most $d$ such that
at least $k$ vertices in $V(\uom)$ have at least two neighbors in 
$V(G)$. Then $(G,\Omega)$ 
has either a fan with $t$ blades, or 
a goose bump of strength $b$.
\end{lemma}

\proof
By Lemma~\ref{EPbumps} we may assume that there exists a set
$X\subseteq V(G)$ of size at most $(b-1)d$ such that
$(G\backslash X,\uom\backslash X)$ has no bump. There are at least
$(t-1)\binom{(b-1)d}2 +1$ vertices in $V(\uom)-X$ with at least two
neighbors in $V(G)$.  Let $v$ be one such vertex, and let $H$
be the component of $G\backslash X$ containing $v$.  Since
$(G\backslash X,\uom\backslash X)$ has no bumps it follows that
$V(H)\cap V(\uom)=\{v\}$.  By the fact that $v$ has at least two
neighbors in $G$ (if $V(H)=\{v\}$) or the 3-connectivity of
$(G,\uom)$ (if $V(H)\ne \{v\}$) it follows that $H$ has at least two
neighbors in $X$. Thus there exist distinct vertices $z_1,z_2$ such that
for at least $t$ vertices of $v\in V(\uom)-X$ the component of 
$G\backslash X$ containing $v$ has $z_1$ and $z_2$ as neighbors.
It follows that $(G,\uom)$ has a fan with $t$ blades, as desired.~\qed

\section{Intrusions, invasions and wars}

Let $\Omega$ be a cyclic permutation. A {\em base} in $\Omega$ is
a pair $(X,Y)$ of subsets of $V(\Omega)$ such that $|X\cap Y|=2$,
$X\cup Y=V(\Omega)$ and for distinct elements $x_1,x_2\in X$ and $y_1,y_2\in Y$
the sequence $(x_1,y_1,x_2,y_2)$ is not clockwise.
Now let  $(G,\Omega)$ be a society. A separation $(A,B)$ of $G$ is
called an {\em intrusion} in  $(G,\Omega)$ if there exists a base
$(X,Y)$ in $\Omega$ such that $X\subseteq A$, $Y\subseteq B$ and there
exist disjoint paths $(P_v)_{v\in A\cap B}$, each with one end in
$X$, the other end in $Y$ and with $v\in V(P_v)$. 
The intrusion $(A,B)$ is {\em \minimal} if there is no intrusion
$(A',B')$ of order $|A\cap B|$ with base $(X,Y)$ such that $A'$ is
a proper subset of $A$.
The paths $P_v$ will be called {\em longitudes}
for the intrusion $(A,B)$. We say that $(A,B)$ is {\em based} at
$(X,Y)$, and that $(X,Y)$ is a base for $(A,B)$.
An intrusion $(A,B)$ in  $(G,\Omega)$ is an {\em invasion} 
if $|A\cap B\cap V(\Omega)|=2$.

\begin{lemma}
\label{existintrusion}
\showlabel{existintrusion}
Let $d$ be a positive integer, and let $(G,\Omega)$ be a
society of depth at most $d-1$. Then for every base $(X,Y)$ in $\Omega$
there exists an intrusion of order at most $2d$ based at $(X,Y)$.
\end{lemma}
                                                                                   
\proof
Let $(t_1,t_2,\dots, t_n)$ and $(X_1,X_2,\dots, X_n)$ be a linear
decomposition of $(G,\uom)$ of depth at most $d-1$, and let
$X\cap Y=\{t_i,t_j\}$.  Let $i',j'\in \{1,2,\dots, n\}$ be such that
$|i-i'| = |j-j'|=1$, and let $Z:=(X_i\cap X_{i'})\cup
(X_j\cap X_{j'})\cup \{t_i,t_j\}$.  It follows from the axioms of a linear
decomposition that $|Z|\le 2d$ and that $Z$ separates $X$ from $Y$ in $G$.
Thus there exists a separation $(A,B)$ of $G$ of order at most $2d$ with
$X\subseteq A$ and $Y\subseteq B$.  Any such separation $(A,B)$ with
$|A\cap B|$ minimum is as desired by Menger's theorem.~\qed

An intrusion $(A,B)$  in a society $(G,\Omega)$ is {\em $t$-separating} 
if $(G,\Omega)$ has goose bumps $P_1,P_2,\ldots,P_t$ and
$Q_1,Q_2,\ldots,Q_t$ such that $V(P_i)\subseteq A-B$ and
$V(Q_i)\subseteq B-A$ for all $i=1,2,\ldots,t$.

\begin{lemma}
\label{existtsepintrus}
\showlabel{existtsepintrus}
Let $d,s,t$ be  positive integers, and let $(G,\Omega)$ be a
society of depth at most $d-1$ with a goose bump of strength $t(s+2d)$.
Then there exist $s$-separating {\minimal} intrusions 
$(A_1,B_1),(A_2,B_2),\ldots,(A_t,B_t)$ of order at most $2d$ such that
$A_i\cap A_j\subseteq B_i\cap B_j$ for all pairs of distinct indices
$i,j=1,2,\ldots,t$.
\end{lemma}
                                                                                   
\proof
Let ${\cal P}$ be the set of paths comprising a goose bump of strength
$t(s +2d)$.  Thus there exist bases $(X_1,Y_1),(X_2,Y_2),\dots
(X_t,Y_t)$ such that the sets $X_i$ are pairwise disjoint and for each
$i=1,2,\dots, t$ exactly $s+2d$ of the paths in ${\cal P}$ have
both ends in $X_i$.  By Lemma~\ref{existintrusion} there exists, 
for each $i=1,2,\dots, t$, an intrusion $(A_i,B_i)$ of order at 
most $2d$ based at $(X_i,Y_i)$.

Let us choose, for each $i=1,2,\dots, t$, an intrusion $(A_i,B_i)$ of order
at most $2d$ based at $(X_i,Y_i)$ in such a way that
\begin{equation}\label{eq:min}
\sum^t_{i=1}|A_i|\mbox{ is minimum.}\end{equation}
We claim that $A_i\cap A_j\subseteq B_i\cap B_j$.  To prove the
claim suppose to the contrary that say $x\in A_1\cap A_2-B_1\cap B_2$.  Let
\begin{align*}
A'_1&=A_1\cap B_2,\\
B'_1&=A_2\cup B_1,\\
A'_2&=A_2\cap B_1,\\
B'_2&=A_1\cup B_2.\end{align*}
Then $(A'_1,B'_1)$ and $(A'_2,B'_2)$ are separations of $G$ with
$X_1\subseteq A'_1$, $Y_1\subseteq B'_1$, $X_2\subseteq A'_2$ and
$Y_2\subseteq B'_2$.  We have
$$|A_1\cap B_1|+|A_2\cap B_2| = |A'_1\cap B'_1| + |A'_2\cap B'_2|.$$
Furthermore, since each longitude for $(A_1,B_1)$ intersects
$A'_1\cap B'_1$ we deduce that $|A'_1\cap B'_1| \ge |A_1\cap B_1|$,
and similarly $|A'_2\cap B'_2|\ge |A_2\cap B_2|$.  Thus the last two 
inequalities hold with equality, and hence the longitudes for $(A_1,B_1)$
are also longitudes for $(A'_1,B'_1)$, and the longitudes for $(A_2,B_2)$
are longitudes for $(A'_2,B'_2)$.  It follows that for $i=1,2$ the
separation $(A'_i,B'_i)$ is an intrusion in $(G,\uom)$ based at
$(X_i,Y_i)$ of order $|A_i\cap B_i|$.  Since $A_1\cap A_2-(B_1\cap B_2)=
(A_1\cap A_2-B_1)\cup (A_1\cap A_2-B_2)$ we may assume that $x\in A_1-B_2$.
But then replacing $(A_1,B_1)$ by $(A'_1,B'_1)$ produces a set of intrusions
that contradict \eqref{eq:min}.  This proves our claim that
$A_i\cap A_j\subseteq B_i\cap B_j$ for all distinct integers 
$i,j=1,2,\dots, t$.

Since at most $2d$ of the paths in ${\cal P}$ with ends in $X_i$ can intersect
$A_i\cap B_i$, we deduce that each intrusion $(A_i,B_i)$ is $s$-separating.
Moreover, each $(A_i,B_i)$ is clearly minimal by \eqref{eq:min}.~\qed

We need a lemma about subsets of a set.

\begin{lemma}
\label{subsets}
\showlabel{subsets} 
Let $d$ and $t$ be positive integers, and let $\cal F$ be a family
of $2^{d+1\choose 2}t^d$ distinct subsets of a set $S$,
where each member of ${\cal F}$ has size at most $d$. Then there exist
a set $X\subset S$ of size at most ${d+1\choose 2}$ and a family
${\cal F}'\subseteq\cal F$ of size at least $t$ such that 
$F\cap F'\subseteq X$ for every two distinct sets $F,F'\in{\cal F}'$.
\end{lemma}

\proof
We proceed by induction on $d+t$.  If $d=1$ or $t=1$, then the lemma clearly
holds, and so we may assume that $d,t>1$.  Let $F_0\in{\cal F}$ be minimal
with respect to inclusion.  If ${\cal F}$ has a subfamily ${\cal F}_1$ of
at least $2^{d+1\choose 2}(t-1)^d$ sets disjoint from $F_0$, then the
result follows from the induction hypothesis applied to ${\cal F}_1$ and
by adding $F_0$ to the family thus obtained.  If the family 
${\cal F}_2=\{F-F_0: F\in {\cal F}, F\cap F_0\ne\emptyset\}$ includes at
least $2^{\binom{d}2} t^{d-1}$ distinct sets, then the result follows
from the induction hypothesis applied to ${\cal F}_2$ by adding $F_0$ to
the set thus obtained.  Thus we may assume neither of the two cases holds.
Thus
$$|{\cal F}|\le 2^{\binom {d+1}2} (t-1)^d-1+2^d2^{\binom{d}2}t^{d-1}-1+1
<2^{\binom {d+1}2}t^d,$$
a contradiction.~\qed

\begin{lemma}
\label{existdisjintrus}
\showlabel{existdisjintrus}
Let $d,s,t$ be positive integers, and let $(G,\Omega)$ be a
society of depth at most $d-1$ with a goose bump of strength 
$2^{{2d+1}\choose 2}t^{2d}(s+2d)$.
Then there exist a set $X\subseteq V(G)$ of size at most ${2d+1\choose 2}$
and $s$-separating intrusions
$(A_1,B_1),(A_2,B_2),\ldots,(A_t,B_t)$ in 
$(G\backslash X,\Omega\backslash X)$ such that
$A_i\cap A_j=\emptyset$ for all pairs of distinct indices
$i,j=1,2,\ldots,t$.
\end{lemma}
                                                                                   
\proof
Let $T=2^{{2d+1}\choose 2}t^{2d}$. By Lemma~\ref{existtsepintrus} there
exist $s$-separating {\minimal} intrusions 
$(A_1,B_1),\allowbreak (A_2,B_2),\ldots,(A_T,B_T)$ of order at most $2d$ such that
$A_i\cap A_j\subseteq B_i\cap B_j$ for all pairs of distinct indices
$i,j=1,2,\ldots,t$. By Lemma~\ref{subsets} applied to the sets 
$A_i\cap B_i$ there exist a set $X\subseteq \bigcup_{i=1}^T (A_i\cap B_i)$ 
of size at most
${{2d+1}\choose2}$ and a subset of $t$ of those intrusions, say
$(A_1,B_1),(A_2,B_2),\ldots,(A_t,B_t)$, such that 
$A_i\cap B_i\cap A_j\cap B_j\subseteq X$ for all distinct integers
$i,j=1,2,\ldots,t$. It follows that $(A_i-X,B_i-X)$ are as required for
$(G\backslash X,\Omega\backslash X)$.
\qed

Our next objective is to prove, albeit with weaker bounds, that the conclusion of
Lemma~\ref{existdisjintrus} can be strengthened to assert that the
intrusions $(A_i,B_i)$ therein are actually invasions.

Let $(A,B)$ be an intrusion in a society $(G,\Omega)$ based at $(X,Y)$. 
A path $P$ in $G[A]$ is a {\em meridian} for $(A,B)$ if its ends are
the two vertices of $X\cap Y$.
If $P$ is a meridian for $(A,B)$ and $(L_v)_{v\in A\cap B}$ are longitudes
for $(A,B)$, then the graph $(P\cup\bigcup_{v\in A\cap B} L_v)\backslash(B-A)$
is called a {\em frame} for $(A,B)$.

\begin{lemma}
\label{existtsepinv}
\showlabel{existtsepinv}
Let $\lambda$ and $s$ be positive integers, let $s'=(s-1)(\lambda-1)+1$, 
let $(G,\Omega)$ be a cross-free society,
and let $(A,B)$ be an $s'$-separating {\minimal} intrusion in $(G,\Omega)$ 
of order at most $\lambda$. Then there exists an 
$s$-separating {\minimal} invasion $(C,D)$
in $(G,\Omega)$ of order at most $\lambda$ with a frame $F$ such that 
$V(F)-V(\Omega)\subseteq A$.
\end{lemma}

\proof
We may assume that

\claim{1}{%
  there is no integer $\lam'\le \lam$ and an 
  $((s-1)(\lam'-1)+1)$-separating minimal intrusion 
  $(A',B')$ in $(G,\uom)$ of order
  at most $\lam'$ with $A'$ a proper subset of $A$,
}

\noindent for if $(A',B')$ exists, and it satisfies the conclusion of
the lemma, then so does $(A,B)$.  We first show that $(A,B)$ has a meridian.
Indeed, suppose not.  Let $(X,Y)$ be a base of $(A,B)$ and let 
$X\cap Y=\{u,v\}$; then $G[A]$ has no $u$-$v$ path.  Since $(G,\uom)$ is
cross-free it follows that $G[A]$ has a separation $(A_1,A_2)$ of order
zero such that both $X_1=X\cap A_1$ and $X_2=X\cap A_2$ are 
intervals in $\uom$. It follows that there exist $Y_1,Y_2$ such that
$(X_1,Y_1)$ and $(X_2,Y_2)$ are bases.  Thus $(A_1, A_2\cup B\cup
(X_1\cap Y_1))$ and $(A_2, A_1\cup B\cup (X_2\cap Y_2))$ are 
minimal 
intrusions, 
and one of them violates (1).  This proves
that $(A,B)$ has a meridian.

Let $M$ be a meridian in $(A,B)$, let $(L_v)_{v\in A\cap B}$ be a
collection of longitudes for $(A,B)$ and let 
$F=M\cup \bigcup_{v\in A\cap B} (L_v\backslash(B-A))$.  
By the same argument that justifies (1) we may assume that

\myclaim{2} there is no integer $\lam' <\lam$ and an 
$((s-1)(\lam'-1)+1)$-separating minimal intrusion $(A',B')$ in $(G,\uom)$
of order at most $\lam'$ with frame $F'$ such that $F'\backslash 
V(\uom)$ is a subgraph of $F$.\par
\medskip

We claim that $|A\cap B\cap V(\uom)|=2$.  We first prove that
$A\cap B\cap X=\{u,v\}$.  To this end suppose for a contradiction that
$w\in A\cap B\cap X-\{u,v\}$; then $w$ divides $X$ into 
two cyclic intervals $X_1$ and $X_2$ with ends $u,w$ and $w,v$,
respectively.  Let $Y_1$ and $Y_2$ be the complementary cyclic intervals
so that $(X_1,Y_1)$ and $(X_2,Y_2)$ are bases.  

For $i=1,2$ let $A_i$ consist of $w$ and all vertices $a\in A$ such that
there exists a path in $G[A]\backslash w$ with one end $a$ and the
other end in $X_i-\{w\}$, and let $A_3=A-A_1-A_2$. It follows
that $A_1\cap A_2=\{w\}$, for if $P$ is a path in $G[A]\backslash w$ with
one end in $X_1$ and the other end in $X_2$, then $(P,P_w)$ is a cross
in $(G,\uom)$, a contradiction. Thus $(A_1,A_2\cup A_3\cup B)$ and 
$(A_2, A_1\cup A_3\cup B)$ are minimal intrusions based on $(X_1,Y_1)$
and $(X_2,Y_2)$, respectively, with $A_1,A_2\subseteq A$.  Thus one of 
them violates (2).

Next we show that $|A\cap B\cap Y|=2$, and so we suppose for a
contradiction that there exists $z\in A\cap B\cap Y-\{u,v\}$.  We define
$B_1,B_2,B_3, X_1,Y_1,X_2,Y_2$ analogously as in the previous 
paragraph, but with the roles of $A$ and $B$ reversed.  Similarly
we find that one of $(A\cup B_1\cup B_3,B_2)$ and 
$(A\cup B_2\cup B_3,B_1)$ is an
$((s-1)(\lam'-1)+1)$-separating minimal%
\REMARK{needs justification using longitudes} 
intrusion in $(G,\uom)$ of order at most
$\lam'$, for some $\lam'<\lam$, and so from the symmetry we may assume
that $(A\cup B_1\cup B_3,B_2)$ has this property.  Since $(M,P_z)$ is not a
cross in $(G,\uom)$ it follows that $M$ and $P_z$ intersect.  Thus
$M\cup P_z$ includes a meridian for $(A\cup B_1\cup B_3,B_2)$.  Finally,
since $Z=B_2\cap (A\cup B_1\cup B_3)\subseteq A\cap B$, the paths 
$(L_v)_{v\in Z}$ form longitudes for $(A\cup B_1\cup B_3,B_2)$, contrary 
to (2).

Thus we have shown that $A\cap B\cap V(\uom)=\{u,v\}$.  Let $Z$ be the set
of all vertices $z\in A$ such that there is no path in $G[A]$ with
one end $z$ and the other end in $X$, let $C=A-Z$ and $D=B\cup Z$.
Then $(C,D)$ is an intrusion with $C\cap D=A\cap B$ and $F$ is a frame
for $(C,D)$ with $V(F)-V(\uom)\subseteq C$.  Since the order of 
$(C,D)$ is at least two, it satisfies the conclusion of the lemma.~\qed

We are ready to deduce the main result of this section.
By a {\em war} in a society $(G,\Omega)$ we mean a set $\cal W$ 
of minimal invasions such that each invasion in $\cal W$ has a meridian,
and $A\cap A'=\emptyset$ for every two distinct
invasions $(A,B),(A',B')\in\cal W$. 
We say that the war $\cal W$ is
{\em $s$-separating} if each invasion in $\cal W$ is $s$-separating,
we say $\cal W$ has {\em order at most $\lambda$} if each member
of $\cal W$ has order at most  $\lambda$, and we say that $\cal W$
is a {\em war of intensity $|{\cal W}|$}.

\begin{lemma}
\label{existinvasions}
\showlabel{existinvasions}
Let $s$, $t$ and $d$ be positive integers, and let 
$b=2^{{2d+1}\choose 2}(2dt)^{2d}(s(2d-1)+2)$. Then 
if a cross-free society  $(G,\Omega)$ of depth at most $d-1$ has
a goose bump of strength $b$, then it has a set $X$ of at most
${2d+1}\choose 2$ vertices such that the society 
$(G\backslash X,\Omega\backslash X)$
has an $s$-separating war of intensity $t$ and order order at most $2d$.
\end{lemma}

\proof
Let $s'=(2d-1)(s-1)+1$.  By Lemma~\ref{existdisjintrus} there exist a set
$X\subseteq V(G)$ with at most $\binom{2d+1}2$ elements and $s'$-separating
intrusions $(A_1,B_1), (A_2,B_2),\dots, (A_{2dt},B_{2dt})$ in $(G\backslash
X,\uom\backslash X)$ of order at most $2d$ such that $A_i\cap A_j=
\emptyset$ for every pair $i,j=1,2,\dots, 2dt$ of distinct integers.
By $2dt$ applications of Lemma~\ref{existtsepinv} there exist, for each
$i=1,2,\dots, 2dt$, and $s$-separating minimal invasion $(C_i,D_i)$ in
$(G\backslash X,\uom\backslash X)$ of order at most $2d$ with a frame
$F_i$ such that $V(F_i)-V(\uom)\subseteq V(A_i)$.  Let $M_i$ be a meridian
for $(C_i,D_i)$, and let $(X_i,Y_i)$ be the base for $(C_i,D_i)$.
Since $(G,\uom)$ has depth at most $d$ there exists a set 
$I\subseteq \{1,2,\dots, 2dt\}$ of size $t$ such that the sets
 $\{X_i\}_{i\in I}$ are pairwise disjoint.  By symmetry we may assume that
$I=\{1,2,\dots, t\}$.  We claim that $(C_1,D_1), (C_2,D_2),\dots ,
(C_t, D_t)$ are as desired.  To prove the claim suppose for a 
contradiction that say $x\in C_i\cap C_j$.  Since $(C_i,D_i)$
is an invasion there exists a path in $G[C_i]$ from $x$ to $X_i\subseteq
Y_j$; therefore this path intersects $C_j\cap D_j$.  Thus there exists a
vertex $v\in C_j\cap D_j\cap C_i$; let $L$ be the longitude of $F_j$
that includes $v$.  But $L$ connects $v\in C_i$ to a vertex of
$X_j\subseteq Y_i\subseteq D_i$, and hence intersects $C_i\cap D_i
\subseteq V(F_i)$.  Thus $F_i$ and $F_j$ intersect.  But
$V(F_i)\cap V(F_j)-V(\uom)\subseteq A_i\cap A_j=\emptyset$ and
$V(F_i) \cap V(F_j)\cap V(\uom)\subseteq X_i\cap X_j=\emptyset$, a 
contradiction.  Thus $(C_1,D_1), (C_2,D_2),\dots, (C_t,D_t)$
satisfy the conclusion of the lemma.~\qed

\section{Using wars}\label{usewars}

\label{sec:wars}
\showlabel{sec:wars}

\begin{lemma}
\label{2nbrs}
\showlabel{2nbrs}
Let $l,t,r$ be positive integers such that $r\ge(t-1){l\choose 2}+1$,
let $(G,\Omega)$ be a connected society, and let $Z\subseteq V(G)$ be a set of
size at most $l$ such that the society $(G\backslash Z,\Omega\backslash Z)$
has a war $\cal W$ of intensity $r$ such that for every $(A,B)\in\cal W$
at least two distinct members of $Z$ have at least one neighbor in $A$.
Then $(G,\Omega)$ has a fan with $t$ blades.
\end{lemma}

\proof There exist distinct vertices $z_1,z_2\in Z$ and a subset 
$\calw'$ of $\calw$ of size $t$ such that for every $(A,B)\in \calw'$ both 
$z_1$ and $z_2$ have a neighbor in $A$.  Furthermore, since $(A,B)$
is a minimal intrusion, it follows that
for every vertex $a\in A$ there exists a path in $G[A]$ from $a$ to
$V(\uom)$.  It follows that $(G,\uom)$ has a fan with $t$ blades,
as desired.~\qed



Let $(A,B)$ be an invasion in a cross-free society $(G,\uom)$,
based at $(X,Y)$, and let $(L_v)_{v\in A\cap B}$ be longitudes for
$(A,B)$.  Let $\uom'$ be a cyclic permutation in $A$ defined as follows:
for each $u\in Y$, if $u$ is an end of $L_v$, then we replace $u$
by $v$, and otherwise we delete $u$.  Then $(G[A],\uom')$ is a society,
and we will call it the {\em society induced by} $(A,B)$.  Since
$(G,\uom)$ is cross-free the definition does not depend on the choice
of longitudes for $(A,B)$.

Assume now that $(G[A],\uom')$ is rural.  A path $P$ in $G[A]$ is called a
{\em perimeter path} in $(G[A], \uom')$ if $A\cap B\subseteq V(P)$ and
$G[A]$ has a drawing in a disk with vertices of $\uom'$ appearing on
the boundary of the disk in the order specified by $\uom'$ and with
every edge of $P$ drawn in the boundary of the disk.

The next lemma is easy and we omit its proof.

\begin{lemma}
\label{indcrossfree}
\showlabel{indcrossfree}
Let $(A,B)$ be an invasion with longitudes $\{P_v\}_{v\in A\cap B}$
in a cross-free society $(G,\Omega)$.
Then the society induced by $(A,B)$ is cross-free.
\end{lemma}

\begin{lemma}
\label{indplanar}
\showlabel{indplanar}
Let  $(G,\Omega)$ be a $5$-connected society, let $Z\subseteq V(G)$ be such that
$(G\backslash Z,\Omega\backslash Z)$ is cross-free, and let
 $(A,B)$ be an invasion in $(G\backslash Z,\Omega\backslash Z)$.
If at most one vertex of $Z$ has a neighbor in $A$, then the
society induced in $(G\backslash Z,\Omega\backslash Z)$ by $(A,B)$
is rural and has a perimeter path.
\end{lemma}

\proof Let $(G[A],\uom')$ be the society induced in $(G\backslash Z,
\uom\backslash Z)$ by $(A,B)$.  By Lemma~\ref{indcrossfree} it is
cross-free and by Theorem~\ref{2paththm} it is rural.  Thus it has a 
drawing in a disk $\udel$ with $V(\uom')$ drawn on the boundary of $\udel$
in the order specified by $\uom'$.  When $\udel$ is regarded as a subset
of the plane, the unbounded face of $G[A]$ is bounded by a walk $W$.
Let $P$ be a subwalk of $W$ containing $A\cap B$.  If $P$ is not a path,
then it has a repeated vertex, say $x$, and $G[A]$ has a separation
$(C,D)$ with $C\cap D=\{x\}$ and $A\cap B\cap V(\uom)\subseteq C$.
Since $(G[A],\uom')$ is cross-free, the latter inclusion implies
that $D-C$ is disjoint from $V(\uom)$ or from $A\cap B$. However,
the latter is impossible, which can be seen by considering the 
drawing of $G[A]$ in $\udel$.  Thus $(D-C)\cap V(\uom)=\emptyset$,
and since $(A,B)$ has longitudes we deduce that $|(D-C)\cap A\cap B|\le 1$.
Let $z\in Z$ be such that no vertex of $Z-\{z\}$ has a neighbor in $A$.
Since $(G,\uom)$ is $4$-connected, the fact that 
$((D-C)\cap A\cap B)\cup\{x,z\}$
does not separate $G$ implies that $D-C$ consists of a unique vertex, 
say $d$, and $d\in A\cap B$.  Furthermore, the only neighbor of $d$ in
$A$ is $x$.  But then $(A-\{d\}, B\cup\{x\})$ contradicts the 
minimality of $(A,B)$.  This proves that $P$ is a path, and it 
follows that it is a perimeter path for $(G[A],\uom')$.~\qed

Let  $(G,\uom)$ be a society. A set $\cal T$ of bumps in  $(G,\uom)$
is called a {\em transaction in  $(G,\uom)$} if there exist
elements $u,v\in V(\Omega)$ such that each member of $\cal T$ has one
end in $u\Omega v$ and the other end in $V(\Omega)-u\Omega v$.
The first part of the next lemma is easy, and the second part
is proved in~\cite[Theorem (8.1)]{RobSeyGM9}.

\begin{lemma}
\mylabel{depthtrans}
Let  $(G,\Omega)$ be a society, and let $d\ge1$ be an integer.
If  $(G,\Omega)$ has depth $d$, then it has no transaction of cardinality
exceeding $2d$. Conversely, if $(G,\Omega)$ has no transaction of 
cardinality exceeding $d$, then it has depth at most $d$.
\end{lemma}

\begin{lemma}
\label{deldepth}
\showlabel{deldepth}
Let  $(G,\Omega)$ be a society of depth $d$, 
and let $X\subseteq V(G)$. Then the society 
$(G\backslash X,\Omega\backslash X)$ has depth at most $2d$.
\end{lemma}

\proof By Lemma~\ref{depthtrans} the society $(G,\uom)$ has
no transaction of cardinality exceeding $2d$.  Then clearly
$(G\backslash X, \uom\backslash X)$ has no transaction of 
cardinality exceeding $2d$, and hence has depth at most $2d$ by
another application of Lemma~\ref{depthtrans}.~\qed


We need one last lemma before we can prove Theorem~\ref{thm1}.
The lemma we need is concerned with the situation when a society of bounded
depth ``almost" has a windmill with $t$ vanes, except that the
paths $P_i$ are not necessarily disjoint and their ends do not
necessarily appear in the right order.
We begin with a special case when the ends of the paths $P_i$
do appear in the right order.

\begin{lemma}
\label{prewindmill}
\showlabel{prewindmill}
Let $t\ge1$ be an integer, and let $\rho=d(t-1)(t'-1)+1$,
where $t'=d(t-1)^2+t$.
Let  $(G,\Omega)$ be a society of depth $d$, let
$(u_1,z_1,v_1,u_2,z_2,v_2,\ldots,u_\rho,z_\rho,v_\rho)$
be clockwise, let $z\in V(G)$, for $i=1,2,\ldots,\rho$
let $P_i$ be a bump with ends $u_i$ and $v_i$, and let
$Q_i$ be a path of length at least one with ends $z$ and $z_i$
disjoint from $V(\Omega)-\{z,z_i\}$.
Assume that the paths $Q_i$ are pairwise disjoint except for $z$,
and that each is disjoint from every $P_j$.
Then  $(G,\Omega)$ has either a windmill with $t$ vanes, 
or a fan with $t$ blades.
\end{lemma}

\proof
By the proof of Lemma~\ref{EPbumps} applied to the paths $P_i$
either some $t$ of those paths are vertex-disjoint, in which case
$(G,\Omega)$ has a windmill with $t$ vanes, or there exists
a set $X\subseteq V(G)$ of size at most $(t-1)d$
such that each $P_i$ uses at least one
vertex of $X$.
We may therefore assume the latter.
For $i=1,2,\ldots,\rho$ the path $P_i$ has a subpath $P_i'$
with one end $u_i$, the other end $x_i\in X$ and no internal vertex in $X$.
Thus there exist $x\in X$ and a set $I\subseteq\{1,2,\ldots,\rho\}$
of size $t'$ such that $x=x_i$ for all $i\in I$.
Let $H$ be the union of all $P_i'$ over $i\in I$.
By an application of Lemma~\ref{EPbumps} to the graph $H\backslash x$
we deduce that either $H\backslash x$ has a goose bump of strength $t$,
in which case  $(G,\Omega)$ has a windmill with $t$ vanes,
or $H$ has a set $Y$ of size at most $(t-1)d$ such that 
$H\backslash Y\backslash x$ has no bumps.
In the latter case for each $i\in I$ there is a path $P''_i$ in $H$ with
one end $u_i$, the other end $y_i\in Y\cup \{x\}$ and otherwise
disjoint from $Y\cup \{x\}$.
Thus there is a vertex $y\in Y\cup \{x\}$ and a set $J\subseteq I$
of size $t$ such that $y_i=y$ for every $i\in J$.
Since $H\backslash Y\backslash x$ has no bumps it follows that
$P''_j$ and $P''_{j'}$ share only $y$ for distinct $j,j'\in J$.
Thus $(G,\Omega)$ has a fan with $t$ blades, as desired.~\qed

Now we are ready to prove the last lemma in full generality.

\begin{lemma}
\label{semiwindmill}
\showlabel{semiwindmill}
Let $t\ge1$ be an integer, and let $\xi=(d+1)\rho$, where $\rho$ is
as in Lemma~\ref{prewindmill}.
Let  $(G,\Omega)$ be a society of depth $d$, 
let $z\in V(G)$, for $i=1,2,\ldots,\xi$ let
$(u_i,z_i,v_i)$ be clockwise, and let
$(u_1,z_1,u_2,z_2,\ldots,u_\xi,z_\xi)$
be clockwise. 
Let $P_i$ be a bump with ends $u_i$ and $v_i$, and let
$Q_i$ be a path of length at least one with ends $z$ and $z_i$
disjoint from $V(\Omega)-\{z,z_i\}$.
Assume that the paths $Q_i$ are pairwise disjoint except for $z$,
and that each is disjoint from every $P_j$.
Then  $(G,\Omega)$ has either a windmill with $t$ vanes,
or a fan with $t$ blades.
\end{lemma}

\proof
Let $(t_1,t_2,\ldots,t_n)$ be
a clockwise enumeration of $V(\Omega)$, and let
$(X_1,X_2,\ldots,X_n)$ be
a corresponding linear decomposition of $(G,\Omega)$ of depth $d$.
Let us fix an integer $i=1,2,\ldots,\rho$, and let
$I=\{(i-1)(d+1)+1,(i-1)(d+1)+2,\ldots,i(d+1)\}$.
For each such $i$ we will construct paths $P^*_i$ and $Q^*_i$
satisfying the hypothesis of Lemma~\ref{prewindmill}.
In the construction we will make use of the paths $P_j$ and $Q_j$
for $j\in I$.

If $(u_j,z_j,v_j,u_{i(d+1)+1})$ is clockwise for some $j\in I$,
then we put $P^*_i=P_j$ and $Q^*_i=Q_j$. 
Otherwise, letting  $s$ be such that $t_s=u_{i(d+1)}$,
we deduce that $P_j$ intersects $X_{t_s}\cap X_{t_{s+1}}$ for all
$j\in I$. Since $|I|>|X_{t_s}\cap X_{t_{s+1}}|$ it follows that
there exist $j<j'\in I$ such that $P_j$ and $P_{j'}$ intersect.
Let $P_i^*$ be a subpath of $P_j\cup P_{j'}$ with ends $u_j$
and $u_{j'}$, and let $Q^*_i=Q_j$.

This completes the construction. The lemma follows from 
Lemma~\ref{prewindmill}.~\qed

\noindent
{\bf Proof of Theorem~\ref{thm1}.}
Let the integers $d$ and $t$ be given, 
let $\xi$ be as in Lemma~\ref{semiwindmill},
let $\ell = 2(t-1)d+
\binom{4d+2}2$, let $\tau =(t-1)\binom\ell 2 +\left(2(t-1)d+
\binom{8d+2}2\right) (6\xi-1)+1$, let $b$ be as in 
Lemma~\ref{existinvasions} with $s=1$, $t=\tau$ and $d$ replaced
by $4d+1$, and let $k$ be as in Lemma~\ref{existgoose} applied to
$b$, $t$, and $4d$.  We will prove that $k$ satisfies the conclusion of
the theorem.

To that end let $(G,\uom)$ be a $k$-cosmopolitan society of 
depth at most $d$, and let $(G_0,\uom_0)$ be a planar truncation of 
$(G,\uom)$.
Let $S\subseteq V(\Omega_0)$. We say that $S$ is {\em sparse}
if whenever $u_1,u_2\in S$ are such that there does not exist $w\in S$
such that $(u_1,w,u_2)$ is clockwise, then there exist two disjoint
bumps $P_1,P_2$ in $(G_0,\Omega_0)$ such that $u_i$ is an end of $P_i$.
The reader should notice that if $H$ is one of the graphs listed
as outcomes (1)-(3) of Theorem~\ref{thm1}, then $V(H)\cap V(\Omega_0)$
is sparse.
We say that  $(G_0,\uom_0)$ is {\em weakly linked} if for every
sparse set $S\subseteq V(\Omega_0)$ there exist $|S|$ disjoint
paths from $S$ to $V(\Omega)$ with no internal vertex in $V(G_0)$.
Thus if the conclusion of the theorem holds for
some weakly linked truncation of $(G_0,\Omega_0)$, then it holds for
$(G,\Omega)$ as well. 
Thus we may assume that $(G_0,\uom_0)$ is a weakly linked 
truncation of $(G,\Omega)$ with $|V(G_0)|$ minimum.
We will prove that $(G_0,\uom_0)$ satisfies the conclusion of
Theorem~\ref{thm1}.
Since $(G_0,\Omega_0)$ is weakly linked, Lemma~\ref{depthtrans} implies
that  $(G_0,\Omega_0)$ has no transaction of cardinality exceeding $2d$,
and hence has depth at most $2d$ by Lemma~\ref{depthtrans}.

By Lemma~\ref{EPcrosses} 
there exists a set $Z_1\subseteq V(G_0)$ such that $|Z_1|\le 2(t-1)d$
and the society $(G_0\backslash Z_1,\uom_1\backslash Z_1)$ is cross-free.
By Lemma~\ref{deldepth} the society $(G_0\backslash Z_1,\uom_0\backslash Z_1)$
has depth at most $4d$.  By Lemma~\ref{existgoose} we may assume that
$(G_0\backslash Z_1,\uom_0\backslash Z_1)$ has a goose bump of strength
$b$.  By Lemma~\ref{existinvasions} there exists a set $Z_2\subseteq
V(G)-Z_1$ such that $|Z_2|\le \binom{4d+2}2$ and in the society 
$(G_0\backslash Z,\uom_0\backslash Z)$ there exists a $1$-separating
war ${\cal W}$ of intensity $\tau$ and order at most $8d+2$, where
$Z=Z_1\cup Z_2$.  If there 
exist at least $(t-1)\binom\ell 2+1$ invasions $(A,B)\in {\cal W}$
such that at least two distinct vertices of $Z$ have a neighbor in $A$,
then the theorem holds by Lemma~\ref{2nbrs}.  We may therefore assume
that this is not the case, and hence $\calw$ has a subset $\calw'$ of
size at least $|Z|(6\xi-1)+1$ such that for every $(A,B)\in\calw' $
at most one vertex of $Z$ has a neighbor in $A$.

Let $(A,B)\in\calw'$ and let $z\in Z$ be such that no vertex in $Z-\{z\}$ has
a neighbor in $A$.
By Lemma~\ref{indplanar} the society $(G_0[A],\uom')$ induced
in $(G_0\backslash Z,\uom_0\backslash Z)$ by $(A,B)$ is rural and has
a perimeter path $P$.  
It follows that $(A\cup\{z\},B\cup\{z\})$ is a separation of $G_0$.
Let $A\cap B=\{w_0,w_1,\ldots,w_s\}$, and let $L_i$ be the longitude
containing $w_i$. Let the ends of $L_i$ be $u_i\in A$ and $v_i\in B$.
We may assume that $(u_0,u_1,\ldots,u_s)$ is clockwise.
The vertices $w_i$ divide $P$ into paths $P_0,P_1,\ldots,P_s$,
where $P_i$ has ends $w_{i-1}$ and $w_i$.
We claim that no $P_i$ includes all neighbors of $z$.
Suppose for a contradiction that $P_i$ does.
Let $(G,\uom)$ be
the composition of $(G_0,\uom_0)$ with a rural neighborhood
$(G_1,\uom,\uom_0)$.  Let $G'_1=G_1\cup G[A\cup\{z\}]$, let 
$G'_0=G_0\backslash (A-B)$ and let $\uom'_0$
consist of $w_s\uom w_0$ followed by $w_{s-1},w_{s-2},\ldots,w_i$
followed by $z$
followed by $w_{i-1},w_{i-2},\ldots,w_1$.
Since $(G[A],\uom')$ is
rural and all neighbors of $z$ belong to $P_i$,
it follows that $(G'_1, \uom,\uom'_0)$ is a rural neighborhood
and $(G,\uom)$ is the composition of $(G'_0,\uom'_0)$ with this
neighborhood.
Thus $(G_0',\Omega_0')$ is a planar truncation of $(G,\Omega)$.
We claim that $(G_0',\Omega_0')$ is weakly linked.
To prove that let $S'\subseteq V(\Omega_0')$ be sparse.
Since $(A,B)$ is a minimal intrusion there exists a set ${\cal P}'$
of $|S'|$ disjoint paths from $S'$ to $V(\Omega_0)$ with no internal
vertex in $G_0'$; let $S$ be the set of their ends in $V(\Omega_0)$.
Since $S'$ is sparse in $(G_0',\Omega_0')$, it follows that $S$ is sparse
in $(G_0,\Omega_0)$. Since $(G_0,\Omega_0)$ is weakly linked there exists
a set $\cal P$ of $|S|$ disjoint paths in $G$ from $S$ to $V(\Omega)$
with no internal vertex in $G_0$.
By taking unions of members of $\cal P$ and ${\cal P}'$ we obtain a set
of paths proving that $(G_0',\Omega_0')$ is weakly linked, as desired.
Since $\calw$ is 1-separating this contradicts the
minimality of $G_0$, proving our claim that
no $P_i$ includes all neighbors of $z$.
The same argument, but with $G_1'=G_1\cup G[A]$ and $\Omega_0'$ not including
$z$ shows that $z$ has a neighbor in $A-B$.

We have shown, in particular, that exactly one vertex of $Z$
has a neighbor in $A-B$.
Thus there exists a subset $\calw''$ of $\calw'$ of size $6\xi$ and
a vertex $z\in Z$ such that for every $(A,B)\in\calw''$ the vertex $z$
has a neighbor in $A-B$.
Now let $w=(A,B)\in \calw''$, and let the notation be as before.
We will construct paths $P_w$, $Q_w$ such that the 
hypotheses of Lemma~\ref{semiwindmill} will be satisfied for at
least half the members $w\in\calw''$.

The facts that $(A,B)$ is a minimal intrusion and that $z$ has a neighbor
in $A-B$ imply that there exists a path $Q_w$ in $G[A\cup\{z\}]$
from $z$ to $z_w\in V(\Omega_0)\cap A$ 
and a choice of longitudes $(L_v:v\in A\cap B)$
for $(A,B)$ such that $Q_w$ is disjoint from all $L_v$.
Referring to the subpaths $P_i$ of the perimeter path $P$ defined above,
since no $P_i$ includes all neighbors of $z$ it follows that there exists
$v\in A\cap B-V(\Omega_0)$. 
We define $P_w$ to be a path obtained from $L_v$ by suitably
modifying $L_v$ inside $B$ such that $P_w$ intersects
$A'$ for at most one $(A',B')\in\calw''-\{(A,B)\}$.
Such modification is easy to make, using the perimeter path of $(A',B')$.
Let $u_w\in A$ and $v_w\in B$ be the ends of $P_w$.

The set ${\cal W}''$ has a subset ${\cal W}'''$ of size $\xi$ such that,
using to the notation of the previous paragraph,
either $(u_w,z_w,v_w)$ is
clockwise for every $w\in {\cal W}'''$ or $(v_w,z_w,u_w)$ is clockwise
for every $w\in {\cal W}'''$, and
for every $w\in {\cal W}'''$ the path $P_w$ is disjoint from
$A'$ for every $(A',B')\in\calw'''-\{w\}$.
The theorem now follows from Lemma~\ref{semiwindmill}.~\qed

\section{Using lack of near-planarity}
\label{sec:lack}
\showlabel{sec:lack}

In this section we prove Theorems~\ref{thm2} and~\ref{thm:society}.
The first
follows immediately from Theorem~\ref{thm1} and the two lemmas below.

\begin{lemma}
\label{planwindmill}
\showlabel{planwindmill}
Let  $(G,\Omega)$ be a rurally $5$-connected society that
is not nearly rural, and let $t$ be a positive integer. If
 $(G,\Omega)$ has a windmill with $4t+1$ vanes, then it has 
a windmill with $t$ vanes and a cross.
\end{lemma}

\proof
Let $x,u_i,v_i,w_i,P_i,Q_i$ be as in the definition of a windmill $W$
with $4t+1$ vanes. Since  $(G\backslash x,\Omega\backslash\{x\})$
is  rurally $4$-connected and not rural, it has a cross $(P,Q)$
by Theorem~\ref{2paththm}. We may choose the windmill $W$ and
cross $(P,Q)$ in $(G\backslash x,\Omega\backslash\{x\})$ such that
$W\cup P\cup Q$ is minimal with respect to inclusion.
If the cross does not intersect the windmill, then the lemma
clearly holds, and so we may assume that a vane $P_i\cup Q_i$ intersects
$P\cup Q$. Let $v$ be a vertex that belongs to both $P_i\cup Q_i$
and $P\cup Q$ such that some subpath $R$ of  $P_i\cup Q_i$ with one end
$v$ and the other end in $V(\Omega)$ has no  vertex in  $(P\cup Q)\backslash v$.
If $R$ has at least one edge, then  $P\cup Q\cup R$ has a proper
subgraph that is a cross, contrary to the minimality of $W\cup P\cup Q$.
Thus $v$ is an end of $P$ or $Q$. Since $P$ and $Q$ have a total of
four ends, it follows that $P\cup Q$ intersects at most four vanes
of $W$. By ignoring those vanes we obtain a windmill with $4(t-1)+1$
vanes,  and a cross $(P,Q)$ disjoint from it. \REMARK{Make clearer.}
The lemma follows.~\qed

\begin{lemma}
\label{planfan}
\showlabel{planfan}
Let  $(G,\Omega)$ be a rurally $6$-connected society that
is not nearly rural, and let $t$ be a positive integer. If
 $(G,\Omega)$ has a fan with $16t+5$ blades, then it has
a fan with $t$ blades and a cross, or 
a fan with $t$ blades and a jump, or
a fan with $t$ blades and two jumps.
\end{lemma}

\proof
Let $z_1,z_2$ be the hubs of a fan $F_2$ with $16t+5$ blades.
If $(G\backslash \{z_1,z_2\},\uom\backslash \{z_1,z_2\})$ has a cross,
then the lemma follows in the same way as Lemma~\ref{planwindmill}, 
and so we may assume  not.  Since $(G\backslash z_1,\uom\backslash 
\{z_1\})$ has a cross, an argument analogous to the proof of
Lemma~\ref{planwindmill} shows that there exists a subfan $F_1$ of
$F_2$ with $4t+1$ blades (that is, $F_1$ is obtained by ignoring a set of
$12t+4$ blades), and two paths $L_2, S_2$ with ends $a_2,c_2$ and
$b_2,z_2$, respectively, such that $x_1,x_2,\dots, x_{4t+1}$, $a_2,b_2,c_2$
is clockwise in $\uom$ for every choice of $x_1,x_2,\dots, x_{4t+1}$ as
in the definition of a fan, and the graphs $L_2,S_2\backslash z_2,
F_1$ are pairwise disjoint.  By using the same argument and the fact that
$(G\backslash z_2,\uom\backslash \{z_2\})$ has a cross we arrive at
a subfan $F$ of $F_1$ with $t$ blades and paths $L_1,S_1$ satisfying
the same properties, but with the index 2 replaced by 1.  We may
assume that $F, L_1,L_2,S_1,S_2$ are chosen so that $F\cup L_1\cup L_2\cup
S_1\cup S_2$ is minimal with respect to inclusion.  This will be referred to
as ``minimality."

If the paths $L_1,L_2,S_1,S_2$ are pairwise disjoint, except possibly
for shared ends and possibly $S_1$ and $S_2$ intersecting, 
then it is easy to see that the lemma holds, and
so we may assume that an internal vertex of $L_1$ belongs to $L_2
\cup S_2$.  Let $v$ be the first vertex on $L_1$ (in either direction)
that belongs to $L_2\cup S_2$, and suppose for a contradiction that
$v$ is not an end of $L_1$. Let $L'_1$ be a subpath of $L_1$ with 
one end $v$, the other end in $V(\uom)$ and no internal vertex in 
$L_2\cup S_2$.  Then by replacing a subpath of $L_2$ or $S_2$
by $L'_1$ we obtain either a contradiction to minimality, or a cross
that is a subgraph of $L_1\cup L_2\cup S_1\cup S_2\backslash \{z_1,
z_2\}$, also a contradiction.  This proves that $v$ is an end of $L_1$,
and hence both ends of $L_1$ are also ends of $L_2$ or $S_2$. 
In particular, $L_1$ and $L_2$ share at least one end.

Suppose first that one end of $L_1$ is an end of $S_2$.
Thus from the symmetry we may assume that $a_1$ is an end of $L_2$ and
$c_1=b_2$; thus $a_2=a_1$, because $a_2,b_2,c_2$ is clockwise.
But now $c_2$ is not an end of $L_1$ or $S_1$, and so the argument of
the previous paragraph implies that no internal vertex of $L_2$
belongs to $S_1\cup L_1$.  The paths $S_1,S_2,L_2$ now show that
$(G,\uom)$ has a fan with $t$ blades and a jump.

We may therefore assume that $a_1=a_2$ and $c_1=c_2$.  
Let $H$ be the union of $L_1, L_2, S_1\backslash z_1$, $S_2\backslash z_2$,
and $V(\Omega)$.  Then the society
$(H, \Omega)$ is rural, as otherwise $(G\backslash \{z_1, z_2\},\Omega)$
 has a cross. 
Let $\Gamma$ be  a drawing of $(H, \Omega)$ in a disk $\Delta$ such that
the vertices of $V(\Omega)$ are drawn on the boundary of $\Delta$ in the
clockwise order specified by $\Omega$.
Let $\Delta'\subseteq\Delta$ be a disk such that $\Delta'$ includes
every path in $\Gamma$ with ends $a_1$ and $c_1$, and the boundary
of $\Delta'$ includes $a_1\Omega c_1$ and a path $P$ of $\Gamma$ 
from $a_1$ to $c_1$.
Then $L_1$ and $L_2$ lie in $\Delta'$, and since
$L_i$ is disjoint from $S_i\backslash z_i$ it follows that
$S_1\backslash z_1$ and $S_2\backslash z_2$ are inside $\Delta'$
and, in particular, are disjoint from $P$. By considering $P$, $S_1$ and $S_2$
we obtain a fan with $t$ blades and a jump.~\qed




\noindent
{\bf Proof of Theorem~\ref{thm2}.}
Let $d$ and $t$ be integers, let $k$ be an integer such that
Theorem~\ref{thm1} holds for $d$ and $16t+5$, and let
$(G,\Omega)$ be a $6$-connected $k$-cosmopolitan society of depth at most $d$.
We may assume that $(G,\Omega)$ is not nearly rural, for otherwise the
theorem holds.
By Theorem~\ref{thm1} the society $(G,\Omega)$ has $t$ disjoint
consecutive crosses, or 
a windmill with $4t+1$ vanes, or
a fan with $16t+5$ blades.
In the first case the theorem holds, and in the second and third case
the theorem follows from Lemma~\ref{planwindmill} and
Lemma~\ref{planfan}, respectively.~\qed

For the proof of  Theorem~\ref{thm:society} we need one more lemma.
Let us recall that presentation of a neighborhood was defined prior to
Theorem~\ref{cosmopolitan}.

\begin{lemma}
\mylabel{nestedtrunk}
Let $d$ and $s$ be integers,
let  $(G,\Omega)$ be an $s$-nested society, and let $(G',\Omega')$ be a planar
truncation of $(G,\Omega)$ of depth at most $d$.
Then $(G,\Omega)$ has an $s$-nested planar truncation of depth at most $2(d+2s)$.
\end{lemma}

\proof 
By a vortical decomposition of a society $(G,\Omega)$ we mean a collection
$(Z_v:v\in V(\Omega))$ of sets such that
\myitem{(i)} $\bigcup (Z_v:v\in V(\Omega))=V(G)$ and every edge of $G$ has
both ends in  $Z_v$ for some $v\in V(\Omega)$,
\myitem{(ii)} for $v\in V(\Omega)$, $v\in Z_v$, and
\myitem{(iii)} if $(v_1,v_2,v_3,v_4)$ is clockwise in $\Omega$, 
then $Z_{v_1}\cap Z_{v_3}\subseteq Z_{v_2}\cup Z_{v_4}$.

\noindent The {\em depth} of such a vortical decomposition is
$\max |Z_u\cap Z_v|$, taken over all pairs of distinct vertices $u,v\in V(\Omega)$
that are consecutive in $\Omega$,
and the depth of $(G,\uom)$ is the minimum depth of a vortical decomposition
of $(G,\uom)$.
Thus if $(G,\Omega)$ has depth at most $d$, then the corresponding linear
decomposition also serves as a vortical decomposition of depth at most~$d$.

Let $(G,\Omega)$ be an $s$-nested society, and let it be the composition
of a society $(G_0,\Omega_0)$ with a rural neighborhood $(G_1,\Omega,\Omega_0)$,
where the neighborhood has a presentation $(\Sigma, \Gamma_1, \Delta, \Delta_0)$
with an $s$-nest $C_1,C_2,\ldots,C_s$.
Let $\Delta_0,\Delta_1,\ldots,\Delta_s$ be as in the definition of $s$-nest.
Let $(G',\Omega')$ be a planar truncation of $(G,\Omega)$ of depth at most $d$.
Then  $(G,\Omega)$ is the composition of  $(G',\Omega')$ with a rural
neighborhood  $(G_2,\Omega,\Omega')$, and we may assume that 
$(G_2,\Omega,\Omega')$ has a presentation $(\Sigma, \Gamma_2, \Delta, \Delta')$,
where $\Delta_0\subseteq\Delta'$.
We may assume that the $s$-nest $C_1,C_2,\ldots,C_s$ is chosen as follows:
first we select $C_1$ such that $\Delta_0\subseteq \Delta_1$ 
and the disk $\Delta_1$ is as small as possible,
subject to that 
we select $C_2$ such that  $\Delta_1\subseteq \Delta_2$ 
and the disk $\Delta_2$ is as small as possible,
subject to that we select $C_3$, and so on.

Let $\Delta^*$ be a closed disk with $\Delta'\subseteq\Delta^*\subseteq\Delta$.
We say that $\Delta^*$ is {\em normal} if whenever an interior point of
an edge $e\in E(\Gamma_1)$ belongs to the boundary of $\Delta^*$, then
$e$ is a subset of the boundary of $\Delta^*$.
A normal disk $\Delta^*$ defines a planar truncation $(G^*,\Omega^*)$
in a natural way as follows: $G^*$ is consists of all vertices and edges
that of $G$ either belong to $G'$, or their image under $\Gamma_1$ belongs to
$\Delta^*$, and $\Omega^*$ consists of vertices of $G$ whose image under
$\Gamma_1$ belongs to the boundary $\Delta^*$ in the order determined by
the boundary of $\Delta^*$.

Given a normal disk $\Delta^*$ and two vertices $u,v\in V(G)$
we define $\xi_{\Delta^*}(u,v)$, or simply $\xi(u,v)$ as follows.
If $u$ is adjacent to $v$, and the image $e$ under $\Gamma_1$ of the edge $uv$
is a subset of the boundary of $\Delta^*$, and for every internal point
$x$ on $e$ there exists an open neighborhood $U$ of $x$ such that
$U\cap \Delta^*=U\cap\Delta_i$, then we let  $\xi(u,v)=i$.
Otherwise we define  $\xi(u,v)=0$.
A short explanation may be in order. If the image $e$ of $uv$ is a subset
of the boundary of $\Delta^*$, then this can happen in two ways:
if we think of $e$ as having two sides, either $\Delta^*$ and $\Delta_i$
appear on the same side, or on opposite sides of $e$. In the definition of
$\xi$ it is only edges with  $\Delta^*$ and $\Delta_i$ on the same side that
count.

We may assume, by shrinking $\Delta'$ slightly, that the boundary of 
$\Delta'$ does not include an interior point of any edge of $\Gamma_2$.
Then $\Delta'$ is normal, and the corresponding planar truncation
is $(G',\Omega')$. Since a linear decomposition
of  $(G',\Omega')$ of depth at most $d$ may be regarded as a vortical 
decomposition of  $(G',\Omega')$ of depth at most $d$, 
we may select a normal disk $\Delta^*$ that gives rise to a planar truncation
$(G^*,\Omega^*)$ of  $(G,\Omega)$, and we may select a vortical decomposition 
$(Z_v:v\in V(\Omega^*))$ of $(G^*,\Omega^*)$ such that
$|Z_u\cap Z_v|\le d+2\xi(u,v)$ for every pair of consecutive vertices of $\Omega^*$.
Furthermore, subject to this, we may choose $\Delta^*$ such that
the number of unordered pairs $u,v$ of distinct vertices of $G$
with $\xi(u,v)=s$ is maximum, subject to that
the number of unordered pairs $u,v$ of distinct vertices of $G$
with $\xi(u,v)=s-1$ is maximum, subject to that
the number of unordered pairs $u,v$ of distinct vertices of $G$
with $\xi(u,v)=s-2$ is maximum, and so on.

We will show that $(G^*,\Omega^*)$ satisfies the conclusion of the theorem.
Let $(t_1,t_2,\ldots,t_n)$ be an arbitrary clockwise enumeration of 
$V(\Omega^*)$, and let $X_i:= Z_{t_i}\cup(Z_{t_1}\cap Z_{t_n})$.
Then $(X_1,X_2,\ldots,X_n)$ is a linear decomposition of $(G^*,\Omega^*)$ of
depth at most $2(d+2s)$.

To complete the proof we must show that $(G^*,\Omega^*)$ is $s$-nested, and
we will do that by showing that each $C_i$ is a subgraph of $G^*$.
To this end we suppose for a contradiction that it is not the case,
and let $i_0\in\{1,2,\ldots,s\}$ be the minimum integer such that
$C_{i_0}$ is not a subgraph of $G^*$.

If $C_{i_0}$ has no edge in $G^*$, then we can construct a new society
$(G_3,\Omega_3)$, where $\Omega_3$ consists of the vertices of $C_{i_0}$ in order,
and obtain a contradiction to the choice of $(G^*,\Omega^*)$. Since the
construction is very similar but slightly easier than the one we are
about to exhibit, we omit the details.
Instead, we assume that $C_{i_0}$ includes edges of both $G^*$ and
$G\backslash E(G^*)$.
Thus there exist vertices 
$x,y\in V(C_{i_0})\cap V(\Omega^*)$ such that some subpath $P$ of $C_{i_0}$ with
ends $x$ and $y$ has no internal vertex in $V(\Omega^*)$.
Let $B$ denote the boundary of $\Delta^*$. 
There are three closed disks with boundaries contained in $B\cup P$.
One of them is $\Delta^*$; let $D$ be the one that is disjoint from $\Delta_0$.
If the interior of $D$ is a subset of $\Delta_{i_0}$ and includes 
no edge of $C_{i_0}$,
then we say that $P$ is a {\em good segment}.
It follows by a standard elementary argument that there is a good segment.

Thus we may assume that $P$ is a good segment, and that the notation is
as in the previous paragraph.
There are two cases: either $D$ is a subset of $\Delta^*$, or the interiors
of $D$ and $\Delta^*$ are disjoint.
Since the former case is handled by a similar, but easier construction,
we leave it to the reader and assume the latter case.
Let $(s_0,s_1,\ldots,s_{t+1})$ be clockwise in $\Omega^*$ such that
$s_0,s_1,\ldots,s_{t+1}$ are all the vertices that belong to $D\cap \Delta^*$.
Thus $\{s_0,s_{t+1}\}=\{x,y\}$.
Let $r_0=s_0,r_1,\ldots,r_k,r_{k+1}=s_{t+1}$ be all the vertices of $P$,
in order, let $H$ be the subgraph of $G^*$ consisting of all vertices
and edges whose images under $\Gamma_1$ belong to $D$, and let
$X:=\{s_0,s_1,\ldots,s_{t+1},r_0,r_1,\ldots,r_{k+1}\}$.
We can regard $H$ as drawn in a disk with the vertices
$s_0,s_1,\ldots,s_{t+1},r_k,r_{k-1},\ldots,r_1$ drawn on the boundary of the
disk in order.
We may assume that every component of $H$ intersects $X$.
The way we chose the cycles $C_{i_0}$ implies that every path in 
$H\backslash\{s_1,s_2,\ldots,s_k\}$ that joins two vertices of $P$ is
a subpath of $P$.
We will refer to this property as the convexity of $H$.
For $i=0,1,\ldots,k+1$ let $b_i$ be the maximum index $j$ such that
the vertex $s_j$ can be reached from $\{r_0,r_1,\ldots,r_i\}$ by a path in $H$
with no internal vertex in $X$.
We define $b_{-1}:=-1$, and let $R_i$ be the set of all vertices of $H$
that can be reached from $\{r_i,s_{b_{i-1}+1},s_{b_{i-1}+2},\ldots,s_{b_i}\}$
by a path with no internal vertex in $X$.
The convexity of $H$ implies that  for $i<j$ the only possible 
member of $R_i\cap R_j$ is $s_{b_i}$.
We now define a new society $(G^{**},\Omega^{**})$ as follows.
The graph $G^{**}$ will be the union of $G^*$ and $H$,
and the cyclic permutation is defined by replacing the subsequence 
$s_0,s_1,\ldots,s_{t+1}$ of $\Omega^*$ by the sequence 
$r_0,r_1,\ldots,r_k,r_{k+1}$.
We define the sets $Z^{**}_v$ as follows. 
For $v\in V(\Omega^*)-V(\Omega^{**})$ we let $Z^{**}_v:=Z_v$.
If $v=r_i$ and $b_i>b_{i-1}$ we define  $Z^{**}_v$ to be the union of
$R_i\cup\{s_{b_i},r_{i-1}\}$ and all $Z_{s_j}$ for 
$j=b_{i-1}+1,b_{i-1}+2,\ldots,b_i$.
If $v=r_i$ and $b_i=b_{i-1}$ we define  
$Z^{**}_v:=R_i\cup\{s_{b_i},r_{i-1}\}\cup (Z_{s_{b_i}}\cap Z_{s_{b_i+1}})$.
It is straightforward to verify that  $(G^{**},\Omega^{**})$ is a planar
truncation of $(G,\Omega)$ and that
$(Z^{**}_v:v\in V(\Omega^{**}))$
is a vortical decomposition of $(G^{**},\Omega^{**})$.
We claim that $\xi_{\Delta^*}(s_j,s_{j+1})<i_0$ for all $j=0,1,\ldots,t$.
To prove this we may assume that $s_j$ is adjacent to $s_{j+1}$, and
let $e$ be the image under $\Gamma_1$ of the edge $s_js_{j+1}$.
It follows that $e$ is a subset of $\Delta_{i_0}$, and hence if 
$s_js_{j+1}\in E(C_k)$ for some $k$, then $k\le i_0$.
Furthermore, if equality holds, then $\Delta_{i_0}$ and $\Delta^*$
lie on opposite sides of $e$, and hence $\xi_{\Delta^*}(s_j,s_{j+1})=0$.
This proves our claim that $\xi_{\Delta^*}(s_j,s_{j+1})<i_0$.
Since for $i=0,1,\ldots,k$ we have 
$Z^{**}_{r_i}\cap Z^{**}_{r_{i+1}}\subseteq (Z_{s_{b_i}}\cap  Z_{s_{b_i+1}})\cup
\{r_i,s_{b_i}\}$,
and $\xi_{\Delta^{**}}(r_i,r_{i+1})=i_0$, 
we deduce that 
$$|Z^{**}_{r_i}\cap Z^{**}_{r_{i+1}}|\le 
|Z_{s_{b_i}}\cap  Z_{s_{b_i+1}}|+2\le
d+\xi_{\Delta^*}(s_{b_i},s_{b_i+1})\le
d+2\xi_{\Delta^{**}}(r_i,r_{i+1}).$$
Thus the existence of $(G^{**},\Omega^{**})$ contradicts
the choice of $(G^*,\Omega^*)$. 
This completes our proof that $C_1,C_2,\ldots,C_s$ are subgraphs of $G^*$,
and hence  $(G^*,\Omega^*)$ is $s$-nested, as desired.~\qed

\noindent
{\bf Proof of Theorem~\ref{thm:society}.}
Let $d$ be as in Theorem~\ref{gm9}, and let $k$ be
as in Corollary~\ref{thm3} applied to $2(d+2s)$ in place of $d$. 
We claim that $k$ satisfies Theorem~\ref{thm:society}.
To prove that let $(G,\Omega)$ be a $6$-connected $s$-nested $k$-cosmopolitan
society that is not nearly rural.
Since $(G,\Omega)$ is an $s$-nested planar truncation of itself,
by Theorem~\ref{gm9} we may assume that $(G,\Omega)$ has
either a leap of length five, in which case it satisfies
Theorem~\ref{thm:society} by Theorem~\ref{thm:leap},
or it has a planar truncation  of depth at most $d$.
In the latter case it has an $s$-nested planar truncation  $(G',\Omega')$ 
of depth at most
$2(d+2s)$ by Lemma~\ref{nestedtrunk}, and
the theorem follows from Corollary~\ref{thm3}
applied to the society $(G',\Omega')$.~\qed

\section{Finding a planar nest}
\label{sec:nest}
\showlabel{sec:nest}



In this section we prove a technical result that applies in the
following situation. We will be able to guarantee
that some societies $(G,\Omega)$ contain certain configurations consisting of
disjoint trees connecting specified vertices in $V(\Omega)$.
The main result of this section, Theorem~\ref{legs} below, states that 
if the society is sufficiently nested, then we can make sure
that the cycles in some reasonably big nest and the trees of the 
configuration intersect nicely.

A {\em target} in a society $(G,\uom)$ is a subgraph  $F$ of $G$ such
that 
\myitem{(i)} $F$ is a forest and every leaf of $F$ belongs to 
$V(\uom)$, and 
\myitem{(ii)} if $u,v\in V(\uom)$ belong to a component $T$ of
$F$, then there exists a component $T'\ne T$ of $F$ and $w\in V(T')
\cap V(\uom)$ such that $(u,w,v)$ is clockwise.

\noindent We say that a vertex $v\in V(G)$ is $F$-{\em special} if 
either $v$ has degree at least three in $F$, or $v$ has degree at
least two in $F$ and $v\in V(\uom)$.

Now let $F$ be a target in $(G,\uom)$ and let $T$ be a component
of $F$.  Let $P$ be a path in $G\backslash V(\Omega)$ with ends $u,v$
such that
$u,v\in V(T)$ and $P$ is otherwise disjoint from $F$. Let $C$ be
the unique cycle in $T\cup P$, and assume that $C$ has at most one
$F$-special vertex.  If $C\backslash u\backslash v$ has no $F$-special
vertex, then let $P'$ be the subpath of $C$ that is complementary to 
$P$, and if $C\backslash u\backslash v$ has an $F$-special vertex, say 
$w$, then let $P'$ be either the subpath of $C\backslash u$ with ends
$v$ and $w$, or the subpath of $C\backslash v$ with ends $u$ and 
$w$.  Finally, let $F'$ be obtained from $F\cup P$ by deleting all
edges and internal vertices of $P'$.  In those circumstances we say that
$F'$ was obtained from $F$ by {\em rerouting}.

A subgraph $F$ of a rural neighborhood $(G,\uom,\uom_0)$ is 
{\em perpendicular} to an $s$-nest $(C_1,C_2,\dots, C_s)$ if for every
component $P$ of $F$
\myitem{(i)} $P$ is a path with one end in $V(\uom)$ and the other
in $V(\uom_0)$, and
\myitem{(ii)} $P\cap C_i$ is a path for all $i=1,2,\dots, s$.

The complexity of a forest $F$ in a society $(G,\uom)$ is
$$\sum (\mbox{deg}_F(v)-2)^++\sum_{v\in V(\uom)} (\mbox{deg}_F(v)-1)^+,$$
where the first summation is over all $v\in V(G)-V(\uom)$ and $x^+$ denotes
$\max (x,0)$.

The following is a preliminary version of the main result of this section.

\begin{theorem}\label{legsbdedtw}
\showlabel{legsbdedtw}
Let $w,s,k$ be positive integers, and let $s'=2w(k+1)+s$.  Then for every
$s'$-nested society $(G,\uom)$ such that $G$ has tree-width
at most $w$ and for every target $F_0$ in $(G,\uom)$ of complexity at most $k$
there exists a target $F$ in $(G,\uom)$ obtained from $F_0$ by repeated
rerouting such that $(G,\uom)$ can be expressed as a composition of
some society with a rural neighborhood $(G',\uom,\uom')$ that
has a presentation with an $s$-nest $(C_1,C_2,\dots, C_s)$ such that 
$G'\cap F$ is perpendicular to $(C_1,C_2,\dots, C_s)$.
\end{theorem}

\proof Suppose that the theorem is false for some integers $w,s,k$, 
a society $(G,\uom)$ and target $F_0$, and choose these entities with
$|V(G)|+|E(G)|$ minimum.
Let $(G,\uom)$ be the composition
of a society $(G_0,\uom_0)$ with a rural neighborhood $(G_1,\uom,\uom_0)$.
Let $\kappa$ be the complexity of $F\cap G_1$ in the society $(G_1,\uom)$, and
let $s''=2w(\kappa+1)+s$.  
Since $(G,\Omega)$ is $s'$-nested and $s''\le s'$ we may choose 
a presentation $(\Sigma, \Gamma,\udel,\udel_0)$ of $(G_1,\uom,\uom_0)$ 
and an $s''$-nest $(C_1,C_2,\dots, C_{s''})$ for it.  
We may assume that $G_0,\uom_0, G_1,F,\Sigma,
\Gamma, \udel,\udel_0,C_1,C_2,\dots, C_{s''}$ are chosen to minimize
$\kappa$.  The minimality of $G$ implies that $G=C_1\cup C_2\cup\cdots
\cup C_{s'}\cup F$.  Likewise, $C_1\cup C_2\cup\cdots\cup C_{s'}$ is
edge-disjoint from $F$, for otherwise contracting an edge belonging to the
intersection of the two graphs contradicts the minimality of $G$.

By a {\em dive} we mean a subpath of $F\cap G_1$ with both ends in 
$V(\uom_0)$ and otherwise disjoint from $V(\uom_0)$.
Let $P$ be a dive with ends $u,v$, and let $P'$ be the corresponding
path in $\ugam$.  Then $\udel_0\cup P'$ separates $\Sigma$;
let $\udel (P')$ denote the component of $\Sigma-\udel_0-P'$
that is contained in $\udel$, and let $H(P)$ denote the subgraph of $G_1$
consisting of all vertices and edges that correspond to vertices or
edges of $\ugam$ that belong to the closure of $\udel (P')$.  Thus
$P$ is a subgraph of $H(P)$.
We say that a dive $P$ is {\em clean} if $H(P)\backslash V(\uom_0)$
includes at most one $F$-special vertex, and if it includes one, say
$v$, then $v\in V(P)$, and no edge of $E(F)-E(P)$ incident with
$v$ belongs to $H(P)$.
The {\em depth} of a dive $P$ is the maximum integer $d\in \{1,2,\dots, s'\}$
such that $V(P)\cap V(C_d)\ne\emptyset$, or 0 if no such integer exists.
It follows from planarity that $|V(P)\cap V(C_i)|\ge 2$ for all
$i=1,2,\dots, d-1$.

\myclaim{1}
Every clean dive has depth at most $2w$.\par
\medskip

To prove the claim suppose for a contradiction that $P_1$ is a clean
dive of depth $d\ge 2w+1$.  Thus $V(P_1)\cap V(C_d)\ne\emptyset$.
Assume that we have already constructed dives $P_1,P_2,\dots, P_t$ for
some $t\le w$ such that $V(P_i)\cap V(C_{d-i+1})\ne\emptyset$ for all
$i=1,2,\dots, t$ and $H(P_t)\subseteq H(P_{t-1})\subseteq\cdots
\subseteq H(P_1)$.  Since $V(P_t)\cap V(C_{d-t+1})\ne\emptyset$, there
exist distinct vertices $x,y\in V(P_t)\cap V(C_{d-t})$.  Furthermore, it is
possible to select $x,y$ such that one of subpaths of $C_{d-t}$
with ends $x,y$, say $Q$, is a subgraph of $H(P_t)$ and no internal vertex of
$Q$ belongs to $P_t$.

\begin{figure}
\begin{center}
\leavevmode
\includegraphics[scale = 1]{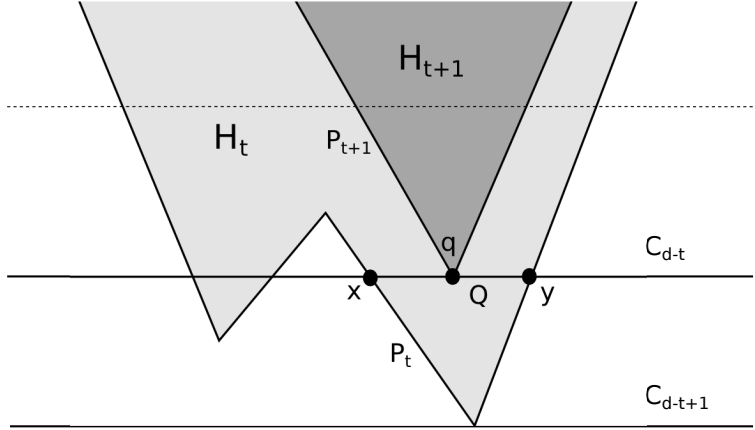}
\end{center}
\caption{Construction of $H(P_{t+1})$.
}
\label{fig:dive}
\end{figure}
\showfiglabel{fig:dive}

We claim that some internal vertex of $Q$ belongs to $F$. Indeed, if not,
then we can reroute $xP_ty$ along $Q$ to produce a target $F'$
and delete an edge of $xP_ty$;
since $P_1$ is clean and $H(P_t)$ is a subgraph of $H(P_1)$ this is
indeed a valid rerouting as defined above.  But this contradicts
the minimality of $G$, and hence some internal vertex of $Q$, say
$q$, belongs to $F$. Since $P_1$ is clean and $H(P_t)$ is a subgraph of
$H(P_1)$ it follows that $q$ belongs to a dive $P_{t+1}$ that is 
a subgraph of $H(P_t)\backslash V(P_t)$.  It follows that
$H(P_{t+1})$ is a subgraph of $H(P_t)$, thus completing the construction. (See Figure~\ref{fig:dive}.)

The dives $P_1,P_2,\dots, P_{w+1}$ just constructed are pairwise
disjoint and all intersect $C_{d-w}$.  Since $d\ge 2w+1$ this implies that
$P_1,P_2,\dots, P_{w+1}$ all intersect each of $C_1,C_2,\dots, C_{w+1}$,
and hence $C_1\cup P_1,C_2\cup P_2,\dots, C_{w+1}\cup P_{w+1}$ is a 
``screen" in $G$ of ``thickness" at least $w+1$.  By 
\cite[Theorem~(1.4)]{SeyThoSearch} the graph $G$ has tree-width
at least $w$, a contradiction.  This proves (1).
\medskip

Our next objective is to prove that $\kappa=0$.  That will take several
steps.  To that end let us define a dive $P$ to be {\em special}
if $P\backslash V(\uom_0)$ contains exactly one $F$-special vertex.
By a {\em bridge} we mean a subgraph $B$ of $G_1\cap F$ consisting of
a component $C$ of $G_1\backslash V(\uom_0)$ together with all edges from
$V(C)$ to $V(\uom_0)$ and all ends of these edges.

\myclaim{2}
If a bridge $B$ includes an $F$-special vertex not in $V(\uom_0)$,
then $B$ includes a special dive.\par
\medskip

To prove Claim (2)
let $B$ be a bridge containing an 
$F$-special vertex not in $V(\uom_0)$.  
For an $F$-special vertex $b\in V(B)-V(\uom_0)$ 
and an edge $e\in E(B)$ incident with
$b$ let $P_e$ be the maximal subpath of $B$ containing $e$ such that  one
end of $P_e$ is $b$ and no internal vertex of $P_e$ is $F$-special
or belongs to $V(\Omega_0)$.
Let $u_e$ be the other end of $P_e$.
The second axiom in the definition of target implies that at most one vertex
of $F$ belongs to $V(\Omega)$.
Since every $F$-special vertex in $V(G_1)-V(\Omega)$ has degree at least three,
it follows that there exists an $F$-special vertex $b\in V(B)-V(\uom_0)$ such that
$u_{e_1},u_{e_2}\in V(\Omega_0)$ for two distinct edges $e_1,e_2\in E(B)$
incident with $b$.
Then $P_{e_1}\cup P_{e_2}$
is as desired.  This proves (2).
\medskip

By (2) we may select a special dive $P$ with
$H(P)$ minimal.  We claim that $P$ is clean.  For let $v\in V(P)-V(\uom_0)$
be $F$-special. If some edge $e\in E(F)-E(P)$ incident with $v$ belongs
to $H(P)$, then there exists a subpath $P'$ of $F$ containing $e$
with one end $v$ and the other end in $V(\uom_0)\cup V(\uom)$.
But $P'$ is a subgraph of $H(P)$, and hence the other end of $P'$ belongs
to $V(\uom_0)$ by planarity.  It follows that $P\cup P'$ includes a
dive 
that contradicts the minimality of $H(P)$.
This proves that the edge $e$ as above does not exist.

It remains to show that no vertex of $H(P)\backslash V(\uom_0)$ except
$v$ is $F$-special.  So suppose for a contradiction that such
vertex, say $v'$, exists.  Then $v'\not\in V(P)$, because $P$ is special,
and hence $v'$ belongs to a bridge $B'\ne B$.  But $B'$ includes a special
dive by (2),
contrary to the choice of $P$.
This proves our claim that $P$ is clean.

By (1)
$P$ has depth at most $2w$.  In particular,
the image under $\ugam$ of some $F$-special vertex belongs to the 
open disk $\udel_{2w+1}$ bounded by the image under $\ugam$ of
$C_{2w+1}$.  Let $G'_0$ consist of $G_0$ and all vertices and edges of 
$G$ whose images under $\ugam$ belong to the closure of $\udel_{2w+1}$,
let $G'_1$ consist of all vertices and edges whose images under 
$\ugam$ belong to the complement of $\udel_{2w+1}$, and let 
$\uom'_0$ be defined by $V(\uom'_0)=V(C_{2w+1})$ and let the cyclic
order of $\uom'_0$ be determined by the order of
$V(C_{2w+1})$. Then $(G,\uom)$ can be regarded as a composition
of $(G'_0,\uom'_0)$ with the rural neighborhood $(G'_1,\uom,\uom'_0)$.
This rural neighborhood has a presentation with a $\sigma$-nest,
where $\sigma=2w\kappa+s$.  On the other hand, the complexity
of $F\cap G'_1$ is at most $\kappa-1$, contrary to the minimality of
$\kappa$.  This proves our claim that $\kappa=0$.

By repeating the argument of the previous paragraph and sacrificing
$2w$ of the cycles $C_i$ we may assume that $(G_1,\Omega,\Omega_0)$
has a presentation with an $s$-nest $C_1,C_2,\ldots,C_s$ and that
there are no dives.
It follows that every component $P$ of $F\cap G_1$ is a path with
one end in $V(\uom)$ and the other in $V(\uom_0)$.  To complete the proof
of the theorem we must show that $P\cap C_i$ is a path for all 
$ i=1,2,\dots, s$.  Suppose for a contradiction that that is not 
the case.  Thus for some $i\in \{1,2,\dots, s\}$ and some component $P$
of $F\cap G_1$ the intersection $P\cap C_i$ is not a path.
Thus there exist distinct vertices $x,y\in V(P\cap C_i)$ such that
$xPy$ is a path with no edge or internal vertex in $C_i$.  Let us choose
$P,i,x,y$ such that, subject to the conditions stated,
$i$ is maximum.  If $i<s$ and $xPy$ intersects $C_{i+1}$, then
$P\cap C_{i+1}$ is not a path, contrary to the choice of $i$.
If $i=1$ or $xPy$ does not intersect $C_{i-1}$, then by rerouting one of the
subpaths of $C_i$ with ends $x,y$ along $xPy$ we obtain
contradiction to the minimality of $G$. Thus we may assume
that $i>1$ and that $xPy$ intersects $C_{i-1}$.

Exactly one of the subpaths of $C_i$ with ends $x,y$, say $Q$,
has the property that the image under $\ugam$ of $xPy\cup Q$ bounds
a disk contained in $\udel$ and disjoint from $\udel_0$.  If no
component of $F\cap G_1$ other than $P$ intersects $Q$, then by
rerouting $F$ along $Q$ we obtain a contradiction to the minimality of
$G$.  Thus there exists a component $P'$ of $F\cap G_1$ other that $P$
that intersects $Q$, say in a vertex $u$.  The vertex $u$ divides $P'$ 
into two subpaths $P_1'$ and $P_2'$.  If both $P_1'$ and $P_2'$ intersect
$C_{i+1}$,  then $P'$ contradicts the choice of $i$. Thus we may
assume that say $P_1'$ does not intersect $C_{i+1}$.  But $P_1'$ 
includes a subpath $P''$ with both ends on $C_i$ and otherwise
disjoint from $C_1\cup C_2\cup\cdots\cup C_s$, and hence by rerouting $C_i$ 
along $P''$ we obtain a contradiction to the minimality of $G$. This
completes the proof of the theorem.~\qed

Before we state the main result of this section we need the
following deep result from \cite{RobSeyGM21}.
A {\em linkage} in a graph $G$ is a subgraph of $G$, every component of
which is a path.  A linkage $L$ in a graph $G$ is {\em vital} if
$V(L)=V(G)$ and there is no linkage $L'\ne L$ in $G$ such that for
every two vertices $u,v\in V(G)$, the vertices $u,v$ are the ends of a
component of $L$ if and only if they are the ends of a component of $L'$.

\begin{theorem}\label{vitallinkage}
\showlabel{vitallinkage}
For every integer $p\ge 0$ there exists an integer $w$ such that
every graph that has a vital linkage with $p$ components has
tree-width at most $w$.
\end{theorem}

Now we are ready to state and prove the main theorem of this
section.  If $F$ is a target in a society $(G,\uom)$ we say that a
vertex $v\in V(G)$ is {\em critical} for $F$ if $v$ is either
$F$-special or a leaf of $F$.  We say that two targets $F,F'$
are {\em hypomorphic} if they have the same set of critical vertices,
say $X$, and $u,v\in X$ are joined by a path in $F$ with no internal
vertices in $X$ if and only if they are so joined in $F'$.

\begin{theorem}\label{legs}
\showlabel{legs}
For every two positive integers $s,k$ there exists an integer
$s'$ such that for every $s'$-nested society $(G,\uom)$
and for every target $F$ in $(G,\uom)$ of complexity at most
$k$ there exists a target $F$ in $(G,\uom)$ obtained from a target
hypomorphic to $F_0$ by repeated rerouting such that $(G,\uom)$ 
can be expressed as a composition of some society with a rural
neighborhood $(G', \uom,\uom')$ that has a presentation 
with an $s$-nest $(C_1,C_2,\dots, C_s)$ such that $G'\cap F$ is
perpendicular to $(C_1,C_2,\dots, C_s)$.
\end{theorem}

\proof We proceed by induction on $|V(G)| + |E(G)|$.  Let $p=k+2$, 
and let $w$ be the bound guaranteed by Theorem~\ref{vitallinkage}.
By hypothesis $(G,\uom)$ is the composition of a society
$(G_0,\uom_0)$ with a rural neighborhood $(G_1,\uom,\uom_0)$,
where $(G_1,\uom,\uom_0)$ has a presentation $(\Sigma,\ugam,\udel,
\udel_0)$ and an $s'$-nest $(C_1,C_2,\dots, C_{s'})$.  Let
$X$ be the set of all vertices critical for $F$, and
let $L=F\backslash X$.  Then $L$ is a linkage in $G\backslash X$.
If it is vital, then $G$ has tree-width at most $|X|+w\le 2k+1+w$,
and hence the theorem follows from Theorem~\ref{legsbdedtw}.

Thus we may assume that $L$ is not vital.  Assume first that there
exists a vertex $v\in V(G)-V(L)$.  If $v\in V(C_i)$ for some
$i\in \{1,2,\dots, s'\}$, then the theorem follows by induction
applied to the graph obtained from $G$ by contracting one of the
edges of $C_i$ incident with $v$; otherwise, the theorem follows
by induction applied to the graph $G\backslash v$.

Thus we may assume that $V(L)=V(G)$, and hence there exists a linkage
$L'\ne L$ linking the same pairs of terminals.  Thus there exists an
edge $e\in E(L)-E(L')$.  If $e\in E(C_i)$ for some $i\in \{1,2,\dots, s'\}$,
then the theorem follows by induction by contracting the edge $e$;
otherwise it follows by induction by deleting $e$, because the
linkage $L'$ guarantees that $G\backslash e$ has a target
hypomorphic to $F$.~\qed

\section{Chasing a turtle}
\label{sec:turtle}
\showlabel{sec:turtle}

In this section we prove Theorem~\ref{thm:largejorgensen},
but first we need
the following two theorems.

\begin{theorem}
\label{thm4}
\showlabel{thm4}
There is an integer $s$ such that if an $s$-nested society 
$(G,\Omega)$ has a turtle, then $G$ has a $K_6$ minor.
\end{theorem}

\begin{figure}
\begin{center}
\leavevmode
\includegraphics[scale = .7]{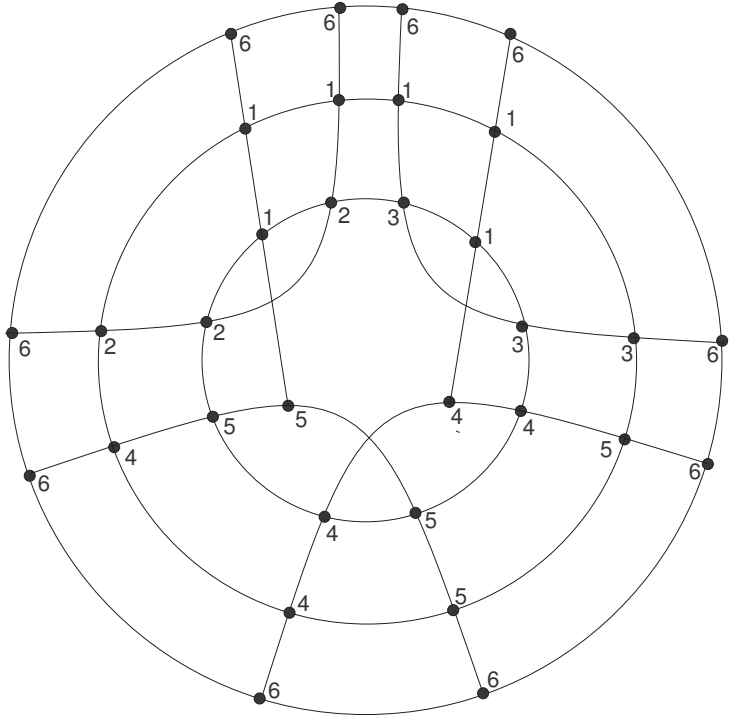}
\end{center}
\caption{
A turtle giving rise to a $K_6$ minor.
}
\label{fig:turtlek6}
\end{figure}

\proof
Let $k$ be the maximum complexity of a turtle, let $s=3$,
and let $s'$ be as in Theorem~\ref{legs}.
We claim that $s'$ satisfies the theorem.
Indeed, let $(G,\Omega)$ be an $s'$-nested society that has a turtle.
Since every turtle is a target, and every target obtained from
a target hypomorphic to a turtle is again a turtle,
we deduce from Theorem~\ref{legs}%
\REMARK{check for problems with rerouting, also for gridlet, etc.}
that $(G,\Omega)$ 
has a turtle $F$ and can be
expressed as a composition of a society with a rural neighborhood
$(G',\Omega,\Omega')$
that has a presentation  with a $3$-nest $(C_1,C_2,C_3)$ such that
$G'\cap F$ is perpendicular to $(C_1,C_2,C_3)$.
It is now fairly straightforward to deduce that $G$ has a $K_6$ minor.
The argument is illustrated in Figure~\ref{fig:turtlek6}.~\qed

\showfiglabel{fig:turtlek6}

\begin{theorem}
\label{thm5}
\showlabel{thm5}
There is an integer $s$ such that if an $s$-nested society
$(G,\uom)$ has three crossed paths, a separated doublecross or a gridlet,
then $G$ has a $K_6$ minor.
\end{theorem}

\begin{figure}
\begin{center}
\leavevmode
\includegraphics[scale = .9]{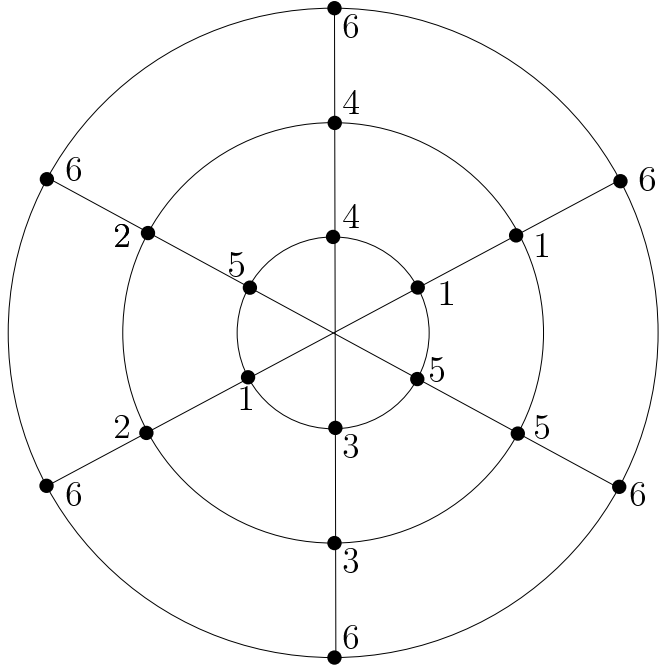}
\end{center}
\caption{
Three crossed paths giving rise to a $K_6$ minor.
}
\label{fig:3xpaths}
\end{figure}

\begin{figure}
\begin{center}
\leavevmode
\includegraphics[scale = .7]{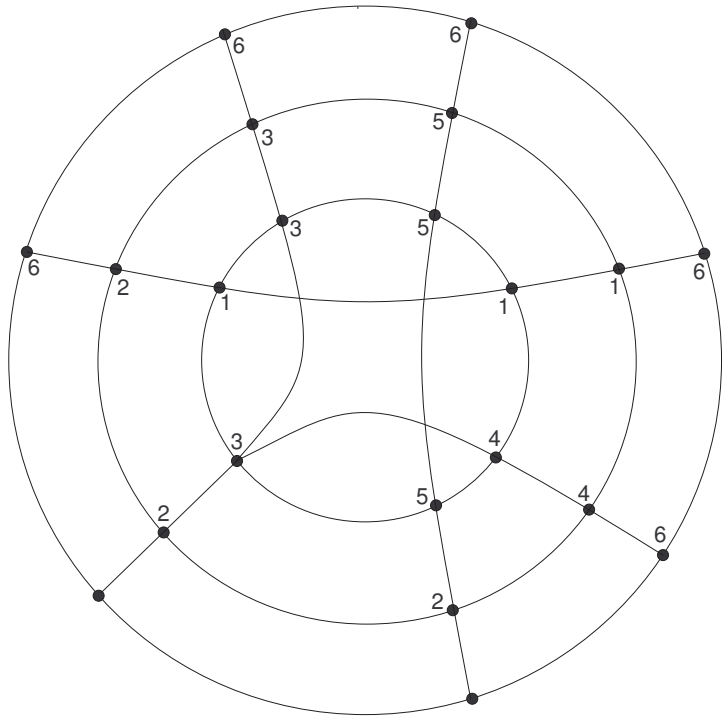}
\end{center}
\caption{
A gridlet giving rise to a $K_6$ minor.
}
\label{fig:gridlet}
\end{figure}

\begin{figure}
\begin{center}
\leavevmode
\includegraphics[scale = .7]{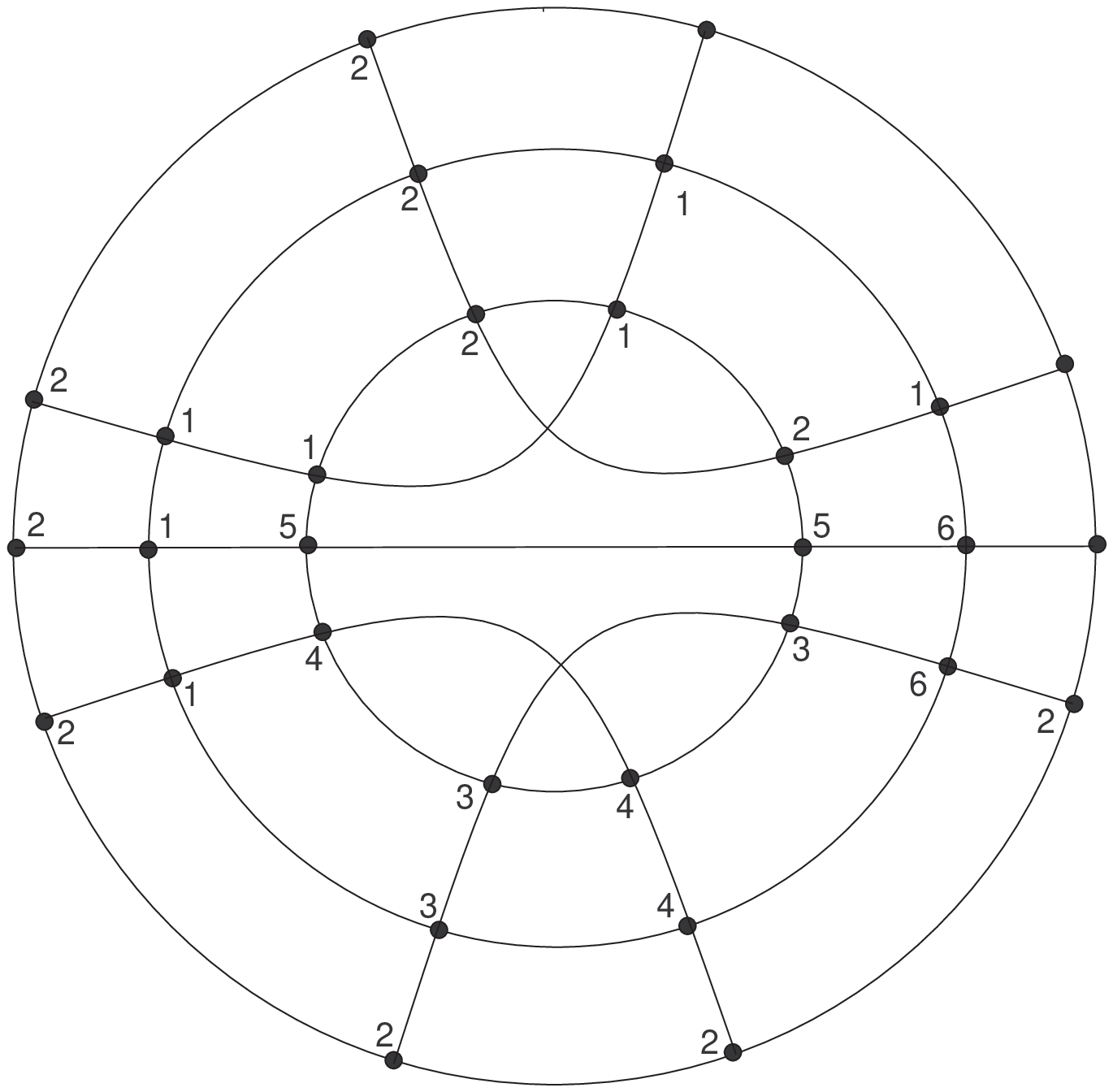}
\end{center}
\caption{
A separated doublecross giving rise to a $K_6$ minor.
}
\label{fig:dblecross}
\end{figure}

\proof
The argument is analogous to the proof of the previous theorem, using
Figures~\ref{fig:3xpaths}, \ref{fig:gridlet}
and~\ref{fig:dblecross} instead. We omit the details.~\qed

\showfiglabel{fig:3xpaths}

\showfiglabel{fig:gridlet}

\showfiglabel{fig:dblecross}

\noindent
{\bf Proof of Theorem~\ref{thm:largejorgensen}.}
Let $s$ be an integer large enough that both Theorem~\ref{thm4} and
Theorem~\ref{thm5} hold for $s$.
Let $k$ be an integer such that Theorem~\ref{thm:society} holds for
this integer.
Let $t$ be such that Theorem~\ref{cosmopolitan} holds for $t$ and
the integer $k$ just defined.
Let $h$ be an integer such that Theorem~\ref{thm:plwall} holds with
$t$ replaced by $t+2s$.
Let $w$ be an integer such that Theorem~\ref{thm:grid} holds for the
integer $h$ just defined.
Finally, let $N$ be as in Theorem~\ref{thm:bdedtwjorgensen}.

Suppose for a contradiction that $G$ is a 6-connected graph on at least
$N$ vertices that is not apex.  By Theorem~\ref{thm:bdedtwjorgensen} 
$G$ has tree-width
exceeding $w$.  By Theorem~\ref{thm:grid} $G$ has a wall of height $h$.
By Theorem~\ref{thm:plwall} $G$ has a planar wall $H_0$ of 
height $t+2s$.  By considering a subwall $H$ of $H_0$
of height $t$ and $s$ cycles of $H_0\backslash V(H)$ we find, by
Theorem~\ref{cosmopolitan}, that the anticompass society $(K,\uom)$ 
of $H$ in $G$ is $s$-nested and $k$-cosmopolitan.  
By Theorem~\ref{thm:society}
 the society $(K,\uom)$ has a turtle,
three crossed paths, a separated doublecross, or a gridlet.  
By Theorems~\ref{thm4} and \ref{thm5} the graph
$G$ has a $K_6$ minor, a contradiction.~\qed

\section*{Acknowledgment}

We would like to acknowledge the contributions of Matthew DeVos and
Rajneesh Hegde, who worked
with us in March 2005 and contributed to this paper,
but did not want to be included as a coauthors.

\baselineskip 11pt
\vfill
\noindent
This material is based upon work supported by the National Science Foundation
under Grants No.~DMS-0200595, DMS-0354742, DMS-0701033, and DMS-0701077.
Any opinions, findings, and conclusions or recommendations expressed in
this material are those of the authors and do not necessarily reflect
the views of the National Science Foundation.

\end{document}